\documentclass[12pt]{amsart}
\usepackage{amssymb,amsthm}
\usepackage[foot]{amsaddr}
\usepackage{amsmath}
\usepackage{tikz}
\usepackage{graphicx}
\usepackage{subcaption}
\usepackage[shortlabels]{enumitem}
\usetikzlibrary{arrows}
\usepackage{float}
\usepackage{graphicx}
\usepackage{subcaption}
\usepackage{hyperref}
\usepackage{pgfplots}

\graphicspath{{MR_figs_final/}}
\usetikzlibrary{decorations.pathmorphing}
\tikzset{snake it/.style={decorate, decoration=snake}}

\makeatletter
\newcommand*{\rom}[1]{\expandafter\@slowromancap\romannumeral #1@}
\makeatother

\numberwithin{equation}{section}

\theoremstyle{plain}
\newtheorem{theorem}{Theorem}
\numberwithin{theorem}{section}

\newtheorem{corollary}[theorem]{Corollary}
\theoremstyle{definition}

\theoremstyle{remark}

\theoremstyle{remark}

\theoremstyle{remark}

\newcommand{\paren}[1]{\left(#1\right)}               

\newcommand{\tblue}[1]{{\color{blue}#1}}

\newcommand{\vb}{\mathbf{b}}

\newcommand{\vu}{{\mathbf{u}}}
\newcommand{\vv}{{\mathbf{v}}}

\newcommand{\bpm}{\begin{pmatrix}}
\newcommand{\epm}{\end{pmatrix}}

\newcommand{\vx}{{\mathbf{x}}}

\newcommand{\vy}{{\mathbf{y}}}

\newcommand{\be}{\begin{equation}}
\newcommand{\ee}{\end{equation}}
\newcommand{\bea}{\begin{eqnarray}}
\newcommand{\eea}{\end{eqnarray}}
\newcommand{\bean}{\begin{eqnarray*}}
\newcommand{\eean}{\end{eqnarray*}}

\newcommand{\bel}[1]{\begin{equation}\label{#1}}

\newcommand{\eel}[1]{{\label{#1}\end{equation}}}

\DeclareMathOperator*{\argmin}{arg\,min}
\DeclareMathOperator{\sinc}{sinc}

\usepackage{fullpage}

\usepackage{datetime}
\usetikzlibrary{quotes, angles}

\providecommand{\keywords}[1]
{
  \small	
  \textbf{\textit{Keywords---}} #1
}

\title[short]{Multi-coil MRI by analytic continuation\\{\footnotesize\ddmmyyyydate\today~\currenttime}}
\author{James W. Webber\textsuperscript{$\dagger$}}
\address[James W. Webber (corresponding author)]{Department of Obstetrics and Gynecology, Brigham and Womens Hospital, 221 Longwood Ave. Boston, MA 02115}
\email[A1]{jwebber5@bwh.harvard.edu\textsuperscript{$\dagger$}}

\begin{document}
\maketitle
\begin{abstract}
We present novel reconstruction and stability analysis methodologies for two-dimensional, multi-coil MRI, based on analytic continuation ideas.
We show that the 2-D, limited-data MRI inverse problem, whereby the missing parts of $\textbf{k}$-space (Fourier space) are lines parallel to either $k_1$ or $k_2$ (i.e., the $\textbf{k}$-space axis), can be reduced to a set of 1-D Fredholm type inverse problems. The Fredholm equations are then solved to recover the 2-D image on 1-D line profiles (``slice-by-slice" imaging). The technique is tested on a range of medical {\it{in vivo}} images (e.g., brain, spine, cardiac), and phantom data. 
  Our method is shown to offer optimal performance, in terms of structural similarity,  when compared against similar methods from the literature, and when the $\textbf{k}$-space data is sub-sampled at random so as to simulate motion corruption. In addition, we present a Singular Value Decomposition (SVD) and stability analysis of the Fredholm operators, and compare the stability properties of different $\textbf{k}$-space sub-sampling schemes (e.g., random vs uniform accelerated sampling).
\end{abstract}
\keywords{{\it{\textbf{Keywords}}} - analytic continuation, multi-coil MRI, Fredholm integral equations, SVD analysis}
\section{Introduction} 
In this paper, we introduce a novel MRI reconstruction and stability analysis methodology, based on the theory of \cite{natterer}. We generalize the theory of \cite{natterer} (applied in that paper to quantitative susceptibility mapping) to limited data, multi-coil MRI, whereby the regions of missing $\textbf{k}$-space are lines parallel either $k_1$ or $k_2$. Without loss of generality, we consider missing lines of $\textbf{k}$-space parallel to $k_1$. The literature considers similar limited data problems in MRI \cite{ESP,ESP1,SENSE,GRAPPA,GRAPPA1,GRAPPA2,TV1,TV2,TV3,feng,motion1,motion2}. The data limitations in those papers are due to, e.g., accelerated imaging \cite{ESP,SENSE} (i.e., deliberate sub-sampling of $\textbf{k}$-space to speed up reconstruction) and movement error \cite{GRAPPA,motion1,motion2} (i.e., when some lines of $\textbf{k}$-space are corrupted due to movement).

In \cite{natterer}, the author introduces a new method to reduce streak artifacts in limited-data MRI reconstruction using analytic continuation, where the missing region of $\textbf{k}$-space is a neighborhood of a cone-shaped surface. Specifically, the author shows that the limited-data MRI problem is equivalent to a 3-D Fredholm equation of the first kind. The Fredholm equation is then solved via repeated approximations and regularized using truncation. The technique is shown to offer significant artifact suppression on simulated phantoms. Similar analytic continuation ideas have been proposed in \cite{AC1,AC2}, for example, to ``fill in" missing X-ray CT and MRI data in \cite{AC1}.

In \cite{ESP}, the authors present a new method, ESPIRiT, for sensitivity map estimation using the SVD of the calibration matrix. The approach is combined with regularized least-squares reconstruction ideas and is shown to offer similar performance to GRAPPA \cite{GRAPPA2} on {\it{in vivo}} (e.g., brain and knee) images. The proposed sensitivity map estimations are shown to be accurate up to absolute value for varying levels of measurement noise, but do not encode the phase sensitivities accurately (the phase is selected at random). Absolute phase estimation has more recently been addressed in \cite{ESP2}, where the authors introduce a new post-processing step to explicitly calculate the coil sensitivities that include the absolute phase of the image. 

In \cite{feng}, the authors introduce a new algorithm for sensitivity encoding and ESPIRiT reconstruction with sparsity constraints. Specifically, the nonlinear, numerically intensive objective function is broken down into two, simpler, sub-problems. The first sub-problem is a linear, least-squares objective, and can be solved efficiently using any appropriate least-squares solver (e.g., Conjugate Gradient Least Squares (CGLS) \cite{CGLS}). The second sub-problem amounts to Total Variation (TV) plus wavelet image denoising, for which there exist a number of highly efficient, high stability methods (e.g., see \cite{Bregman}). The algorithm is shown to outperform other methods from the literature, such as GRAPPA \cite{GRAPPA2} and SENSE \cite{SENSE}, in terms of Root Mean Squared Error (RMSE) on brain image examples.


We propose a new reconstruction methodology for MRI, which is a generalization of the ideas of \cite{natterer} to multi-coil MRI.
Specifically, we show that the (harder to solve) two-dimensional, limited-data MRI inverse problem can be reduced to a set of (easier to solve) one-dimensional Fredholm type inverse problems. The Fredholm equations are then solved to reconstruct the 2-D image on a line-by-line (slice-by-slice) basis. The 1-D formulation has many advantages when compared to the 2-D formulation, including increased stability (i.e., there are less unknowns to recover), efficiency, and more localized reconstruction flexibility. For example, if a radiologist wished to extract a specific line profile of the image to help identify, e.g., a tumor, our method could be applied to extract the line profile efficiently.

To solve the Fredholm equations, the operators are discretized and a least-squares objective is formulated. The inversion is regularized using the smoothed TV penalty of \cite{pet-mri}, and the sensitivity maps are  approximated using ESPIRiT \cite{ESP}. The technique can be thought of as a new ESPIRiT variation, given the least-squares formulation and sensitivity map estimation. In light of this, we compare our method to three other ESPIRiT variations, namely TV ESPIRiT \cite{feng}, $L^1$ ESPIRiT \cite{ESP,ESP1,feng}, and CG (Tikhonov regularized) ESPIRiT \cite{SENSE,ESP}. The proposed technique is shown to offer optimal structural similarity across of range of MRI examples (e.g., brain, spine and phantom images) when compared to the similar ESPIRiT variations listed above, and when the locations of the missing $\textbf{k}$-space lines are selected at random, so as to simulate motion corruption \cite{motion1,motion2}. In the spirit of \cite{natterer}, we denote our method AC, for Analytic Continuation, as the main reconstruction ideas of \cite{natterer} are based on analytic continuation.

In addition to the novel reconstruction methodology, we apply our theorems to stability analysis. Specifically, we calculate the SVD of the 1-D Fredholm operators and present an SVD analysis. For example, in section \ref{coils}, we investigate how the number of coils effects the problem stability using singular value plots, and, in section \ref{subsamp}, we compare the stability properties of the operator for different missing regions of $\textbf{k}$-space (e.g., accelerated vs random sampling). Random undersampling is used to simulate data corruption due to motion, as is, for example, done in \cite{motion1,motion2}. Other studies consider the \emph{geometry factor} (denoted ``$g$-factor" for short) \cite{G1,SENSE} to compare the Signal-to-Noise-Ratio (SNR) of uniform accelerated sampling schemes (e.g., $R=2,3,4\ldots$). The $g$-factor is calculated using the image noise matrix (derived in \cite[appendix A]{SENSE}), which is dependent on a specific Gaussian noise decomposition model. We present an SVD analysis of the operator matrix (i.e., the matrix which is inverted to reconstruct the image). Our analysis is independent of the noise modeling and can be used to gauge the level of noise amplification (e.g., using the singular values) independent of the noise model. We expect our new analysis to work in conjunction with the $g$-factor, rather than in competition, to help inform problem stability.
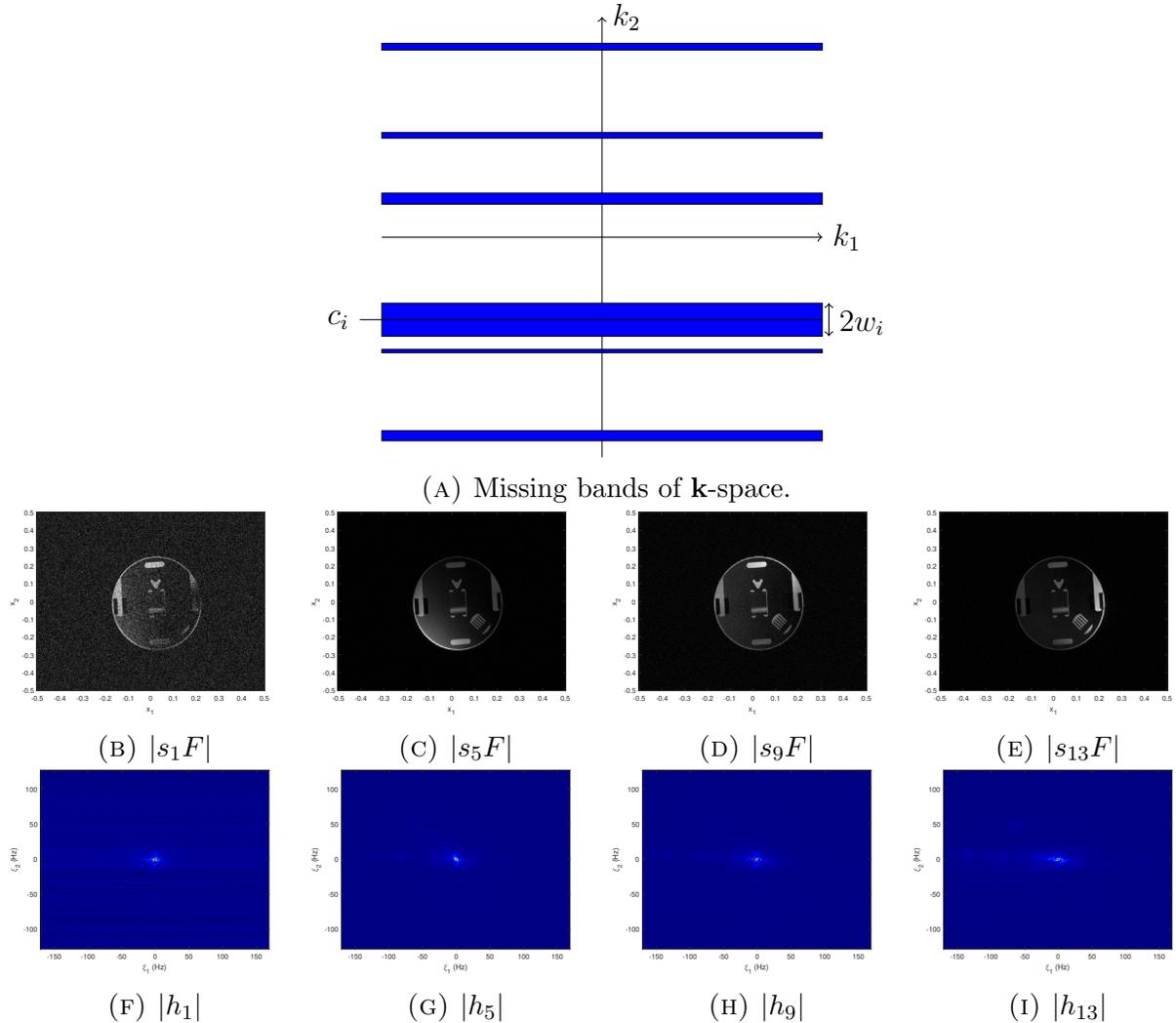
\begin{figure}[!h]
\centering
\begin{subfigure}{1\textwidth}
\centering
\begin{tikzpicture}[scale=3]
\draw [->] (0,-1)--(0,1)node[right]{$k_2$};
\draw [->] (-1,0)--(1,0)node[right]{$k_1$};
\draw [fill=blue] (-1,0.15) rectangle (1,0.2);
\draw [fill=blue] (-1,0.45) rectangle (1,0.475);
\draw [fill=blue] (-1,-0.45) rectangle (1,-0.3);
\draw (-1.1,-0.375)node[left]{$c_i$}--(1,-0.375);
\draw [<->] (1.03,-0.3)--(1.03,-0.45);
\node at (1.18,-0.39) {$2w_i$};
\draw [fill=blue] (-1,-0.525) rectangle (1,-0.51);
\draw [fill=blue] (-1,-0.925) rectangle (1,-0.88);
\draw [fill=blue] (-1,0.85) rectangle (1,0.88);
\end{tikzpicture}
\subcaption{Missing bands of $\textbf{k}$-space.}\label{fig1}
\end{subfigure}
\begin{subfigure}{0.24\textwidth}
\includegraphics[width=0.9\linewidth, height=3.2cm, keepaspectratio]{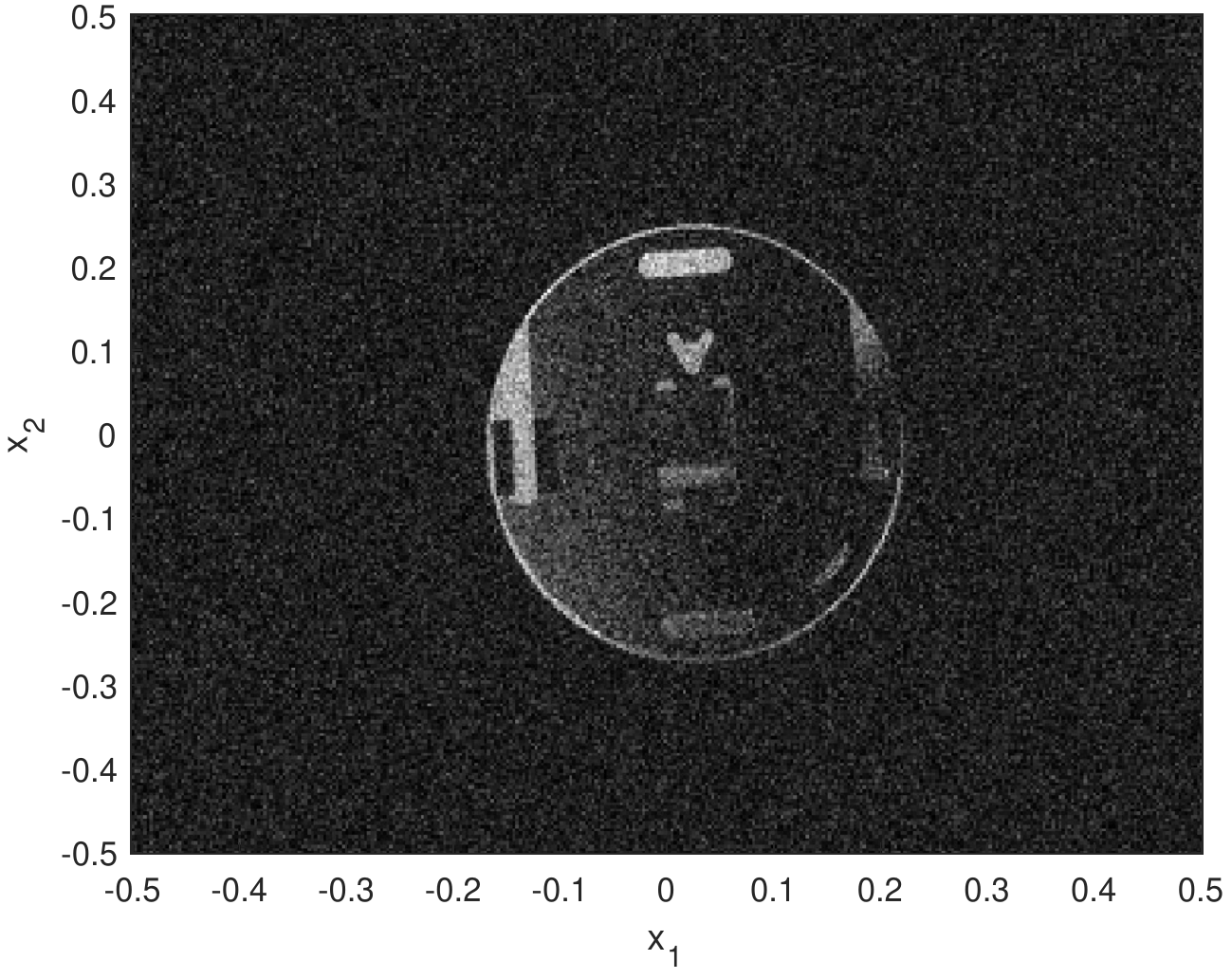}
\subcaption{$|s_1F|$} \label{s1F}
\end{subfigure}
\begin{subfigure}{0.24\textwidth}
\includegraphics[width=0.9\linewidth, height=3.2cm, keepaspectratio]{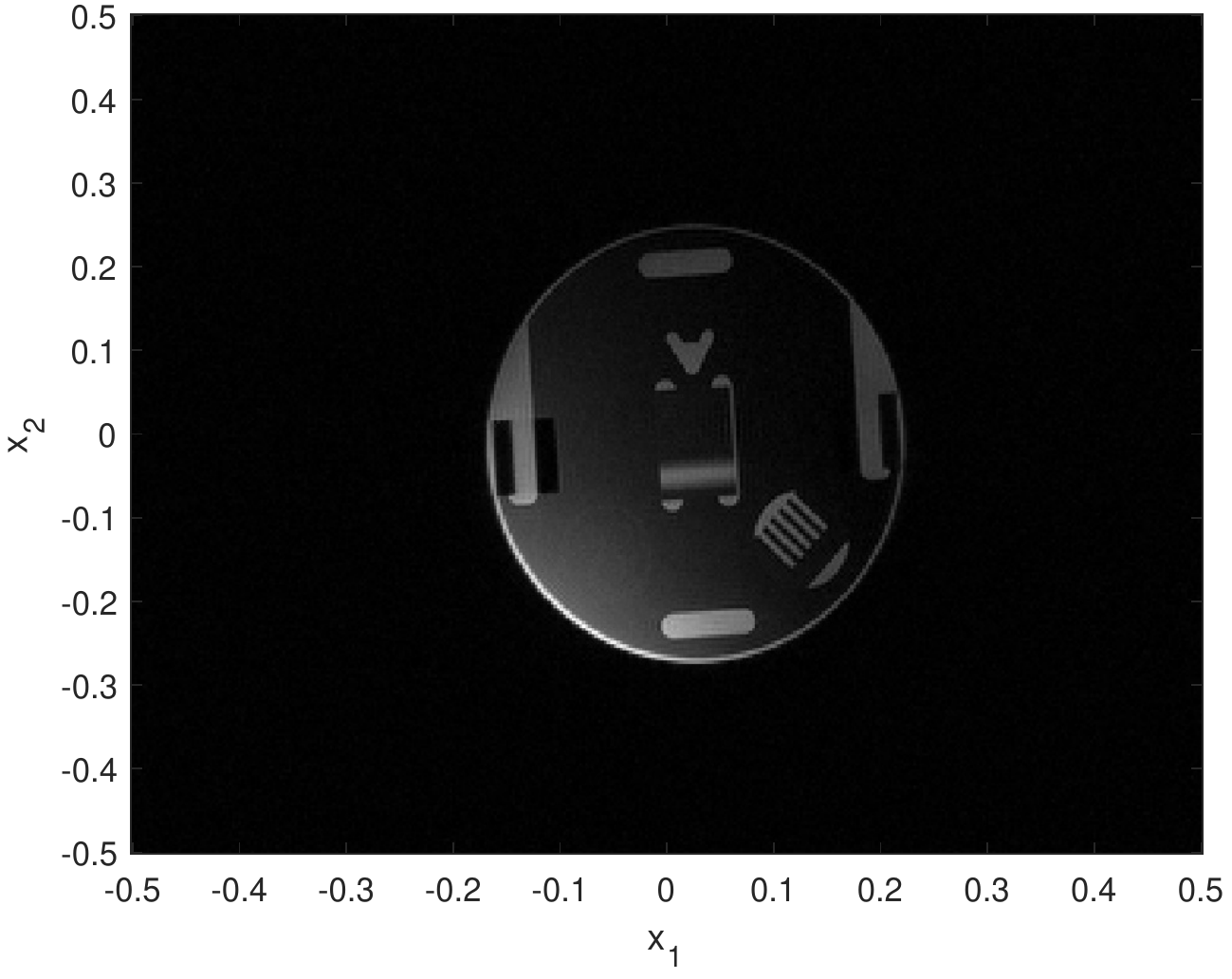}
\subcaption{$|s_5F|$}
\end{subfigure}
\begin{subfigure}{0.24\textwidth}
\includegraphics[width=0.9\linewidth, height=3.2cm, keepaspectratio]{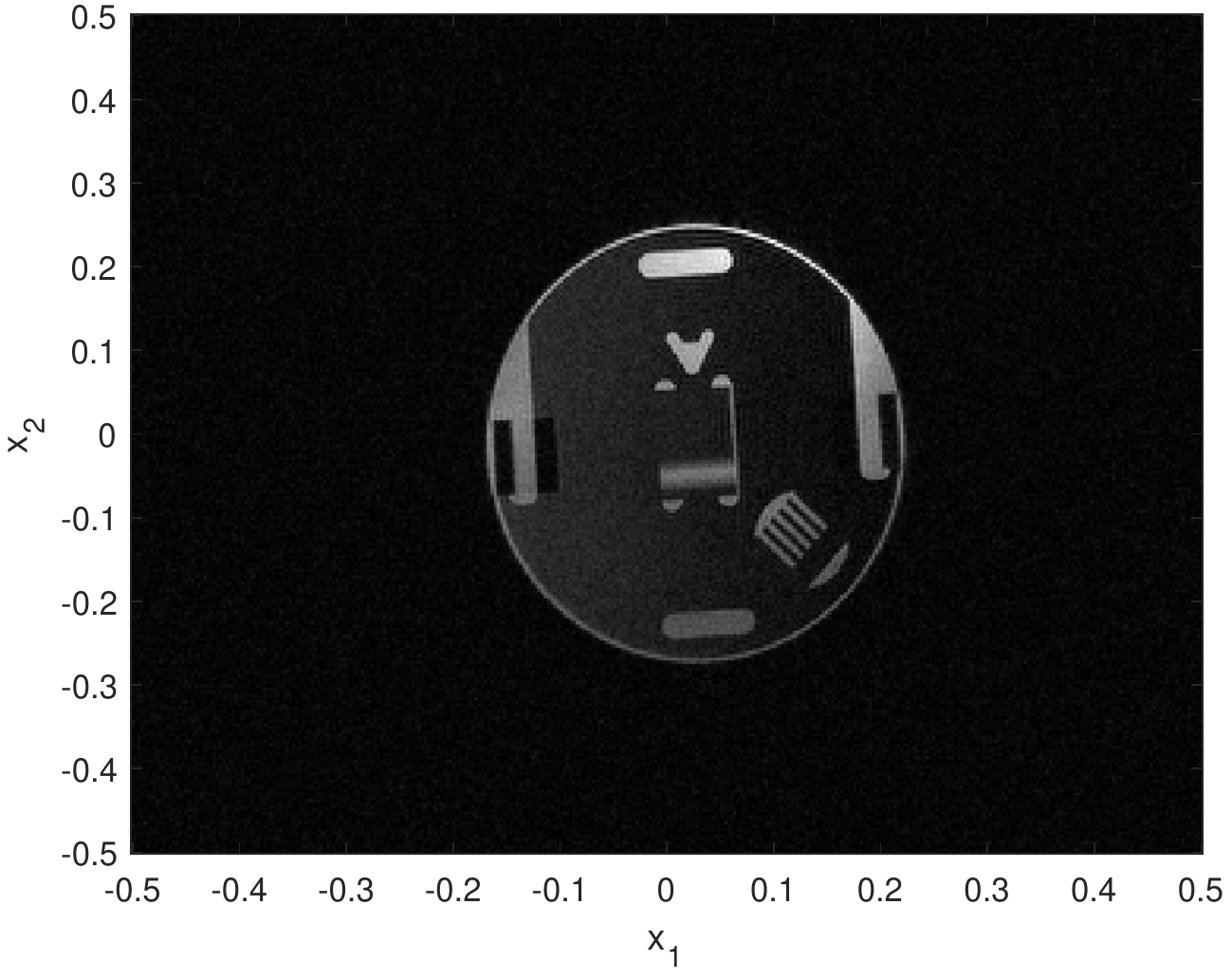}
\subcaption{$|s_9F|$}
\end{subfigure}
\begin{subfigure}{0.24\textwidth}
\includegraphics[width=0.9\linewidth, height=3.2cm, keepaspectratio]{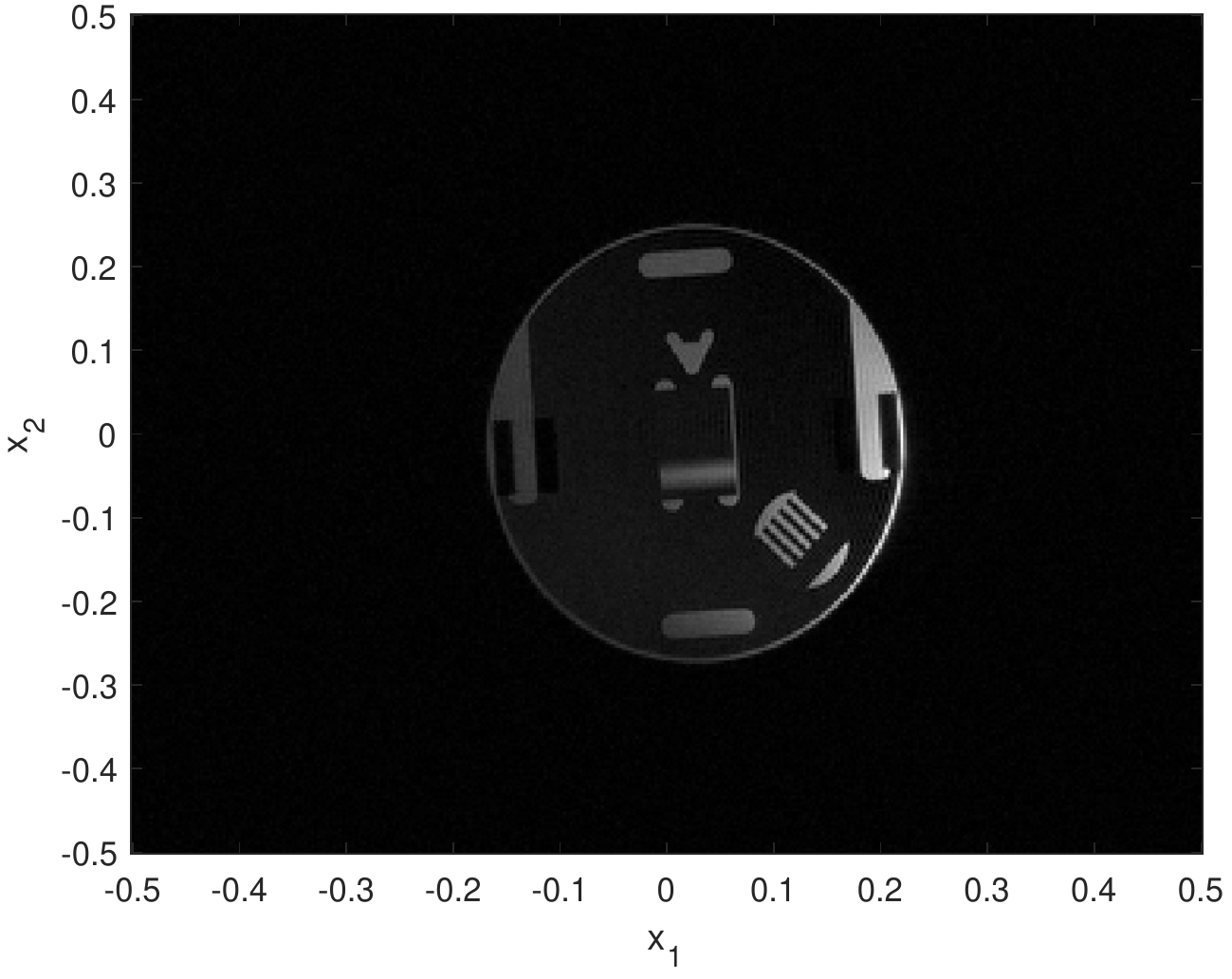}
\subcaption{$|s_{13}F|$}
\end{subfigure}
\begin{subfigure}{0.24\textwidth}
\includegraphics[width=0.9\linewidth, height=3.2cm, keepaspectratio]{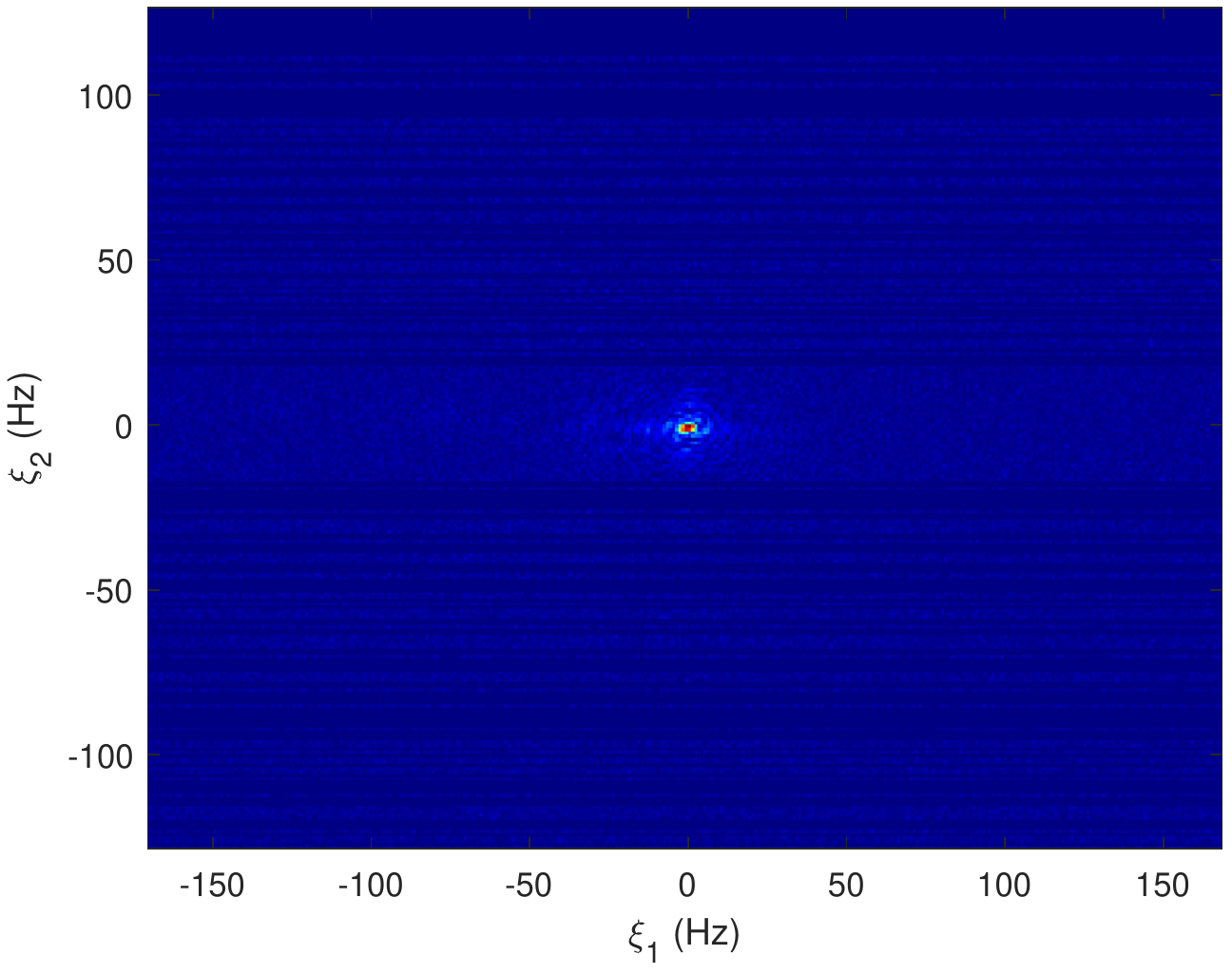}
\subcaption{$|h_1|$}
\end{subfigure}
\begin{subfigure}{0.24\textwidth}
\includegraphics[width=0.9\linewidth, height=3.2cm, keepaspectratio]{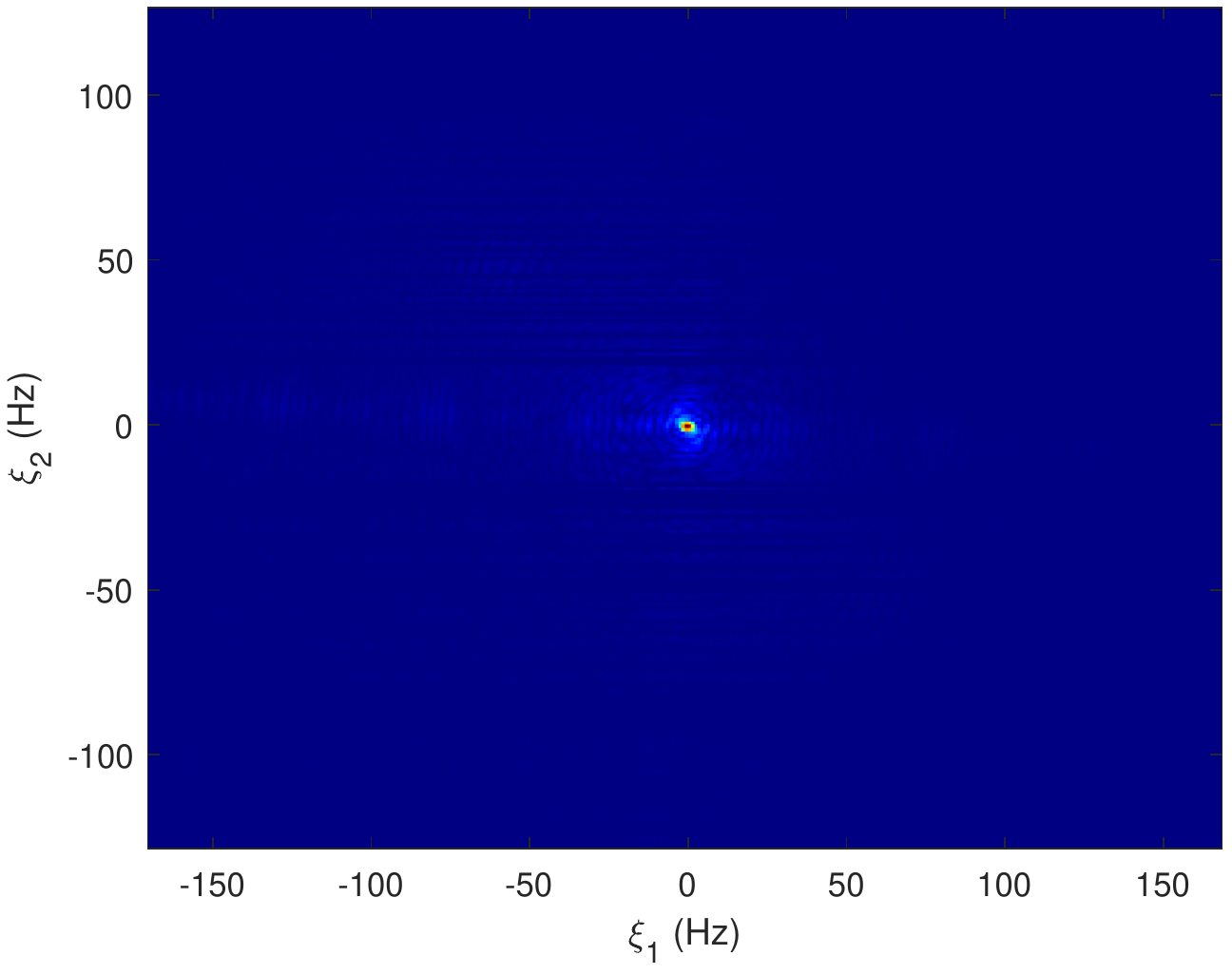}
\subcaption{$|h_5|$}
\end{subfigure}
\begin{subfigure}{0.24\textwidth}
\includegraphics[width=0.9\linewidth, height=3.2cm, keepaspectratio]{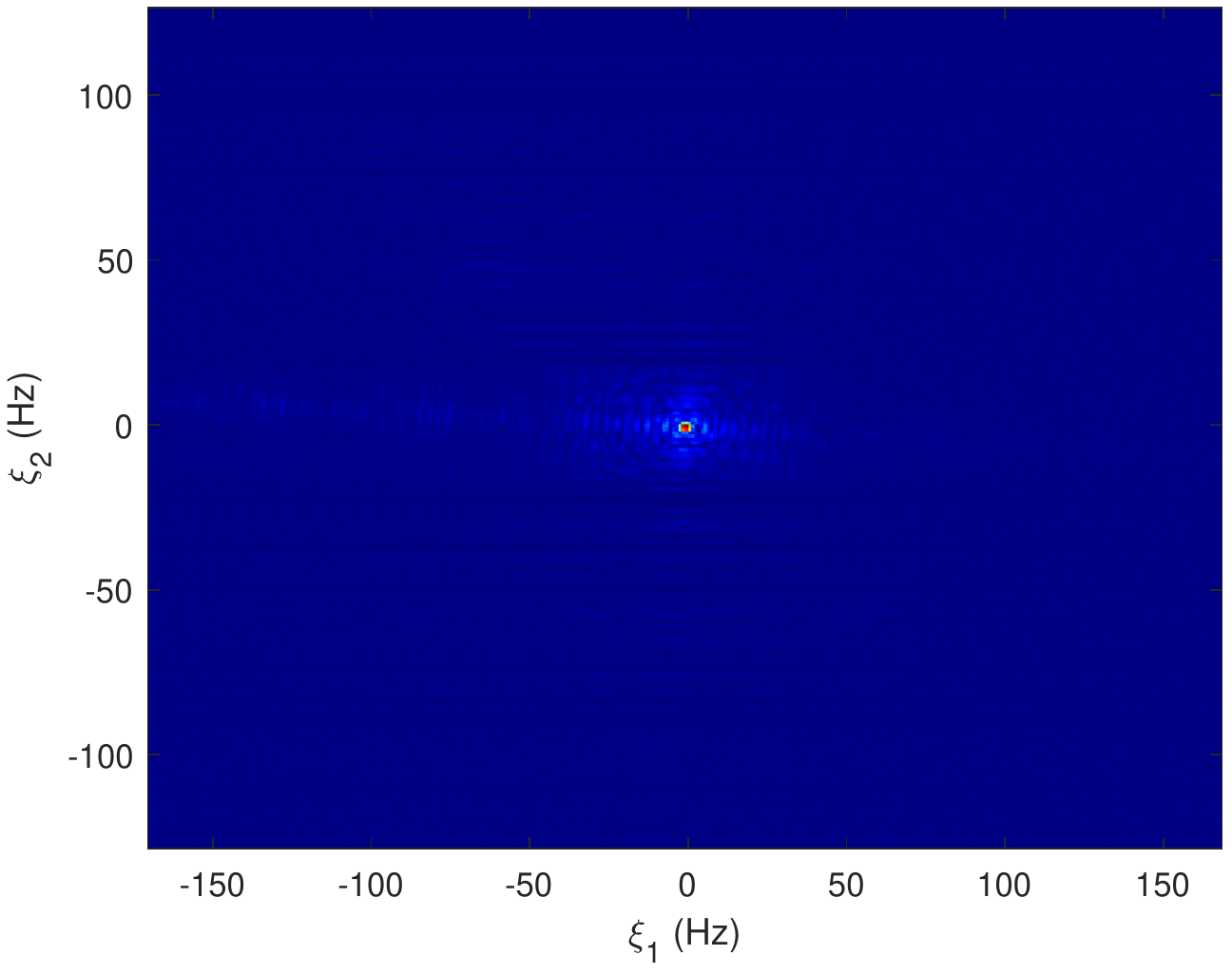}
\subcaption{$|h_9|$}
\end{subfigure}
\begin{subfigure}{0.24\textwidth}
\includegraphics[width=0.9\linewidth, height=3.2cm, keepaspectratio]{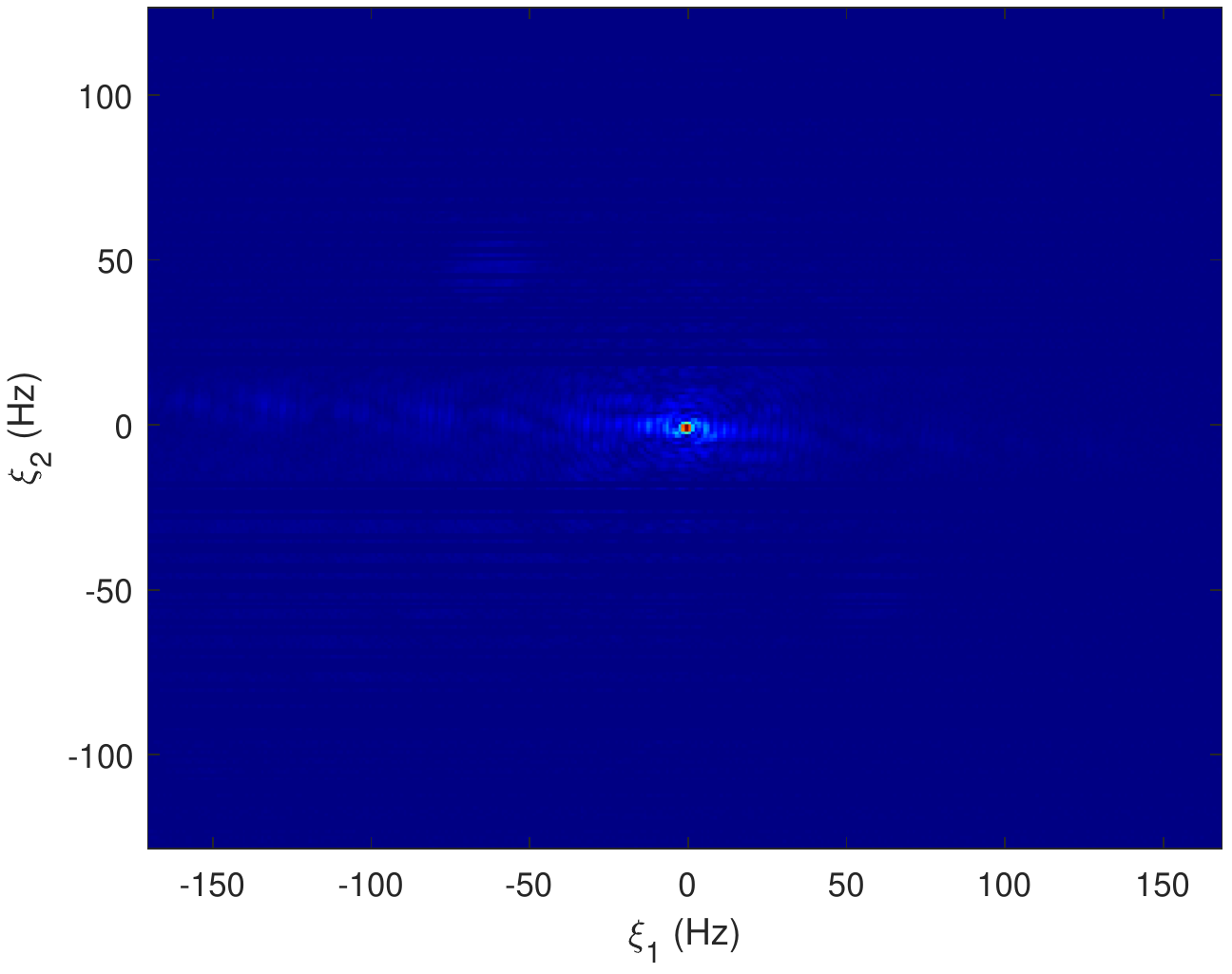}
\subcaption{$|h_{13}|$} \label{h13}
\end{subfigure}
\caption{(A) - Missing bands of $\textbf{k}$-space with center $c_{i}$ and width $2w_i$, displayed as blue rectangles. (B)-(I) - Phantom image data with $K=16$ coils. We show coils $j=1,5,9,13$, going up in steps of 4 from coil 1. Top row - image weighted by sensitivity maps. Bottom row - corresponding $\textbf{k}$-space with missing (zeroed out) lines chosen at random (i.e., $\Lambda$ is chosen at random).}
\label{rp}
\end{figure}

The remainder of this paper is organized as follows. In section \ref{theory}, we present our theory, and show how the weighted, truncated Fourier operator can be decomposed into a set of 1-D Fredholm operators. In section \ref{SVD}, we present an SVD analysis and compare the stability properties of different $\textbf{k}$-space sub-sampling schemes. In section \ref{results}, we introduce our reconstruction method and give a comparison to three similar methods (i.e., CG ESPIRiT, TV ESPIRiT, and $L^1$ ESPIRiT) from the literature.

\section{Theory}
\label{theory}
Let $\Omega\subset\mathbb{R}^2$, let $L^2_c(\Omega)$ denote the set of complex valued square integrable functions with compact support on $\Omega$, and let $C(\Omega)$ be the set of complex valued continuous functions with domain $\Omega$. 

In parallel MRI, the data is modeled by the partial Fourier transform of the image weighted by sensitivity maps
\begin{equation}
\label{ft}
\begin{split}
h_j(\textbf{k})&=\mathcal{F}\paren{s_j F}\paren{\textbf{k}}\\
&=\int_{\left[-\frac{1}{2},\frac{1}{2}\right]^2}s_j(\vx)F(\vx)e^{-i\vx\cdot\textbf{k}}\mathrm{d}\vx,
\end{split}
\end{equation}
where $h_j$ is the signal measured at coil $j\in\{1,2,\ldots,K\}$, $F\in L^2_c([-\frac{1}{2},\frac{1}{2}]^2)$ is the reconstruction target, which is compactly supported on the unit square,  $\vx=(x_1,x_2)$, $\textbf{k}=(k_1,k_2)$, and $s_j\in C([-\frac{1}{2},\frac{1}{2}]^2)$ is the sensitivity map corresponding to coil $j$. Here $K$ is the total number of coils.
The Fourier transform in \eqref{ft} is defined in terms of angular frequency. The data $h_j$ is known for $\textbf{k} \in \mathbb{R}^2\backslash\Lambda$, for each $j$, where $\Lambda=\cup_{i=1}^q\Lambda_i$, a disjoint union, and $\Lambda_i=\mathbb{R}\times [c_i-w_i,c_i+w_i]$, for $1\leq i\leq q$. That is, the missing parts of Fourier space are a set of horizontal, nonintersecting bands with centers $c_i$ and width $2w_i$. See figure \ref{rp}. See figure \ref{fig1} for an illustration of the missing $\textbf{k}$-space blocks, and figures \ref{s1F}-\ref{h13} for some example $s_j F$, $h_j$ pairs on real data. The data is a phantom image from \cite{PULSAR}.

Before we move onto our main theory, we give a short proof of uniqueness of solution in parallel MRI. First we state the Paley-Wiener-Schwartz Theorem \cite[page 22]{hor1}, which will be needed to prove our results.
\begin{theorem}[Paley-Wiener-Schwartz]
Let $\mathcal{E}'(\mathbb{R}^n)$ be the set of distributions of compact support in $\mathbb{R}^n$ and let $f\in\mathcal{E}'(\mathbb{R}^n)$. Then the Fourier transform $\mathcal{F}(f)$ is an entire analytic function.
\end{theorem}
We now have the corollary.
\begin{corollary}
\label{cor1}
Let $F\in L_c^2(\Omega)$, with $\Omega\subset \mathbb{R}^2$, and let $s_1\ldots,s_K\in C(\Omega)$ denote a set of sensitivity maps. Let $\{w_1,\ldots,w_q\}$ and $\{c_1,\ldots,c_q\}$ be a set of widths and centers such that $\Lambda_i \cap \Lambda_j =\emptyset$, for any $i\neq j$, where $\Lambda_i=\mathbb{R}\times [c_i-w_i,c_i+w_i]$. Let $h_j$ be known for $\textbf{k} \in \mathbb{R}^2\backslash \Lambda$, for every $1\leq j\leq K$. Then, the functions $f_j=s_j F$ are uniquely determined, for every $j$.
\begin{proof}
Since $F\in L_c^2(\Omega)$, $f_j\in L_c^2(\Omega)$, for every $j\in\{1,\ldots,K\}$. Thus, $h_j=\mathcal{F}(f_j)$ is an entire analytic function, by the Paley-Wiener-Schwartz Theorem. By definition, $\mathbb{R}^2\backslash \Lambda$ is an open subset of $\mathbb{R}^2$ with nonzero area. Thus, $h_j$ and hence $f_j$ can be determined uniquely on $\mathbb{R}^2$ by analytic continuation, for all $j\in\{1,\ldots,K\}$.
\end{proof}
\end{corollary}

Corollary \ref{cor1} shows that, for any sub-sampling scheme in parallel MRI, the coil images $f_j$ can be recovered uniquely, for all $1\leq j \leq k$. After the coil images are recovered, $F$ is conventionally approximated $F(\vx)\approx\sqrt{\sum_{j=1}^k|f_j|^2}$
by a Sum of Squares (SOS) image. 

We now show that the truncated Fourier operator applied in 2-D parallel MRI operator can be reduced to a set of 1-D Fredholm integral operators, and detail how to reconstruct $F$ on a line-by-line basis.

\subsection{Reduction to a set of 1-D equations}
Let $f_j=s_j F$, as in Corollary \ref{cor1}. Then, by the Fourier inversion formula, we have \cite{natterer}
\begin{equation}
\label{equ1}
\begin{split}
f_j(\vx)&=(2\pi)^{-1}\paren{\int_{\Lambda}\hat{f_j}(\textbf{k})e^{i\vx\cdot\textbf{k}}\mathrm{d}\textbf{k}+\int_{\mathbb{R}^2\backslash\Lambda}\hat{f_j}(\textbf{k})e^{i\vx\cdot\textbf{k}}\mathrm{d}\textbf{k}}\\
&=(2\pi)^{-2}\int_{\Lambda}\underbrace{\left[\int_{\left[-\frac{1}{2},\frac{1}{2}\right]^2}f_j(\vy)e^{-i\vy\cdot\textbf{k}}\mathrm{d}\vy\right]}_{\hat{f_j}(\textbf{k})}e^{i\vx\cdot\textbf{k}}\mathrm{d}\textbf{k}+g_j(\vx)\\
&=\int_{\left[-\frac{1}{2},\frac{1}{2}\right]^2}J(\vx-\vy)f_j(\vy)\mathrm{d}\vy+g_j(\vx),
\end{split}
\end{equation}
where $g_j(\vx)=(2\pi)^{-1}\int_{\mathbb{R}^2\backslash\Lambda}\hat{f_j}(\textbf{k})e^{i\vx\cdot\textbf{k}}\mathrm{d}\textbf{k}$, and
\begin{equation}
J(\vx)=(2\pi)^{-2}\int_{\Lambda}e^{i\vx\cdot\textbf{k}}\mathrm{d}\textbf{k}.
\end{equation}
Here, the Fourier transform $\hat{f_j}$ is defined in terms of angular frequency. 


We can simplify $J$ in the following way
\begin{equation}
\label{K}
\begin{split}
J(\vx)&=(2\pi)^{-2}\int_{\mathbb{R}^2}\sum_{j=1}^m\text{rect}\paren{\frac{k_2-c_j}{2w_j}}e^{i\vx\cdot\textbf{k}}\mathrm{d}\textbf{k}\\
&=(2\pi)^{-2}\int_{\mathbb{R}}e^{ix_1k_1}\mathrm{d}k_1\int_{\mathbb{R}}\sum_{j=1}^m\text{rect}\paren{\frac{k_2-c_j}{2w_j}}e^{ix_2k_2}\mathrm{d}k_2\\
&={\pi}^{-1}\delta(x_1)\sum_{j=1}^mw_j e^{i c_j x_2}\sinc(w_j x_2),
\end{split}
\end{equation}
where
\begin{equation}
\text{rect}(k_2)=\begin{cases}
0 & |k_2|>\frac{1}{2} \\ 
\frac{1}{2} & |k_2|=\frac{1}{2}\\
1 & |k_2|<\frac{1}{2}
\end{cases}
\end{equation}
is the rectangular pulse, and
$$\sinc(x_2)=\frac{\sin(x_2)}{x_2}$$
is the Fourier transform of $\text{rect}$. 
Substituting \eqref{K} into \eqref{equ1} yields
\begin{equation}
\label{equ2}
\begin{split}
g_j(x_1,x_2)&=f_j(x_1,x_2)-{\pi}^{-1}\int_{\tblue{\left[-\frac{1}{2},\frac{1}{2}\right]^2}}\delta(x_1-y_1)\sum_{j=1}^mw_j e^{ic_j (x_2-y_2)}\sinc(w_j (x_2-y_2))f_j(y_1,y_2)\mathrm{d}\vy\\
&=f_j(x_1,x_2)-\int_{\left[-\frac{1}{2},\frac{1}{2}\right]}\underbrace{\left[{\pi}^{-1}\sum_{j=1}^mw_j e^{ic_j (x_2-y_2)}\sinc(w_j (x_2-y_2))\right]}_{L(x_2-y_2)}f_j(x_1,y_2)\mathrm{d}y_2\\
&=f_j(x_1,x_2)-\int_{\left[-\frac{1}{2},\frac{1}{2}\right]}L(x_2-y_2)f_j(x_1,y_2)\mathrm{d}y_2\\
&=(I-L)(f_j)(x_1,x_2)
\end{split}
\end{equation}
where
$$L(x_2)={\pi}^{-1}\sum_{j=1}^mw_j e^{ic_j x_2}\sinc(w_j x_2),$$
$g_j(x_1,x_2)=g_j(\vx)$, and $f_j(x_1,x_2)=f_j(\vx)$. Thus, for each fixed $x_1$, we are left with solving a set of simultaneous 1-D Fredholm integral equations to recover $F(x_1,\cdot)$ (i.e., the vertical line profile of $F$ at $x_1$). The operator $I-L$ of \eqref{equ2} is the sum of two parts, namely the identity map ($I$), which maps $f_j \to f_j$, and $L$, which maps $f_j$ to its convolution with a sinc type kernel. That is, $g_j=(I-L)f_j$ is the original image, $f_j$, plus artifacts induced by $L$.

\subsection{Discretization}
Here we detail the discrete form of the Fredholm operators in equation \eqref{equ2}.  Let $\left[-\frac{1}{2},\frac{1}{2}\right]^2$ (the reconstruction space) be discretized to the $n\times m$ grid
$$\mathcal{U}\times \mathcal{V} = \{u_1,\ldots,u_m\} \times \{v_1,\ldots,v_n\},$$
where $u_i=-\frac{1}{2}+\frac{i}{m}$, and $v_i=-\frac{1}{2}+\frac{i}{n}$. Let $A\in\mathbb{C}^{n\times n}$ be the discretized form of the operator $(I-L)$, introduced in the last line of equation \eqref{equ2}. Let 
$$S^{(j)}_{i}= \text{diag}\paren{s_j(u_i,v_1),\ldots, s_j(u_i,v_n)}.$$
Then, the Fredholm equations of \eqref{equ2} have the discrete formulation
\begin{equation}
\label{1D}
\begin{bmatrix}
AS^{(1)}_{i} \\
AS^{(2)}_{i} \\
\vdots \\
AS^{(k)}_{i} 
\end{bmatrix}
\begin{pmatrix}
F\paren{u_i,v_1} \\
F\paren{u_i,v_2} \\
\vdots \\
F\paren{u_i,v_n}
\end{pmatrix}=
\begin{pmatrix}
\vb^{(1)}_{i} \\
\vb^{(2)}_{i} \\
\vdots \\
\vb^{(k)}_{i}
\end{pmatrix},
\end{equation}
for every $1\leq i\leq m$, where
$$\vb^{(j)}_{i}=\paren{g_j\paren{u_i,v_1},\ldots,g_j\paren{u_i,v_n}}^T.$$
Thus, solving for $F$ on $\mathcal{U}\times \mathcal{V}$ is equivalent to inversion of the block-diagonal matrix
\begin{equation}
\label{M}
M=\begin{bmatrix}
A_1 & 0 & \ldots & 0 \\
0 & A_2 & \ldots & 0 \\
\vdots & \vdots & \ddots & \vdots \\
0 & 0 & \ldots & A_m
\end{bmatrix}\in\mathbb{C}^{nkm\times nm}, \ \text{where} \ 
A_i=\begin{bmatrix}
AS^{(1)}_{u_i} \\
AS^{(2)}_{u_i} \\
\vdots \\
AS^{(k)}_{u_i} 
\end{bmatrix}\in\mathbb{C}^{nk\times n}.
\end{equation}
For the purposes of reconstruction, we do not build $M$. In practice, to reconstruct $F$, we solve \eqref{1D} for every $u_i$, i.e., we recover $F$ slice-by-slice on lines parallel to the $x_2$ axis. We discuss in more detail the reconstruction method in section \ref{results}. The above formulation for $M$ is needed for the next section, and for stability analysis purposes.

\subsection{Singular Value Decomposition analysis}
\label{SVD}
In this sub-section, we present an SVD and stability analysis of $M$, as defined in equation \eqref{M}. Let
$$A_i=U_i\Sigma_iV_i$$
be decomposed into its SVD. Then, $M$ can be decomposed as $M=U\Sigma V$, where
\begin{equation}
\label{Msvd}
U=\begin{bmatrix}
U_1 & 0 & \ldots & 0 \\
0 & U_2 & \ldots & 0 \\
\vdots & \vdots & \ddots & \vdots \\
0 & 0 & \ldots & U_m
\end{bmatrix}, \ \ 
\Sigma=\begin{bmatrix}
\Sigma_1 & 0 & \ldots & 0 \\
0 & \Sigma_2 & \ldots & 0 \\
\vdots & \vdots & \ddots & \vdots \\
0 & 0 & \ldots & \Sigma_m
\end{bmatrix}, \ \
V=\begin{bmatrix}
V_1 & 0 & \ldots & 0 \\
0 & V_2 & \ldots & 0 \\
\vdots & \vdots & \ddots & \vdots \\
0 & 0 & \ldots & V_m
\end{bmatrix}.
\end{equation}
Thus, given the block diagonal form of $M$, we can calculate the SVD of $M$ using the SVD of the $A_i$. The columns of $U$ and $V$ are the left and right singular vectors of $M$, respectively, and the diagonal entries of $\Sigma$ are the singular values of $M$. Through analysis of $U$, $\Sigma$, and $V$, we can gain insight into the inversion stability of $M$, and the magnitude to which noise in the data will be amplified in the reconstruction. 

Let $\sigma=(\sigma_1,\ldots,\sigma_{nm})$ be the vector of singular values of $M$, ordered such that $\sigma_i\geq\sigma_{i+1}$ for all $1\leq i\leq nm-1$.  Then, the pseudoinverse of $M$, $M^{\dagger}$, is defined
\begin{equation}
M^{\dagger}\vb=\sum_{i=1}^{nm}\frac{\vu_i^T\vb}{\sigma_i}\vv_i,
\end{equation}
for some data $\vb$, where $\vu_i$ and $\vv_i$ denotes column $i$ of $U$ and $V$, respectively. Since we divide by the $\sigma_i$ when inverting $M$, we must pay close attention to the size of the $\sigma_i$, since the smaller $\sigma_i$ will more greatly amplify any noise in $\vb$. In particular, we consider the following metrics to assess stability:
\begin{itemize}
\item Condition number - defined as the ratio of the maximum and minimum singular value
$$\kappa(M)=\frac{\sigma_1}{\sigma_{nm}}\geq 1.$$
$\kappa(M)$ close to 1 implies less noise amplification, and vice-versa.
\item Effective null space dimension - for some threshold $t\approx 0$, we define
$$d(M)=nm-\argmin_{1\leq i\leq nm}|\sigma_i-t|.$$ 
$d$ thus measures the number of singular values close to zero. Large $d(M)$ indicates high-dimension null space and an unstable inversion, and conversely for small $d(M)$.
\end{itemize}
We also consider the right singular vectors $\vv_i$, to analyze the noise amplification in different parts of the image.
\begin{figure}
\centering
\begin{subfigure}{0.24\textwidth}
\includegraphics[width=0.9\linewidth, height=3.2cm, keepaspectratio]{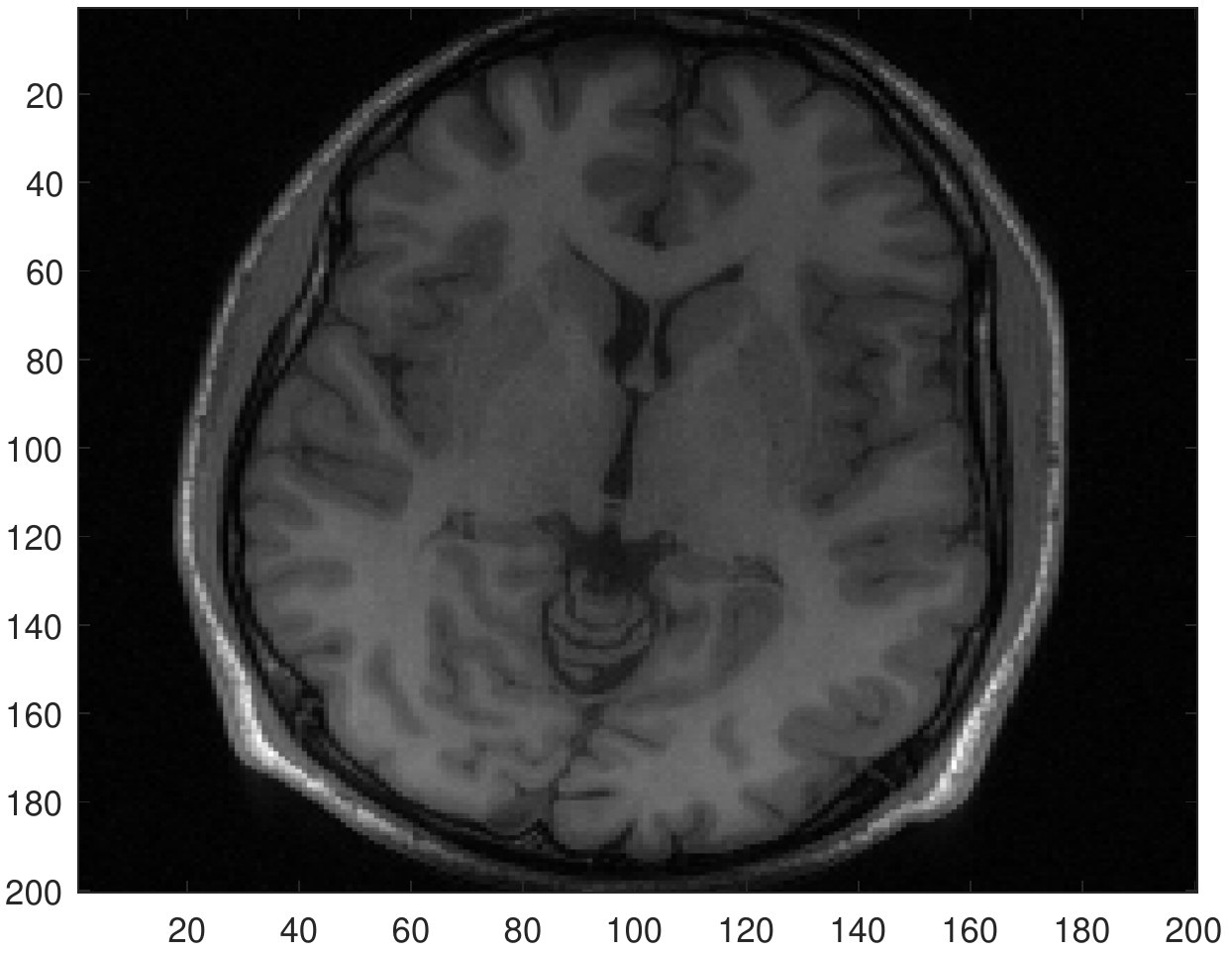}
\subcaption{$|F|$}
\end{subfigure}
\begin{subfigure}{0.24\textwidth}
\includegraphics[width=0.9\linewidth, height=3.2cm, keepaspectratio]{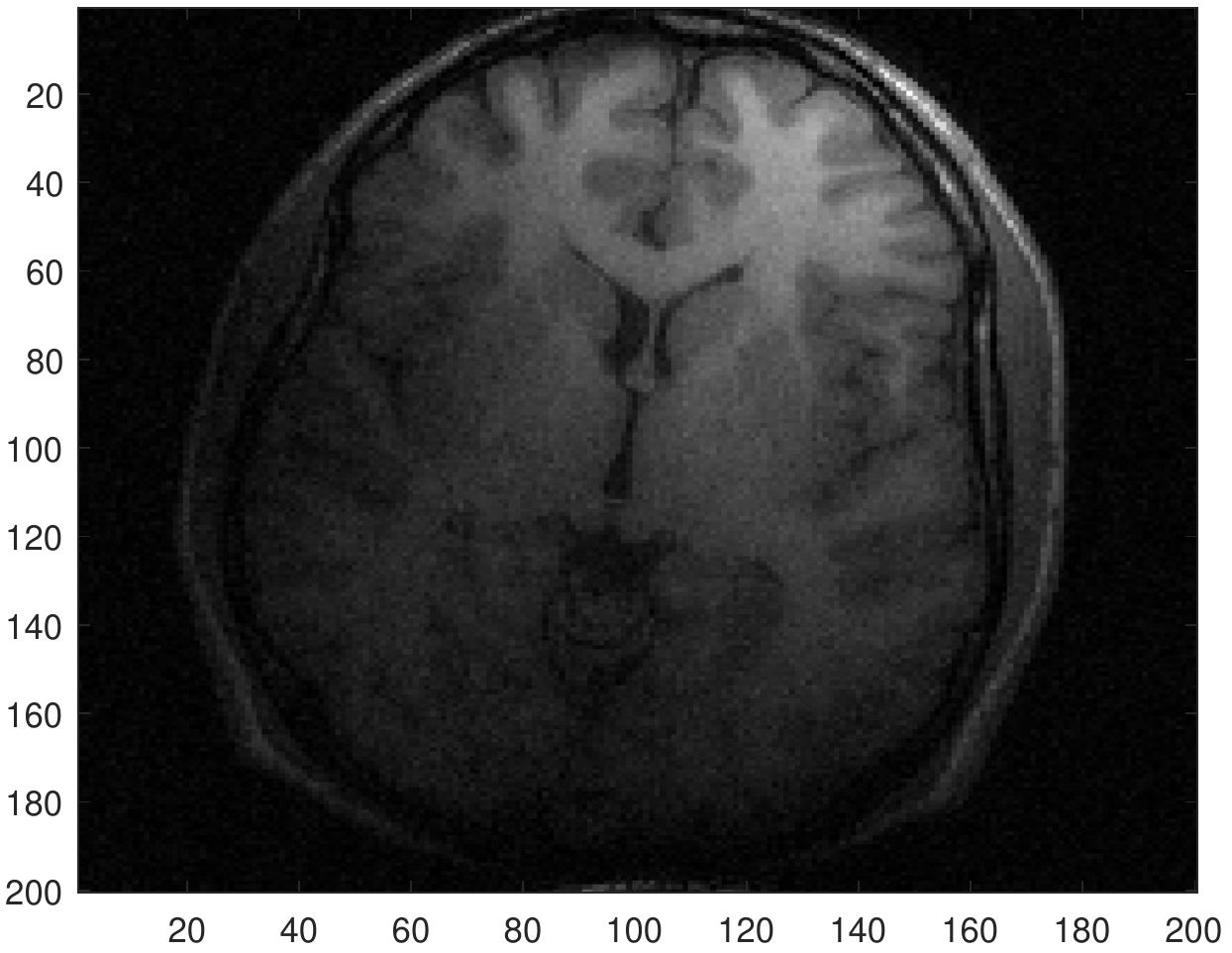}
\subcaption{$|s_1F|$}
\end{subfigure}
\begin{subfigure}{0.24\textwidth}
\includegraphics[width=0.9\linewidth, height=3.2cm, keepaspectratio]{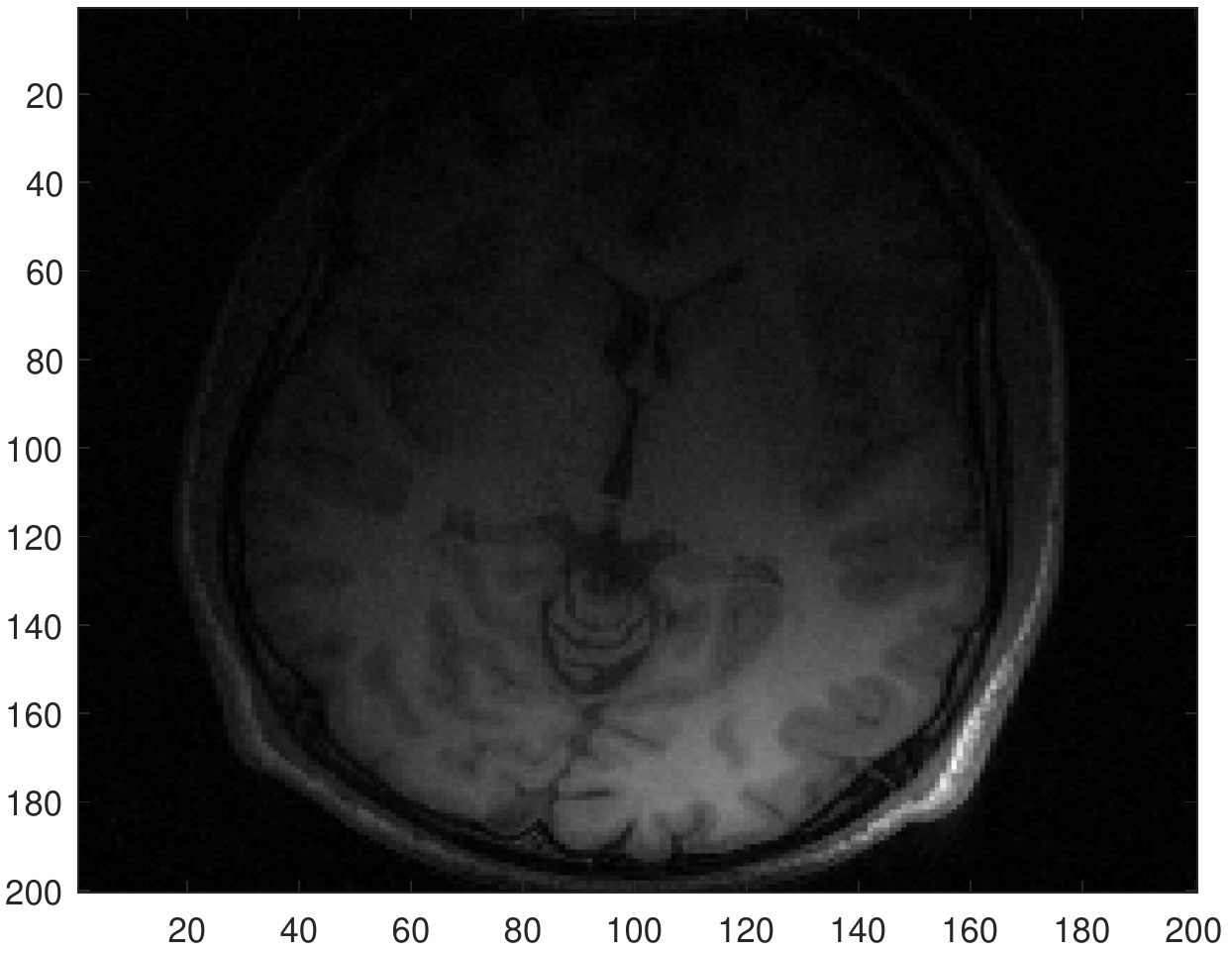}
\subcaption{$|s_4F|$}
\end{subfigure}
\begin{subfigure}{0.24\textwidth}
\includegraphics[width=0.9\linewidth, height=3.2cm, keepaspectratio]{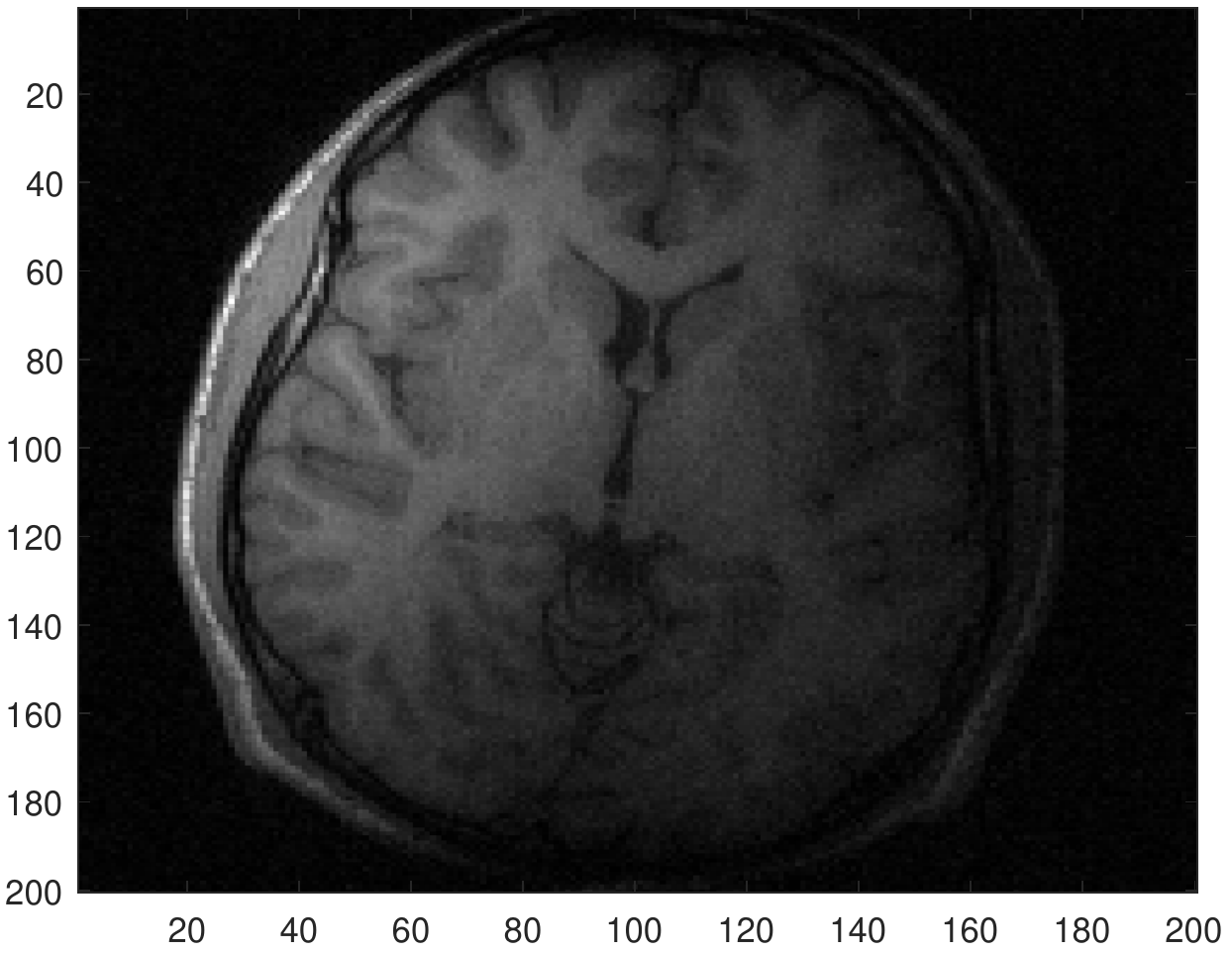}
\subcaption{$|s_{7}F|$}
\end{subfigure}
\caption{Multi-coil brain image. In total, there are $K=8$ coils equally spaced on a ring around the head. The coil sensitivities $s_1,\ldots,s_K$ are ordered clockwise around the head. We show the ground truth (left-hand figure) and three coil images.}
\label{bp}
\end{figure}

\subsection{How the number of coils affects inversion stability}
\label{coils}
In this section, we investigate the relation between the number of coils and inversion stability, using the SVD of $M$. For this example, we consider brain image shown in figure \ref{bp}. The data was downloaded from the BART toolbox \cite{BART}. In this example, $n=m=200$, and the number of coils is $K=8$. We consider accelerated sampling, at rates of $R=2,3,4$. $R=n_s$ acceleration means that 1 out of every $n_s$ $\textbf{k}$-space lines are collected during the scan with uniform spacing, while retaining (fully scanning) a central region of $\textbf{k}$-space of pre-specified width. In this paper, we retain 32 central lines of $\textbf{k}$-space, for calibration. For this analysis, the sensitivity maps $s_1,\ldots,s_K$ are approximated using ESPIRiT \cite{ESP}. 

\begin{figure}
\centering
\begin{subfigure}{0.24\textwidth}
\includegraphics[width=0.9\linewidth, height=3.2cm, keepaspectratio]{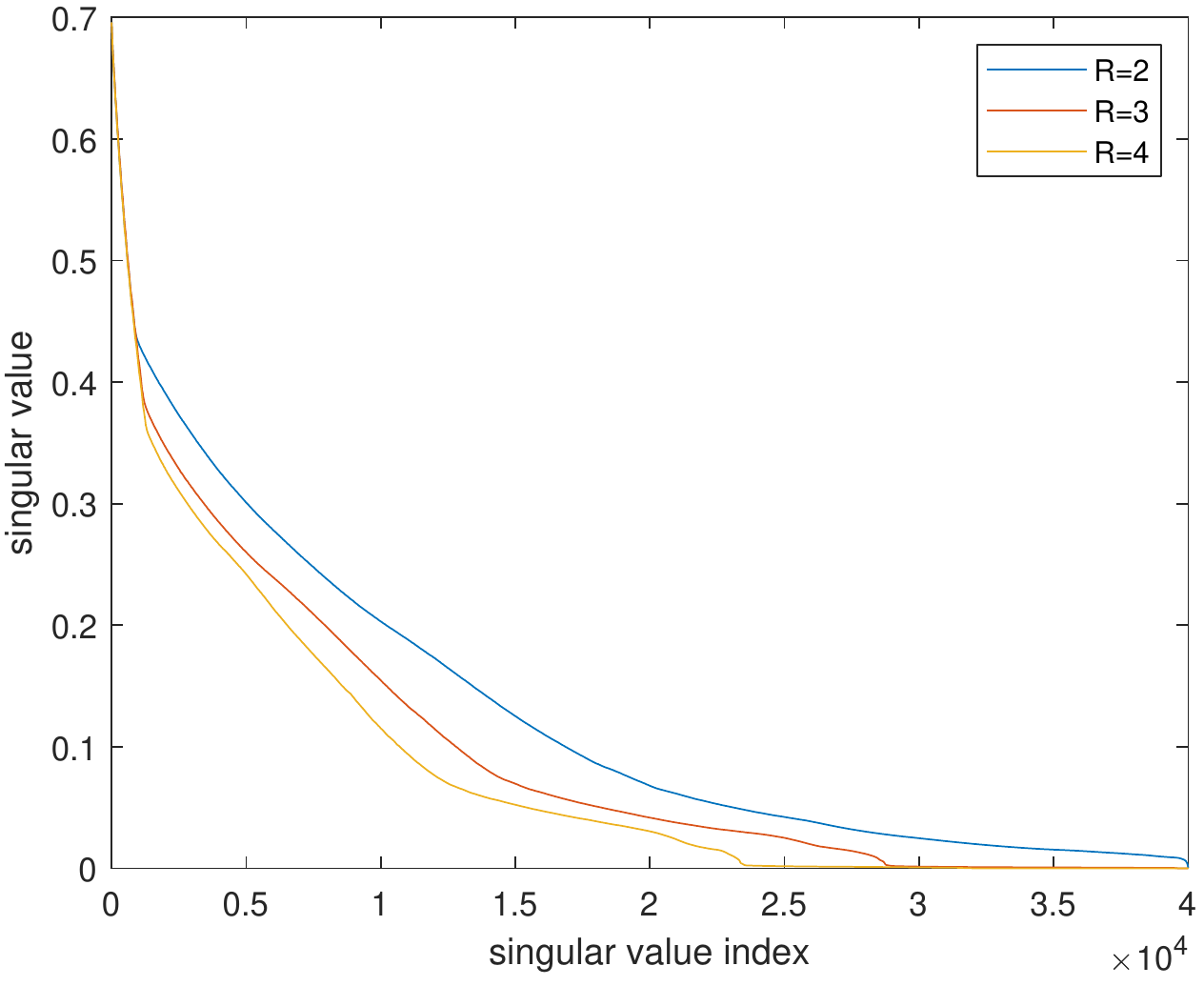}
\subcaption{$C_1$}
\end{subfigure}
\begin{subfigure}{0.24\textwidth}
\includegraphics[width=0.9\linewidth, height=3.2cm, keepaspectratio]{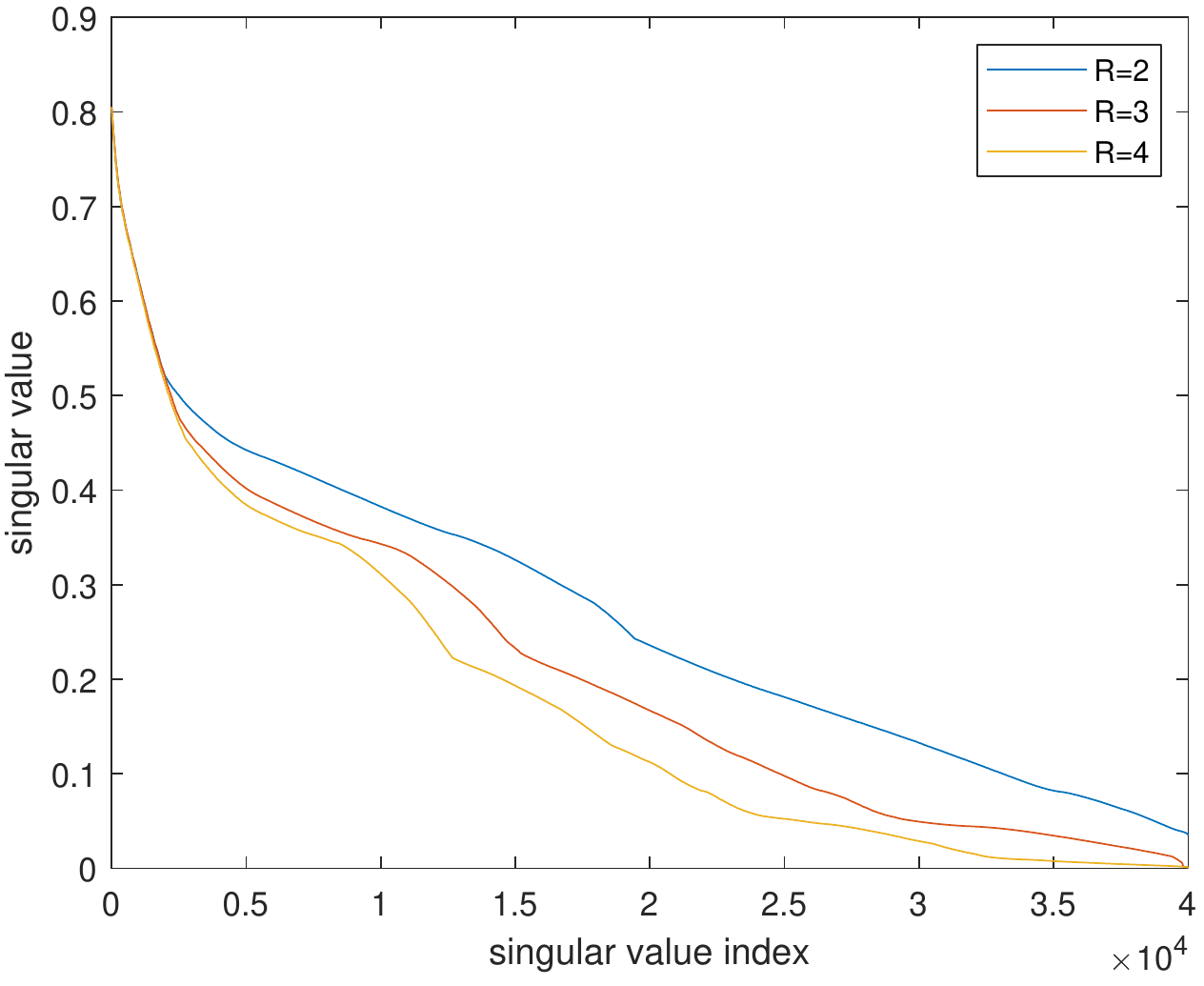}
\subcaption{$C_2$}
\end{subfigure}
\begin{subfigure}{0.24\textwidth}
\includegraphics[width=0.9\linewidth, height=3.2cm, keepaspectratio]{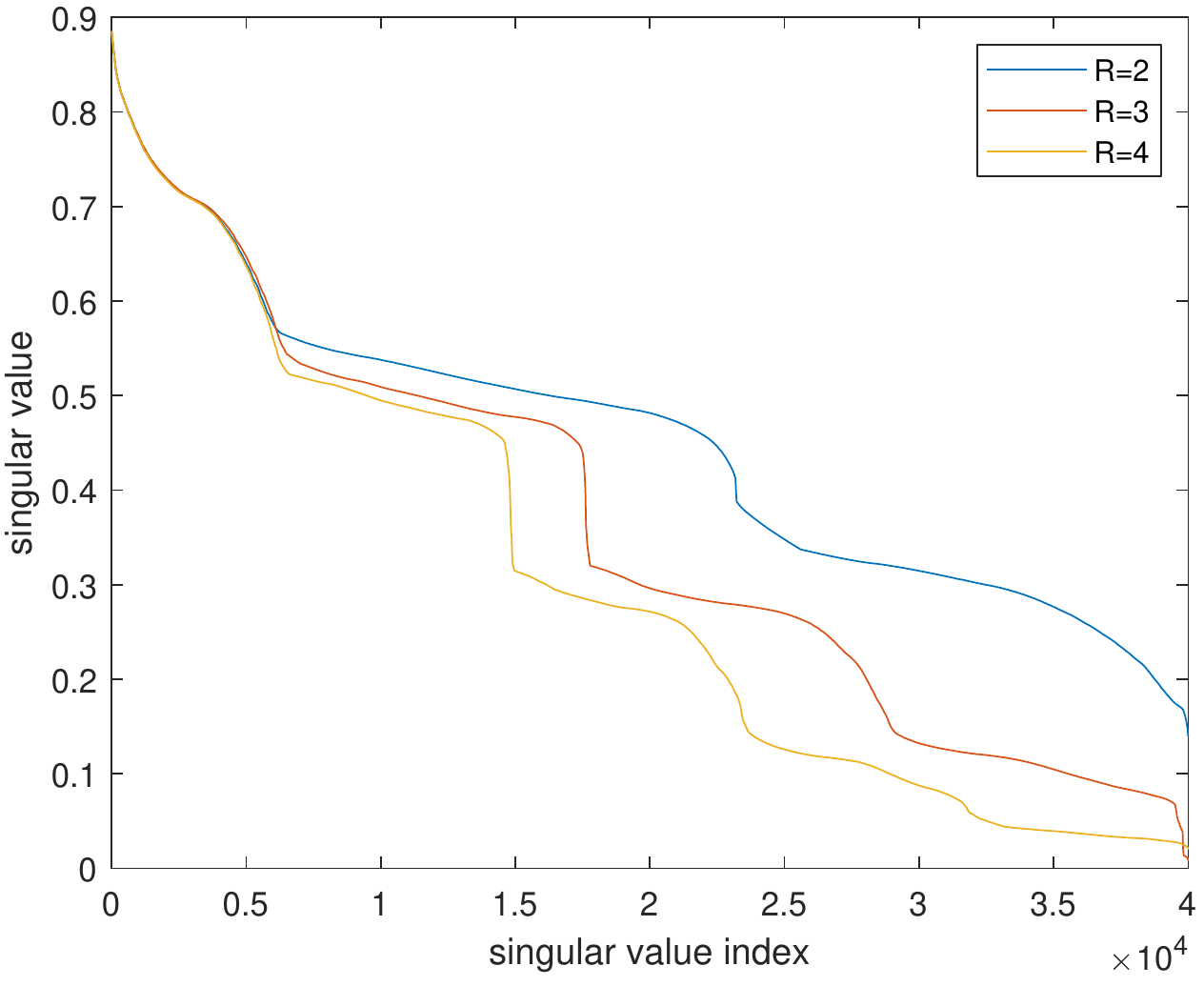}
\subcaption{$C_3$}
\end{subfigure}
\begin{subfigure}{0.24\textwidth}
\includegraphics[width=0.9\linewidth, height=3.2cm, keepaspectratio]{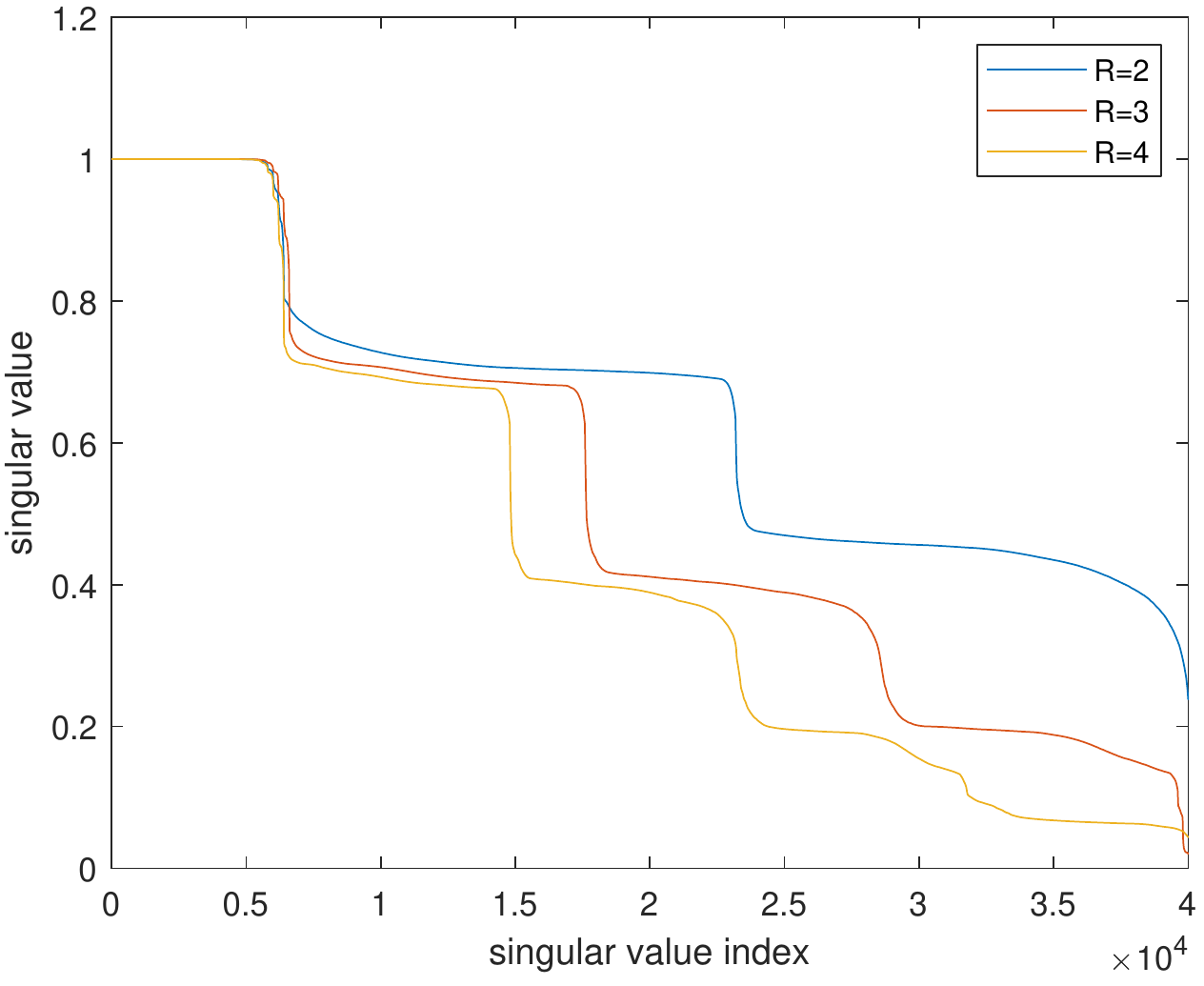}
\subcaption{$C_4$}
\end{subfigure}
\begin{subfigure}{0.31\textwidth}
\includegraphics[ width=1\linewidth, height=1\linewidth, keepaspectratio]{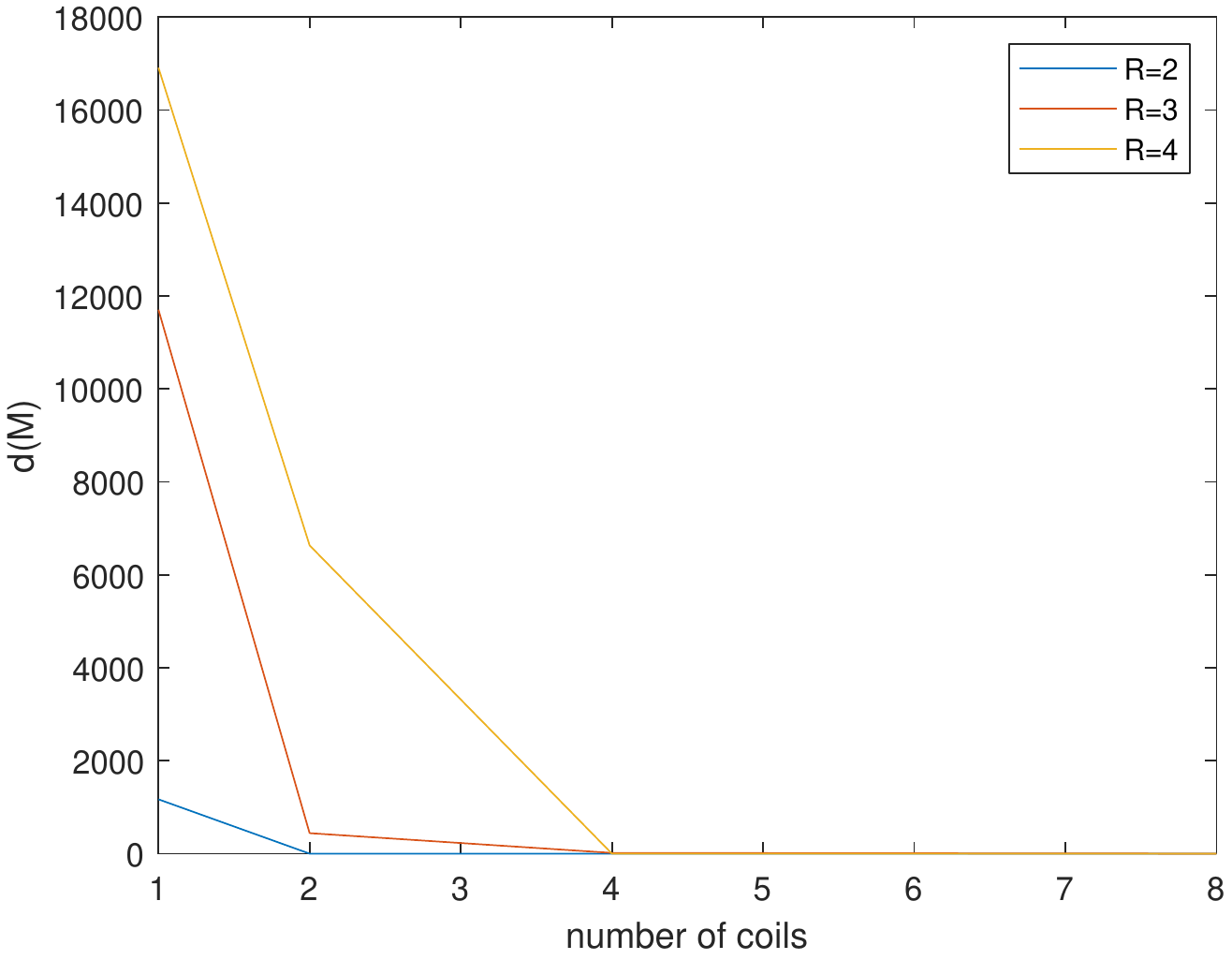} 
\subcaption{$d(M)$} \label{staba}
\end{subfigure}
\hspace{0.3cm}
\begin{subfigure}{0.31\textwidth}
\includegraphics[ width=1\linewidth, height=1\linewidth, keepaspectratio]{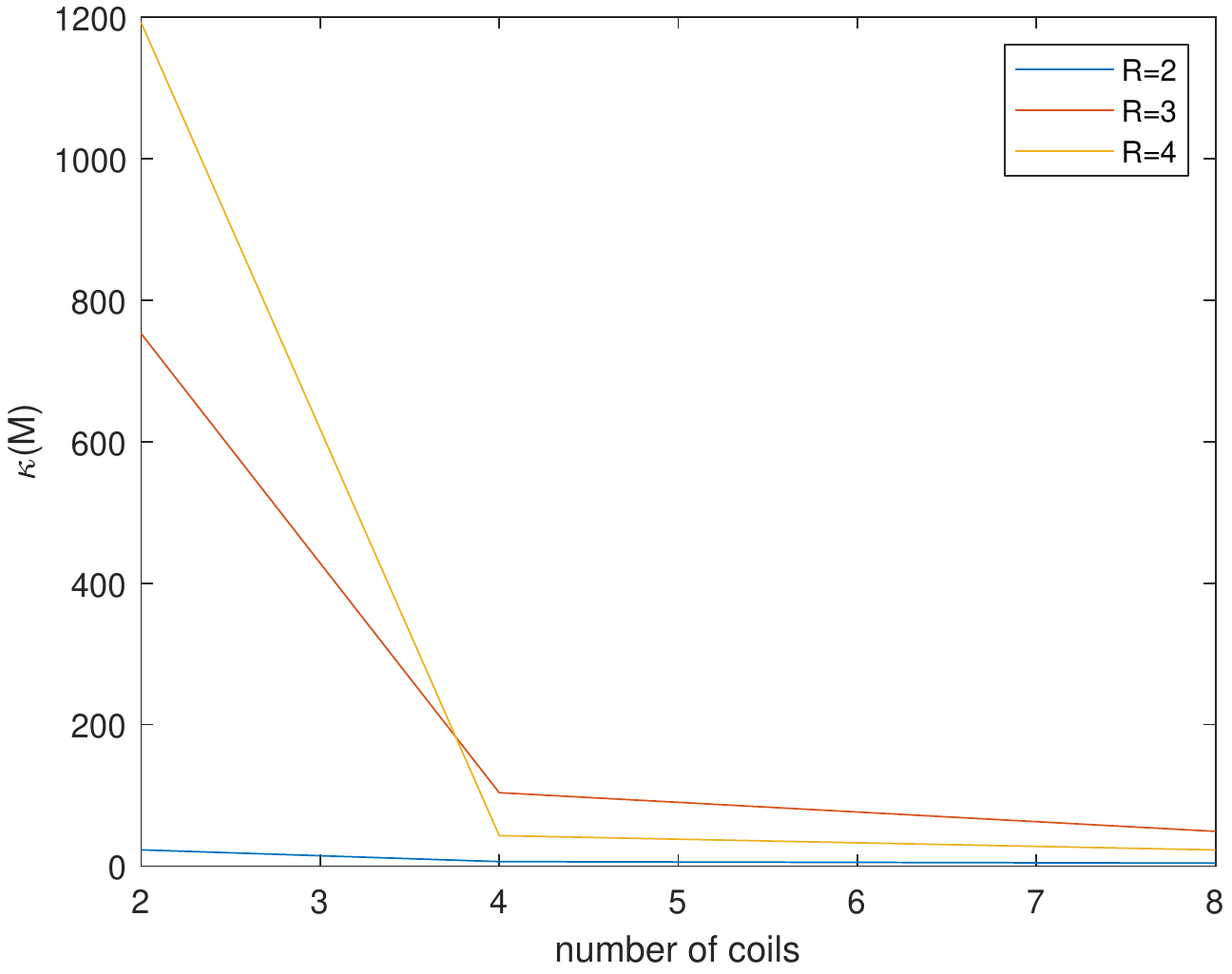} 
\subcaption{$\kappa(M)$} \label{stabb}
\end{subfigure}
\caption{Top row -  Plots of the singular values of $M$ for varying coil subsets and acceleration factors. Bottom row - Stability metric plots for varying $R$ and $K$. The effective null space dimension threshold is set at $t=0.01$.}
\label{singular}
\end{figure}

We consider four subsets of coils such that the coil spacing remains even around the head. Specifically, we consider the subsets $C_1=\{1\}$ (singleton coil), $C_2=\{1,5\}$ (two coils opposite one another), $C_3=\{1,3,5,7\}$ (four coils uniformly spaced), and $C_4=\{1,\ldots,8\}$ (all coils). In figure \ref{singular} (top row), we show plots of the singular values of $M$ for each coil subset, and for the acceleration factors $R=2,3,4$. In figure \ref{singular} (bottom row), we show plots of $d(M)$ and $\kappa(M)$ with the number of coils, for $R=2,3,4$. 

We see a clear increase in the problem stability in terms of $d(M)$ and $\kappa(M)$, as the number of coils increases. We see a similar effect as $R$ decreases, as we would expect, since there is less missing data. As the number of coils decreases, the rate of decay of the singular values appears to increase. This is verified by the $d(M)$ values shown in figure \ref{staba}. In this example, four coils are sufficient to for zero effective null space (i.e., $d(M)=0$), for $R=2,3,4$. When $R=2$ and there is less missing data, two coils are sufficient for zero null space, and the $\kappa(M)$ values small relative to $R=3,4$. See figure \ref{stabb}. When $R=4$ the $\kappa(M)$ values are smaller than $R=3$ when the number of coils is greater than or equal to 4. The dimension of the effective null space is larger, in the $R=4$ case, however.


\subsection{Sub-sampling scheme comparison}
\label{subsamp}
In this section, we use the SVD of $M$ to analyze the problem stability for different $\textbf{k}$-space sub-sampling patterns. Specifically, we compare conventional $R=2,3,4$ acceleration with uniform random undersampling. Random undersampling is used to simulate $\textbf{k}$-space corruption due to motion, as is, for examples, done in \cite{motion1,motion2}. As in section \ref{coils}, we consider the brain image of figure \ref{bp} in the examples presented.

\begin{figure}[!h]
\centering
\begin{subfigure}{0.32\textwidth}
\includegraphics[width=0.9\linewidth, height=3.5cm, keepaspectratio]{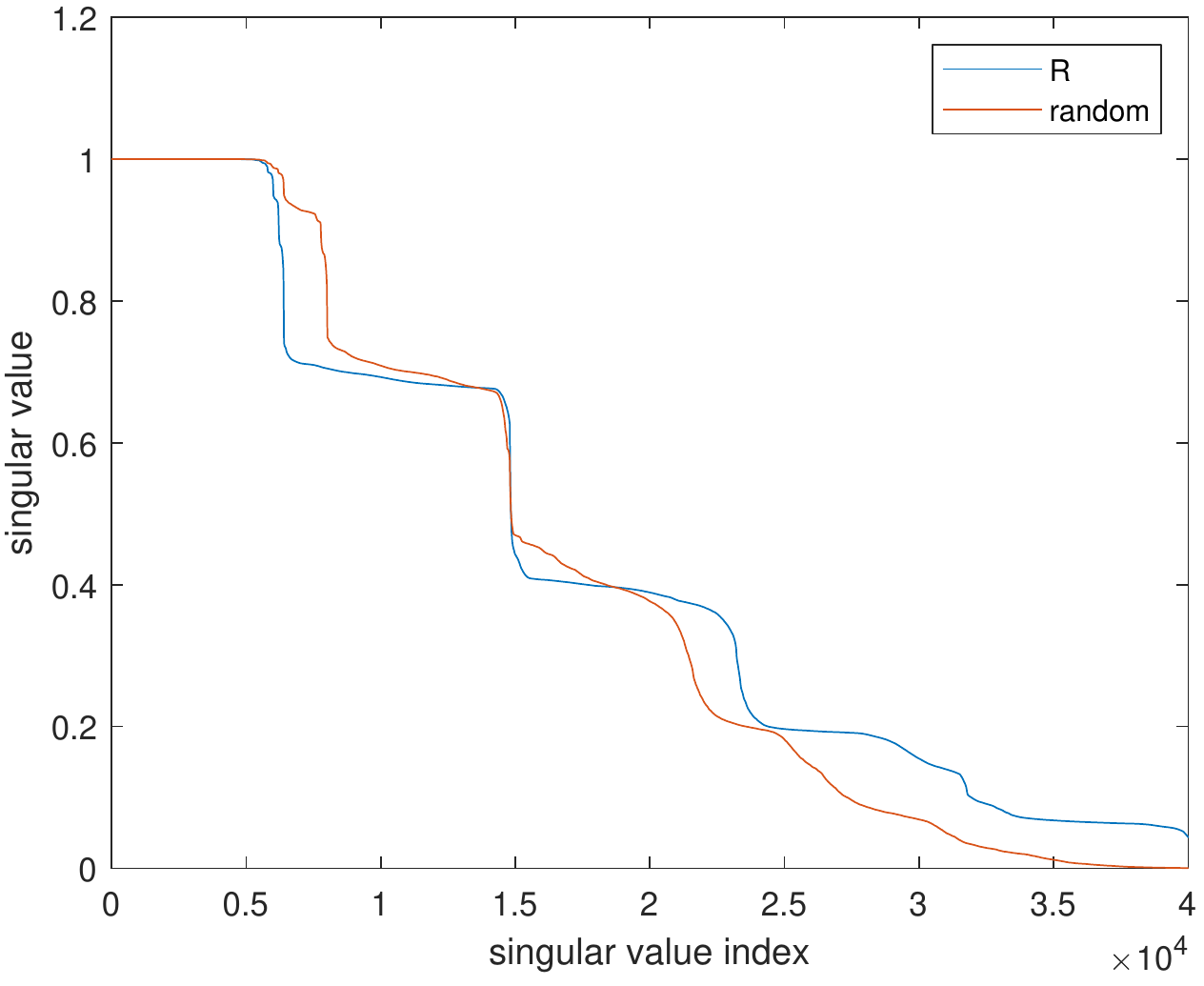}
\subcaption{$\text{scan time} = 37\%$ ($R=4$)}
\end{subfigure}
\begin{subfigure}{0.32\textwidth}
\includegraphics[width=0.9\linewidth, height=3.5cm, keepaspectratio]{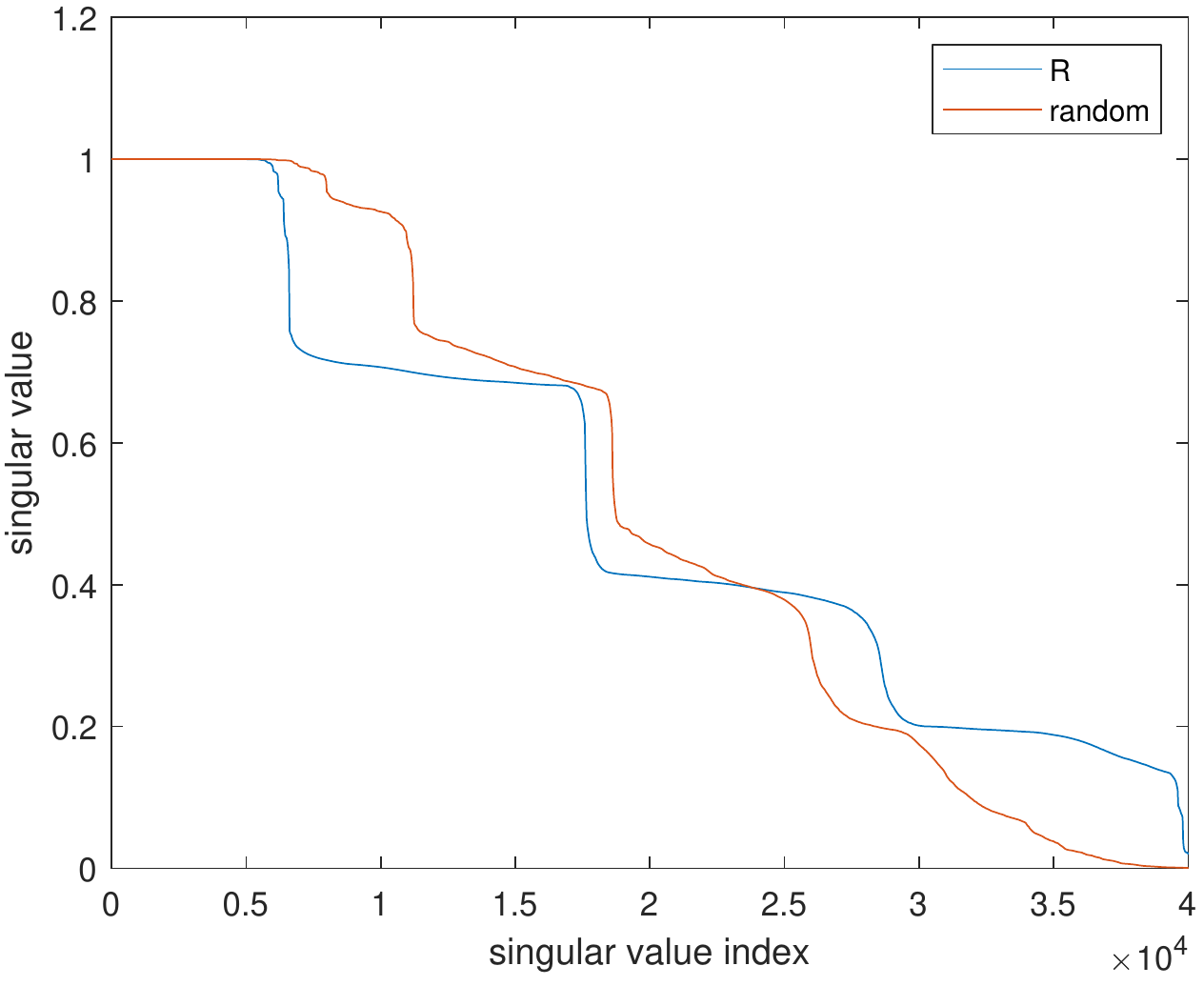}
\subcaption{$\text{scan time} = 44.5\%$ ($R=3$)}
\end{subfigure}
\begin{subfigure}{0.32\textwidth}
\includegraphics[width=0.9\linewidth, height=3.5cm, keepaspectratio]{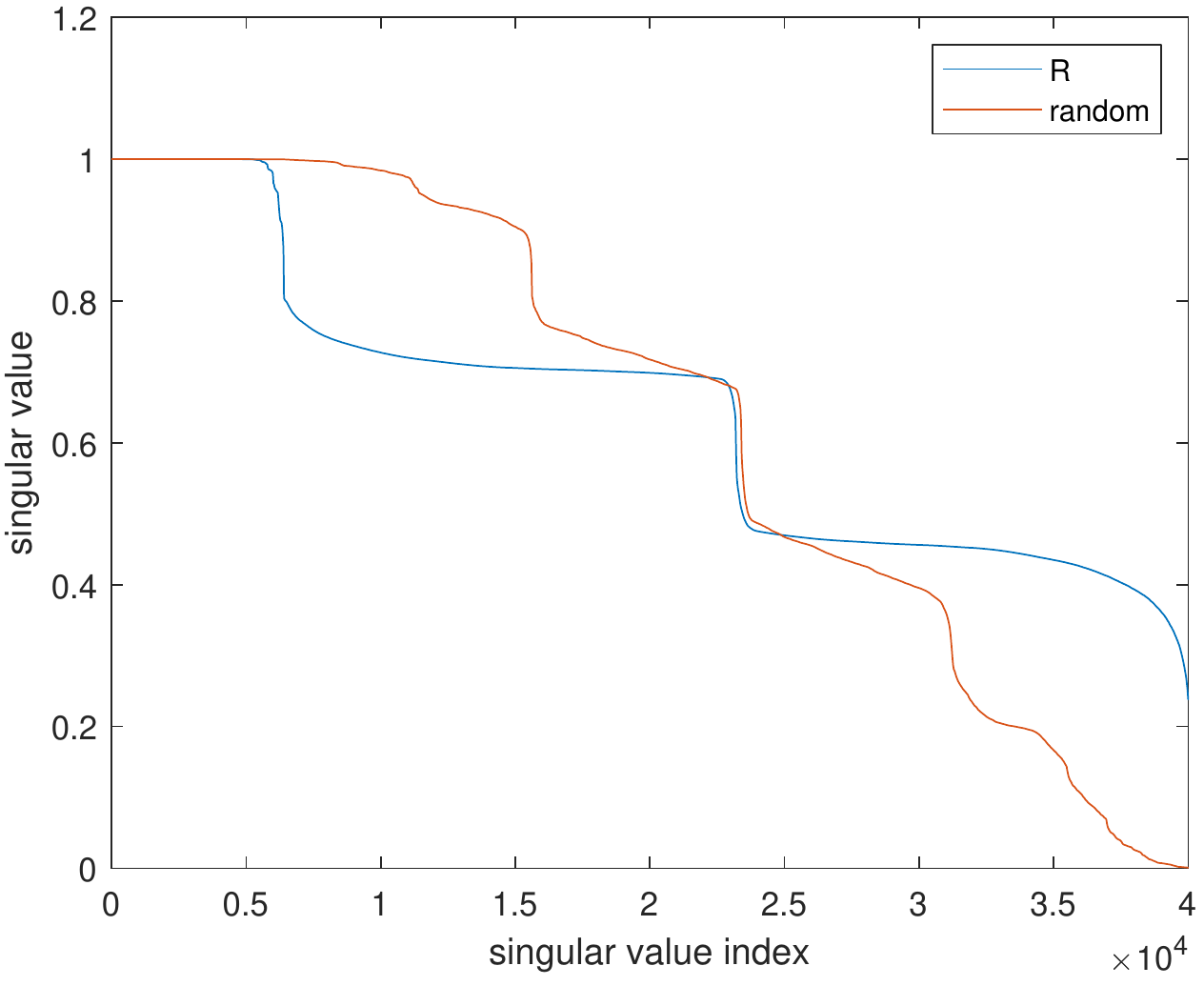}
\subcaption{$\text{scan time} = 58\%$ ($R=2$)}
\end{subfigure}
\begin{subfigure}{0.31\textwidth}
\includegraphics[ width=1\linewidth, height=1\linewidth, keepaspectratio]{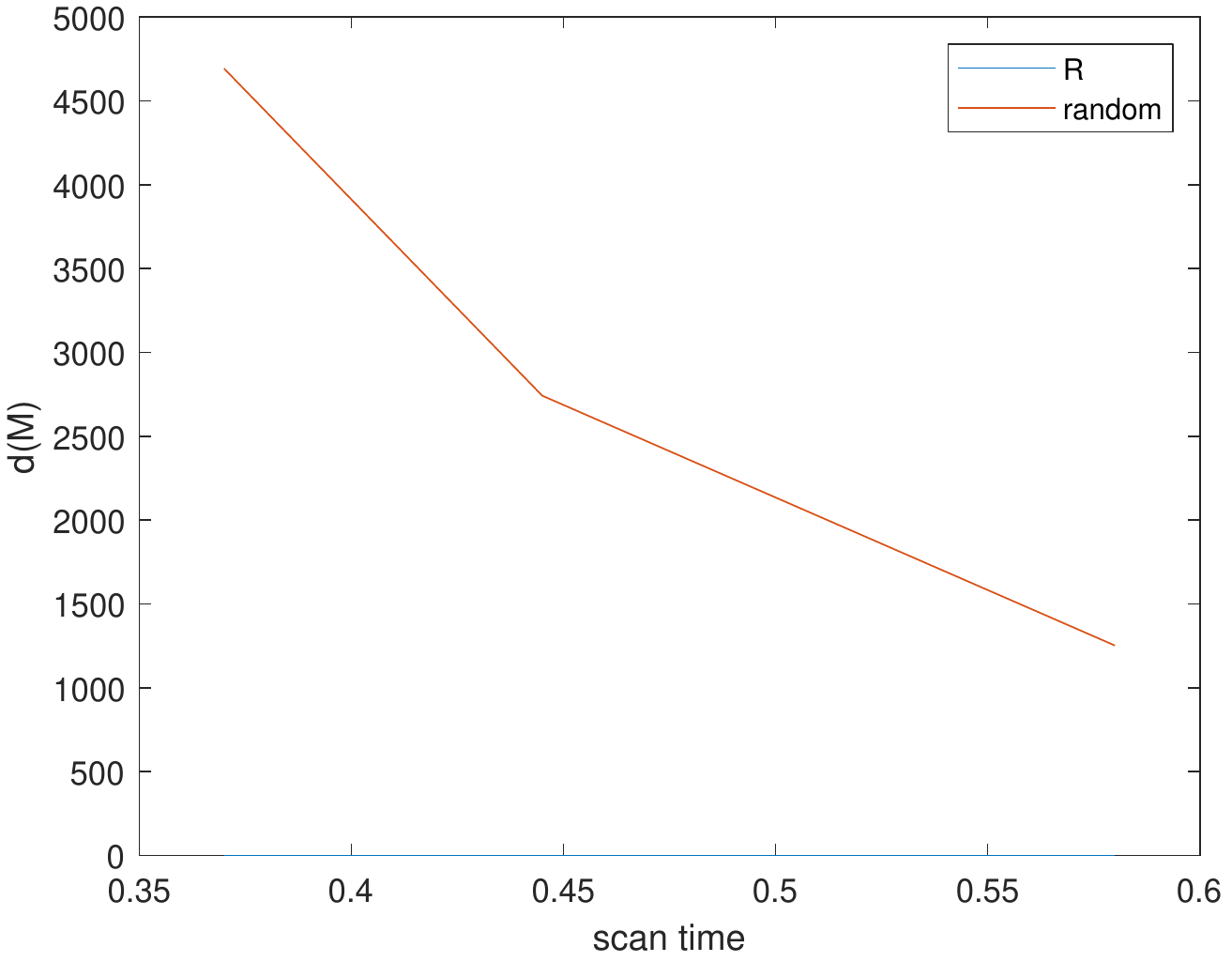}
\subcaption{$d(M)$}
\end{subfigure}
\hspace{0.3cm}
\begin{subfigure}{0.31\textwidth}
\includegraphics[ width=1\linewidth, height=1\linewidth, keepaspectratio]{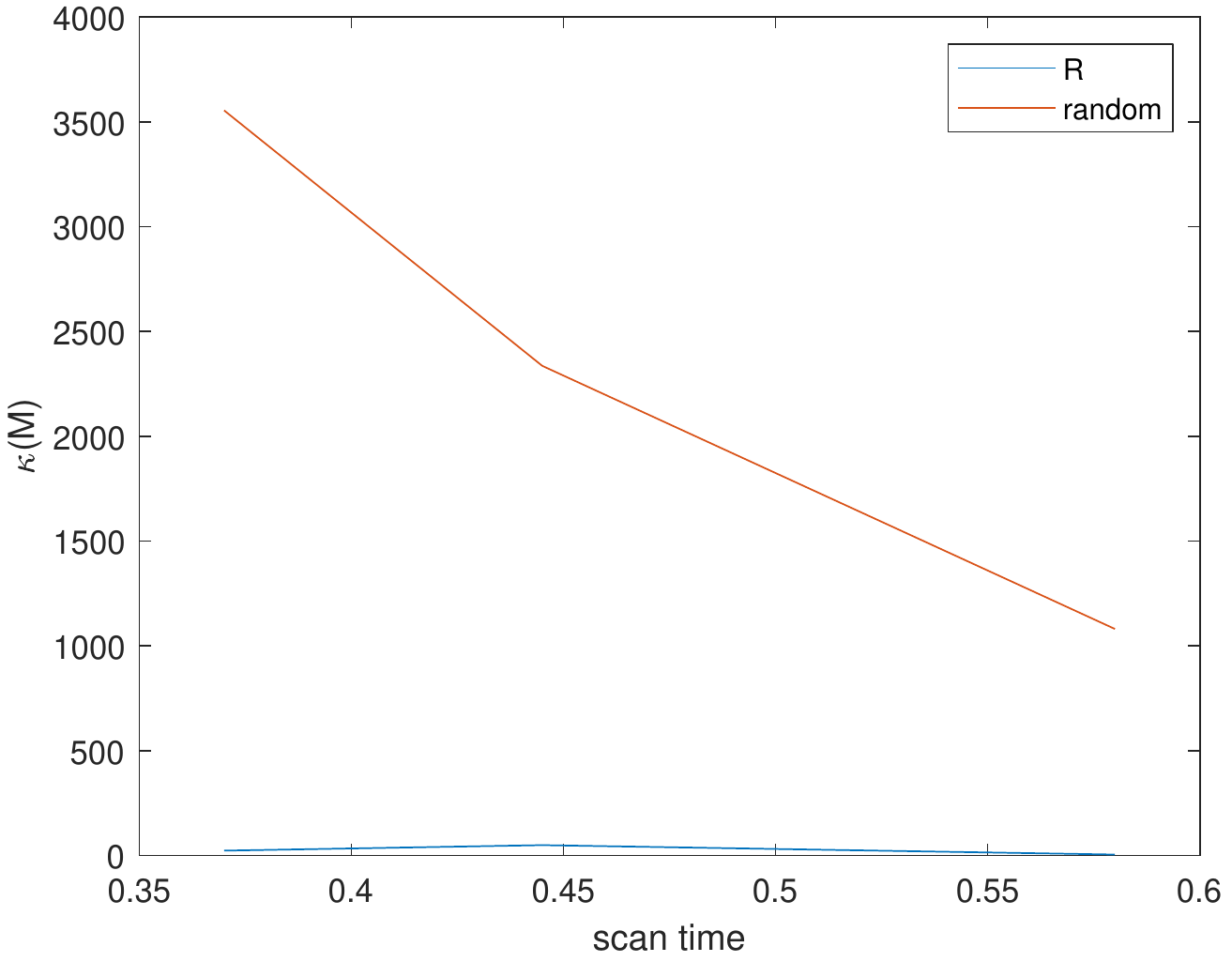}
\subcaption{$\kappa(M)$}
\end{subfigure}
\caption{Top row - Singular value plot comparison for random vs accelerated sub-sampling, for scan times corresponding to $R=2,3,4$. Bottom row - Condition number and null space dimension comparison for random vs accelerated sampling. The scan times on the $x$ axis are those of figure \ref{S_rvsR}. Note, the blue curves (corresponding to $R$) are smaller magnitude than the red curves, and appear close to the $x$ axis in the plots.}
\label{S_rvsR}
\end{figure}

In figure \ref{S_rvsR} (top row),  we present singular value plot comparisons for random vs accelerated sampling, for $R=2,3,4$. For both sub-sampling schemes, we retain 32 central lines for calibration. The scan times are given in the figure sub-caption, and correspond to $R=2,3,4$ acceleration rates. The scan time is the percentage of total lines retained for scanning. Explicitly, the scan time is calculated $\text{scan time}=(32+l)/n$, where $l\geq 0$ is the number of lines retained for scanning. In figure \ref{S_rvsR} (bottom row), we show $d(M)$ and $\kappa(M)$ plots comparing random and accelerated sampling for varying scan times. The plots in figure \ref{S_rvsR} imply greater $M$ inversion stability using accelerated sampling, compared to random sampling, for all scan times considered. Thus, we would expect to see greater overall noise amplification when the $\textbf{k}$-space lines  are sub-sampled at random (e.g., when there is motion corruption), when compared to uniform, accelerated sampling.

For accelerated sampling (e.g., $R=2,3,4$), the missing blocks of $\textbf{k}$-space have uniform width $w_j$ (refer to \eqref{equ2}). When the locations of the missing $\textbf{k}$-space lines are sampled at random from a uniform distribution, and $n$ is small (in this case $n=200$), the missing $\textbf{k}$-space lines are likely to cluster together to form larger $w_j$ blocks, which are harder to recover. 
Noise amplification is not the only aspect one should consider, however, in stability analysis. While the overall error amplification with random sampling is greater, when compared to uniform acceleration, this does not account for the distribution of artifacts within the image. To investigate this further, we analyze the right singular vectors of $M$, specifically those which correspond to the smallest singular values. The right singular vectors of $M$ which correspond to the smallest singular values (e.g., $\sigma_{nm}$) span the effective null space of $M$, and provide insight into the distribution of image artifacts due to null space.

\begin{figure}[!h]
\centering
\begin{subfigure}{0.24\textwidth}
\includegraphics[width=0.9\linewidth, height=3.2cm, keepaspectratio]{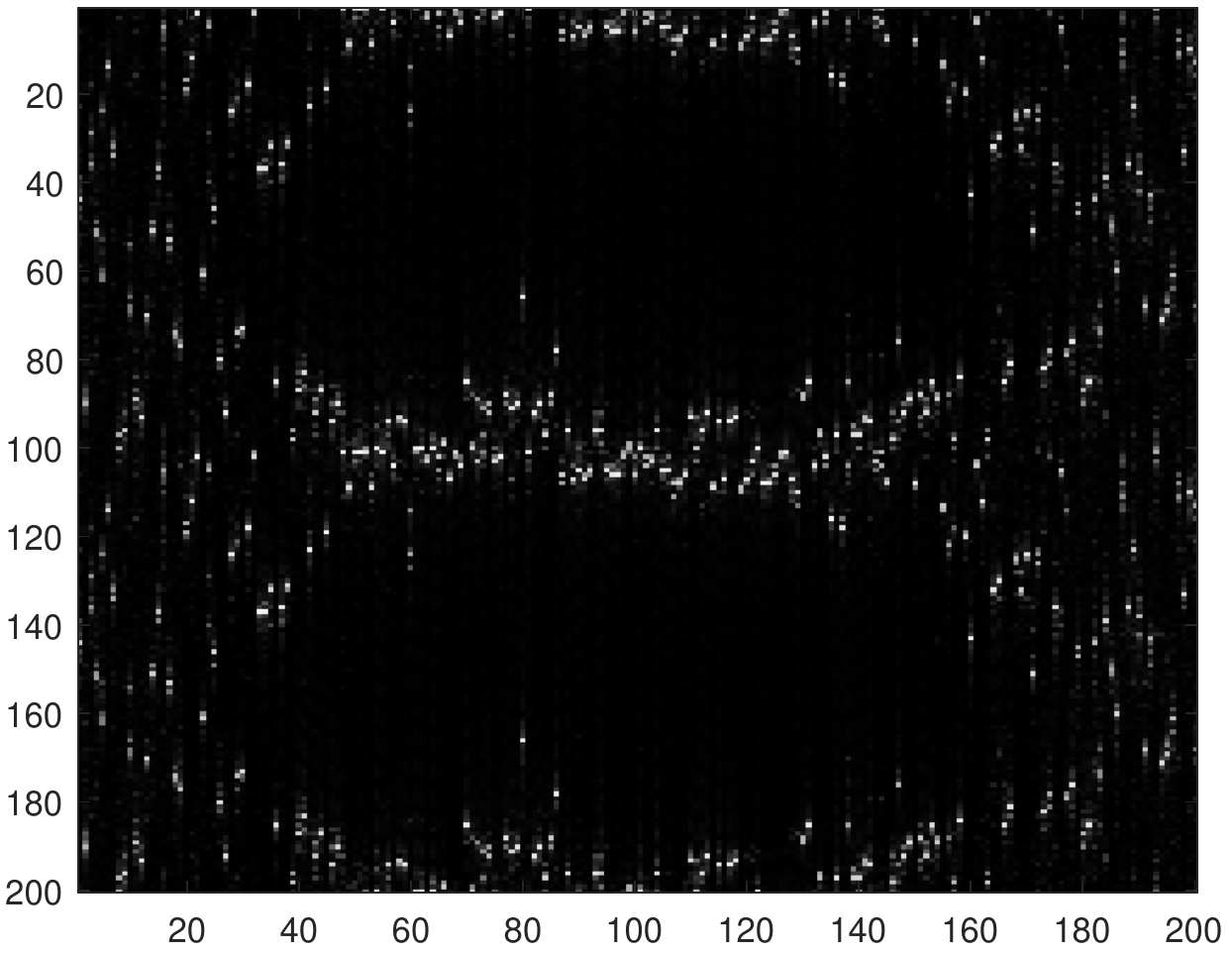}
\subcaption{right singular vectors ($R=2$)}\label{ba1}
\end{subfigure}
\begin{subfigure}{0.24\textwidth}
\includegraphics[width=0.9\linewidth, height=3.2cm, keepaspectratio]{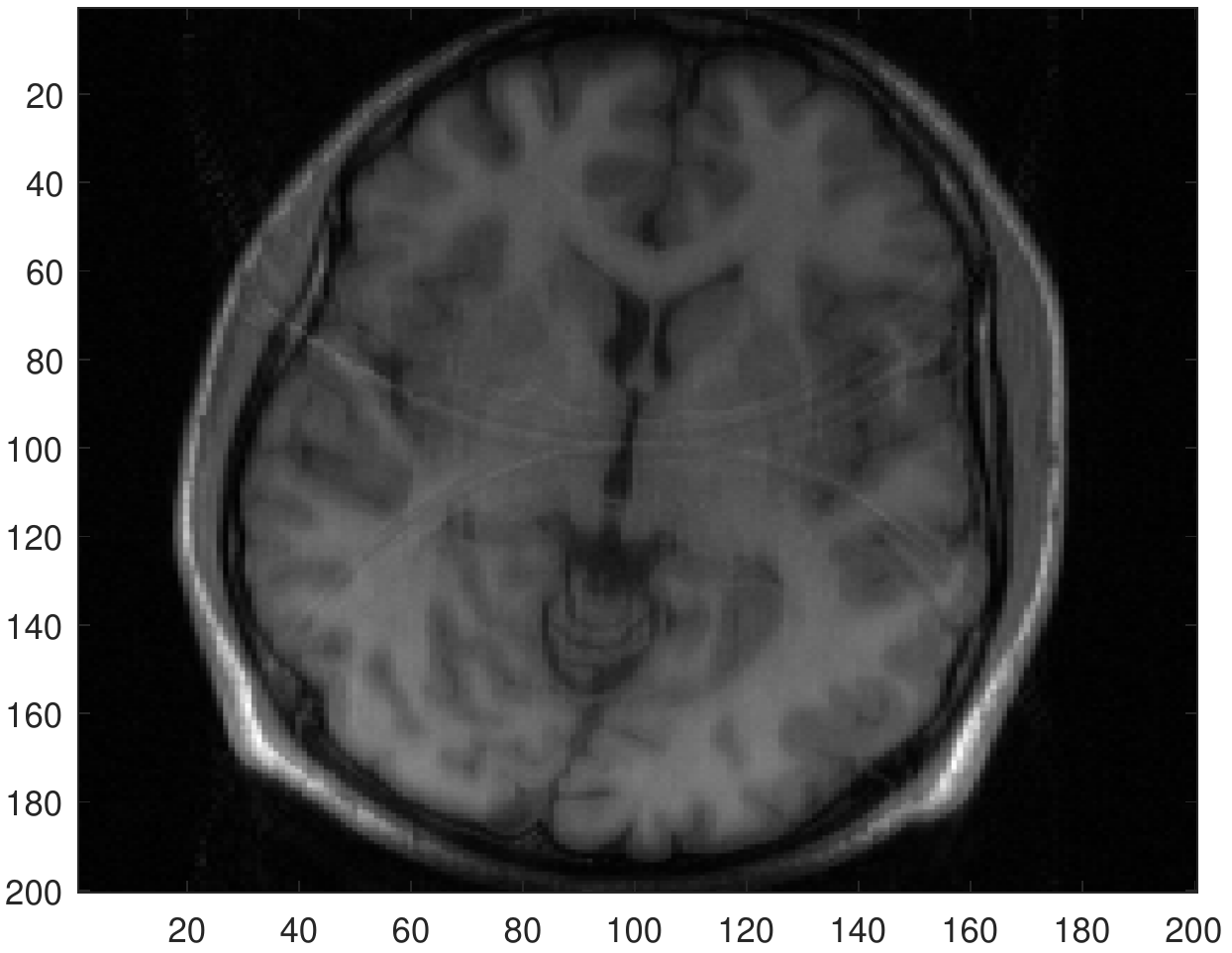}
\subcaption{image reconstruction ($R=2$)} \label{be1}
\end{subfigure}
\begin{subfigure}{0.24\textwidth}
\includegraphics[width=0.9\linewidth, height=3.2cm, keepaspectratio]{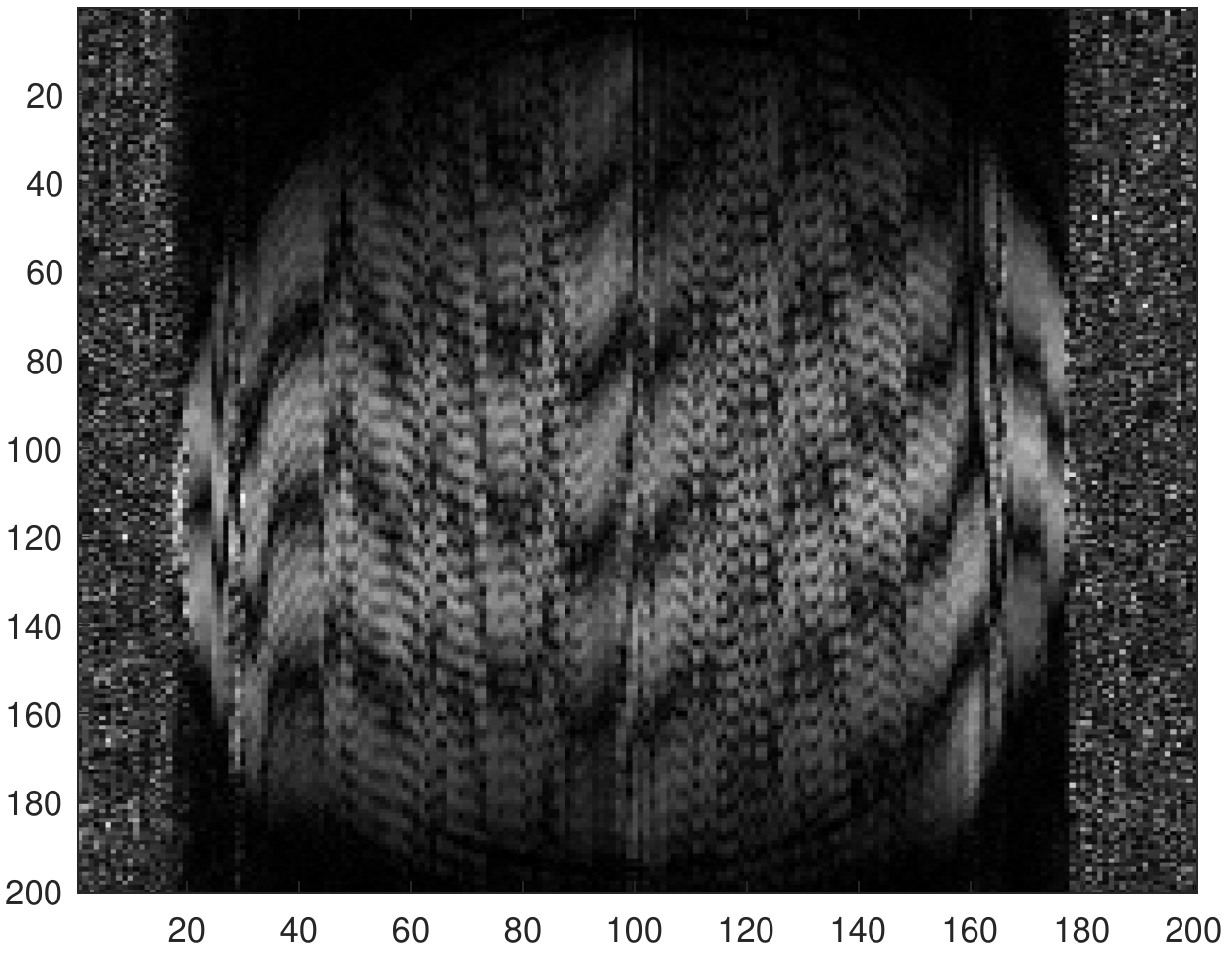}
\subcaption{right singular vectors (random sampling)} \label{bar1}
\end{subfigure}
\begin{subfigure}{0.24\textwidth}
\includegraphics[width=0.9\linewidth, height=3.2cm, keepaspectratio]{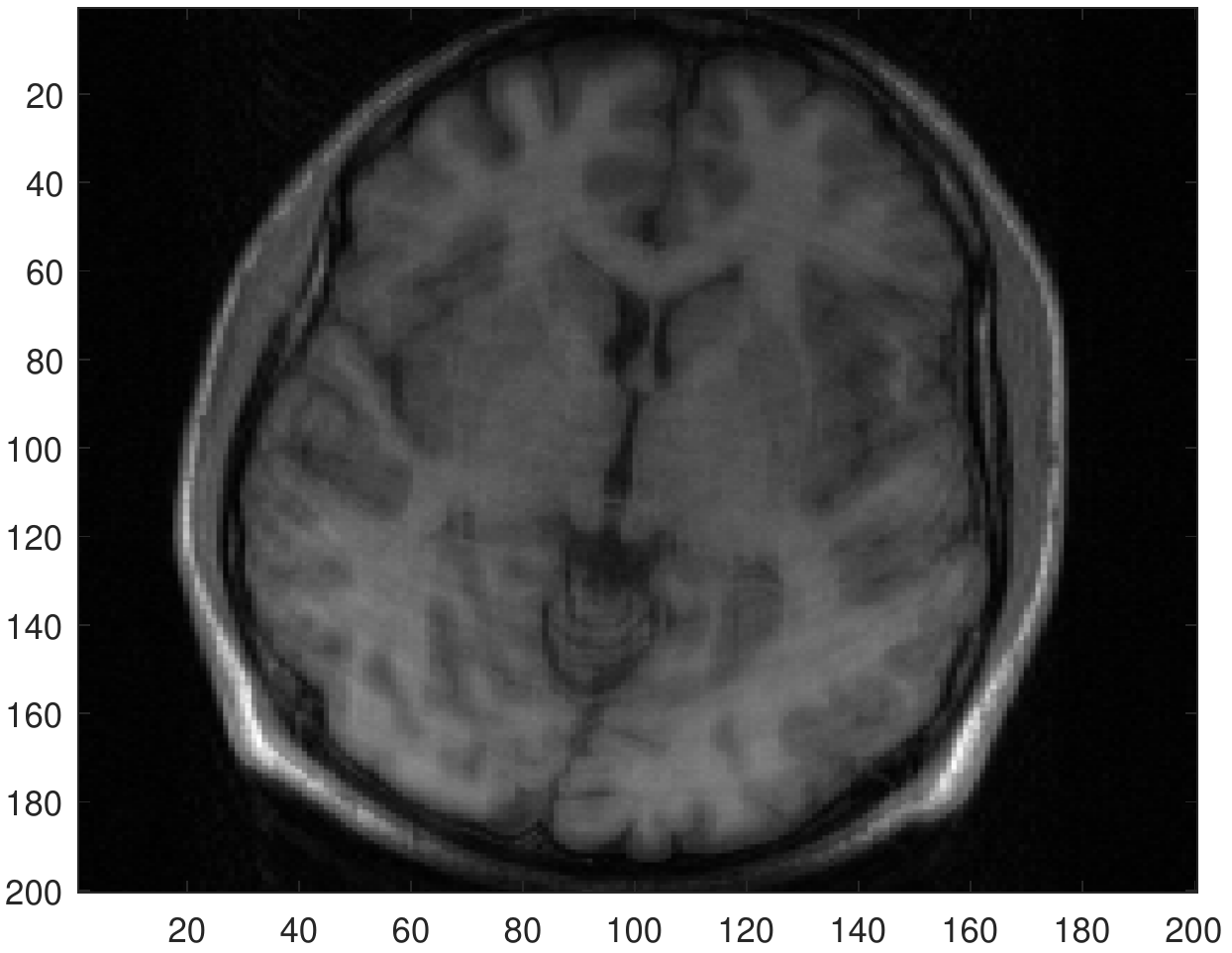}
\subcaption{image reconstruction (random sampling)} \label{ber1}
\end{subfigure}
\caption{Right singular vector images and SOS image reconstructions with missing data, setting missing $\textbf{k}$-space to zero. We consider the case of accelerated sampling with $R=2$ and random sampling with $\text{scan time} = 58\%$ (this is equivalent to $R=2$).}
\label{RsingV_R}
\end{figure}

See figure \ref{RsingV_R} where we have shown right singular vector images and limited-data brain reconstructions for $R=2$ sampling, and random sampling at the same scan time. To clarify, for each $M$ considered, the right singular vector images are the matrices $\paren{\vv_{1n},\ldots,\vv_{mn}}$, where $\vv_{in}$ is the $n_{\text{th}}$ (i.e., the smallest) singular vector of $A_i$. For $R=2$, there are strong aliasing artifacts which appear as sharp curves through the center of the brain. See figure \ref{be1}. The artifacts are highlighted also in the right singular vector image in figure \ref{ba1}. For random under-sampling in figure \ref{ber1},  the artifacts are less focused to a particular spatial region as in the $R=2$ case, and we see a more general blurring (``shaking") effect in the reconstruction, as is often seen in applications with motion error \cite{GRAPPA,motion1,motion2}. The right singular vectors corresponding to the smallest singular values show a similar blurring/shaking effect in figure \ref{bar1}. While the overall noise amplification is reduced for $R=2$ and accelerated sampling, the artifacts are stronger and more localized when compared to random undersampling.

Random and accelerated sampling were chosen as examples of interest here, to highlight the potential benefits of the Fredholm operator SVD to analyze the stability of multi-coil MRI problems. The SVD analysis is not specific to random or accelerated sampling, however, and can be applied to any sub-sampling scheme where the missing lines of $\textbf{k}$-space are parallel to either $k_1$ or $k_2$. 



\section{Methods}
\label{results}
Here we and generalize the discrete formulation of \eqref{1D} to multiple ESPIRiT sensitivity maps and detail our reconstruction method, whereby the 2-D image is reconstructed slice-by-slice on lines parallel to the $x_2$ axis. We test our method on several data sets in {\it{in vivo}} MRI (e.g., brain, spine, cardiac) when the lines of missing $\textbf{k}$-space are sampled at random.

 Let $s^{(1)}_j,\ldots,s^{(p)}_j$ be a set of $p$ ESPIRiT maps corresponding to coil $j$, approximated using the algorithm of \cite{ESP}. Then, we assume that $|F|$ may be decomposed in the form
\begin{equation}
\label{|F|}
|F|=\sqrt{\sum_{q=1}^p|F_q|^2},
\end{equation}
where $F_q$ is the component of $F$ corresponding to map $q$. Let 
\begin{equation}
{\bf{F}}^{(i)}_q=\begin{pmatrix}
F_q\paren{u_i,v_1} \\
F_q\paren{u_i,v_2} \\
\vdots \\
F_q\paren{u_i,v_n}
\end{pmatrix}\in\mathbb{C}^n,\ \ \text{and}\ \ 
\vb_i=\begin{pmatrix}
\vb^{(1)}_{u_i} \\
\vb^{(2)}_{u_i} \\
\vdots \\
\vb^{(k)}_{u_i}
\end{pmatrix}\in\mathbb{C}^{nk}.
\end{equation}
Then, the generalization of \eqref{1D} to multiple sensitivities becomes
\begin{equation}
C_i{\bf{F}}_i=\vb_i,
\end{equation}
where
\begin{equation}
C_i=\begin{bmatrix}
A\left[S^{(1,1)}_{i},\ldots,S^{(1,p)}_{i}\right] \\
A\left[S^{(2,1)}_{i},\ldots,S^{(2,p)}_{i}\right] \\
\vdots \\
A\left[S^{(k,1)}_{i},\ldots,S^{(k,p)}_{i}\right]
\end{bmatrix}\in \mathbb{C}^{nk\times np}, \ \ \text{and} \ \  
{\bf{F}}_i=\begin{pmatrix}
{\bf{F}}^{(i)}_1 \\
{\bf{F}}^{(i)}_2 \\
\vdots \\
{\bf{F}}^{(i)}_p
\end{pmatrix}\in\mathbb{C}^{np},
\end{equation}
where
$$S^{(j,q)}_{i}= \text{diag}\paren{s^{(q)}_j(u_i,v_1),\ldots, s^{(q)}_j(u_i,v_n)}.$$
To recover ${\bf{F}}_i$, we aim to minimize the functional
\begin{equation}
\label{ls}
\left\|\begin{bmatrix}
\Re C_i & -\Im C_i\\
\Im C_i & \Re C_i 
\end{bmatrix}
\begin{pmatrix}
\Re {\bf{F}}_i\\
\Im {\bf{F}}_i
\end{pmatrix}
-\begin{pmatrix}
\Re \vb_i\\
\Im \vb_i
\end{pmatrix}
\right\|_2^2+\alpha\text{TV}_{\beta}(\Re{\bf{F}}_i,\Im {\bf{F}}_i),
\end{equation}
where 
\begin{equation}
\label{TV}
\text{TV}_{\beta}(\vx,\vy)=\sqrt{\paren{\sum_{i=1}^{np-1}(x_{i+1}-x_i)^2}+\paren{\sum_{i=1}^{np-1}(y_{i+1}-y_i)^2}+\beta^2}
\end{equation}
for $\vx=(x_1,\ldots,x_{np})$, and $\vy=(y_1,\ldots,y_{np})$. The regularization penalty \eqref{TV} is the smoothed 1-D TV regularizer of \cite{pet-mri} applied jointly to the real and imaginary parts of ${\bf{F}}_i$. The regularization parameter $\alpha>0$ controls the level of TV regularization. The $\beta>0$ term is included so that the gradient of $\text{TV}_{\beta}$ is defined at $(\vx,\vy)=({\bf{0}},{\bf{0}})$, and thus we can apply ideas from smooth optimization to solve \eqref{ls}. Specifically, to solve the objective in equation \eqref{ls}, we apply the L-BFGS-B code of \cite{byrd}. Finally, to generate the image, we solve \eqref{ls} for every line profile $i\in\{1,\ldots,m\}$. After which, the 2-D images $F_q$ are pieced together from the 1-D slices and the final image, $|F|$, is obtained from \eqref{|F|}. The reconstruction method detailed above will be denoted as Analytic Continuation (AC) for the remainder of this paper.

\subsection{Data sets}
\label{test_images}
\begin{figure}[!h]
\centering
\begin{subfigure}{0.24\textwidth}
\includegraphics[width=0.9\linewidth, height=3.2cm, keepaspectratio]{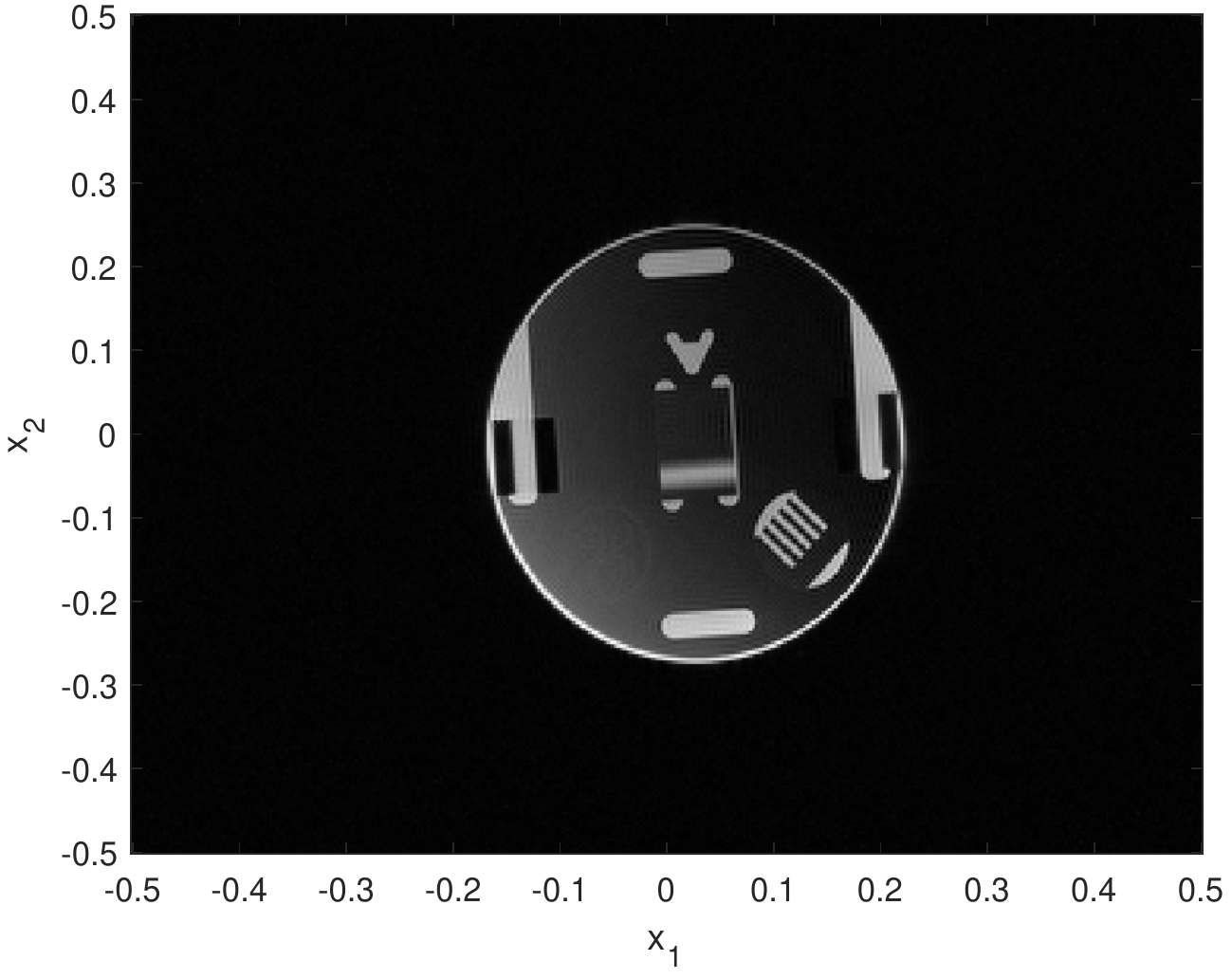}
\subcaption{phantom (16ch)}
\end{subfigure}
\begin{subfigure}{0.24\textwidth}
\includegraphics[width=0.9\linewidth, height=3.2cm, keepaspectratio]{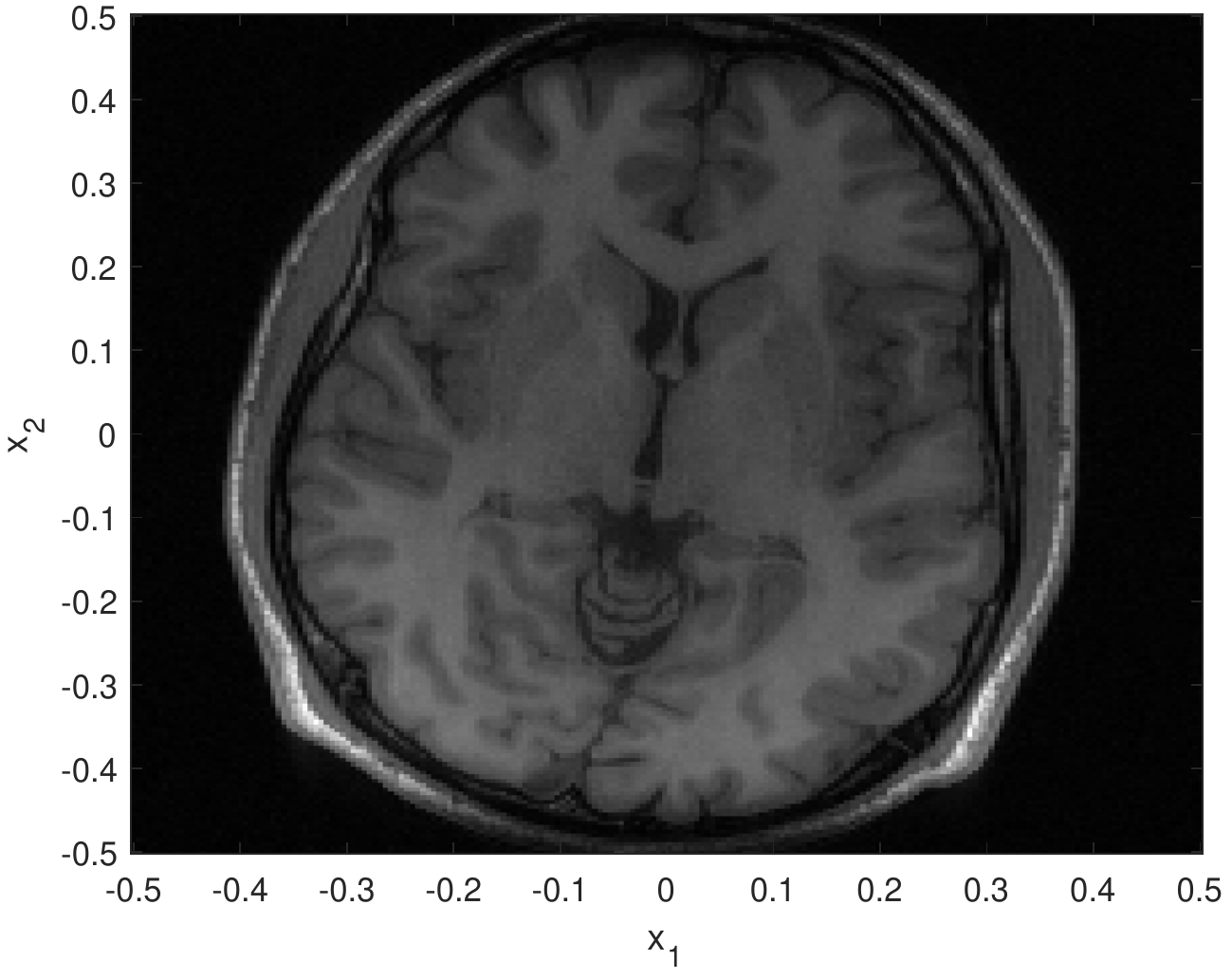}
\subcaption{brain (8ch)}
\end{subfigure}
\begin{subfigure}{0.24\textwidth}
\includegraphics[width=0.9\linewidth, height=3.2cm, keepaspectratio]{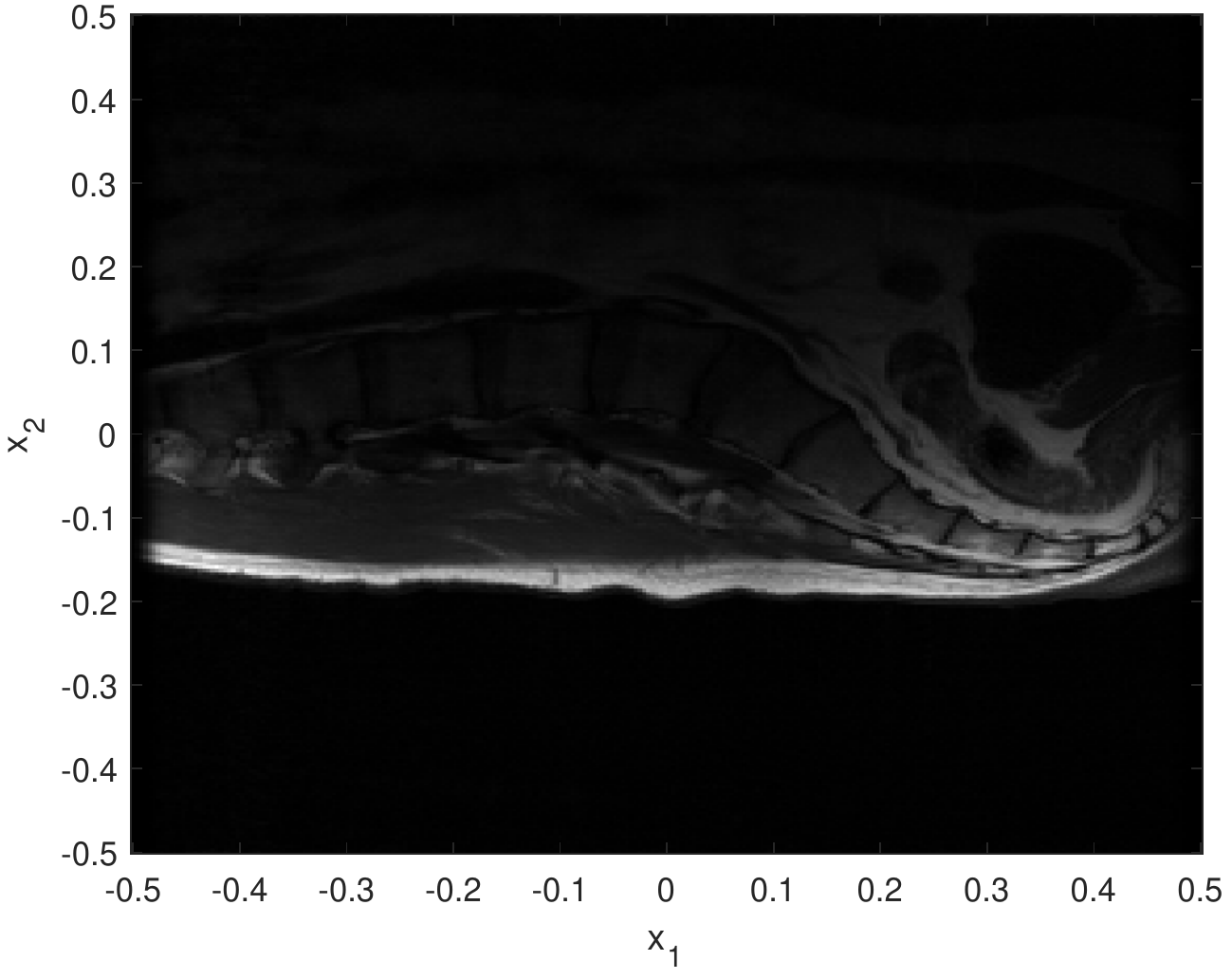}
\subcaption{spine (4ch)}
\end{subfigure}
\begin{subfigure}{0.24\textwidth}
\includegraphics[width=0.9\linewidth, height=3.2cm, keepaspectratio]{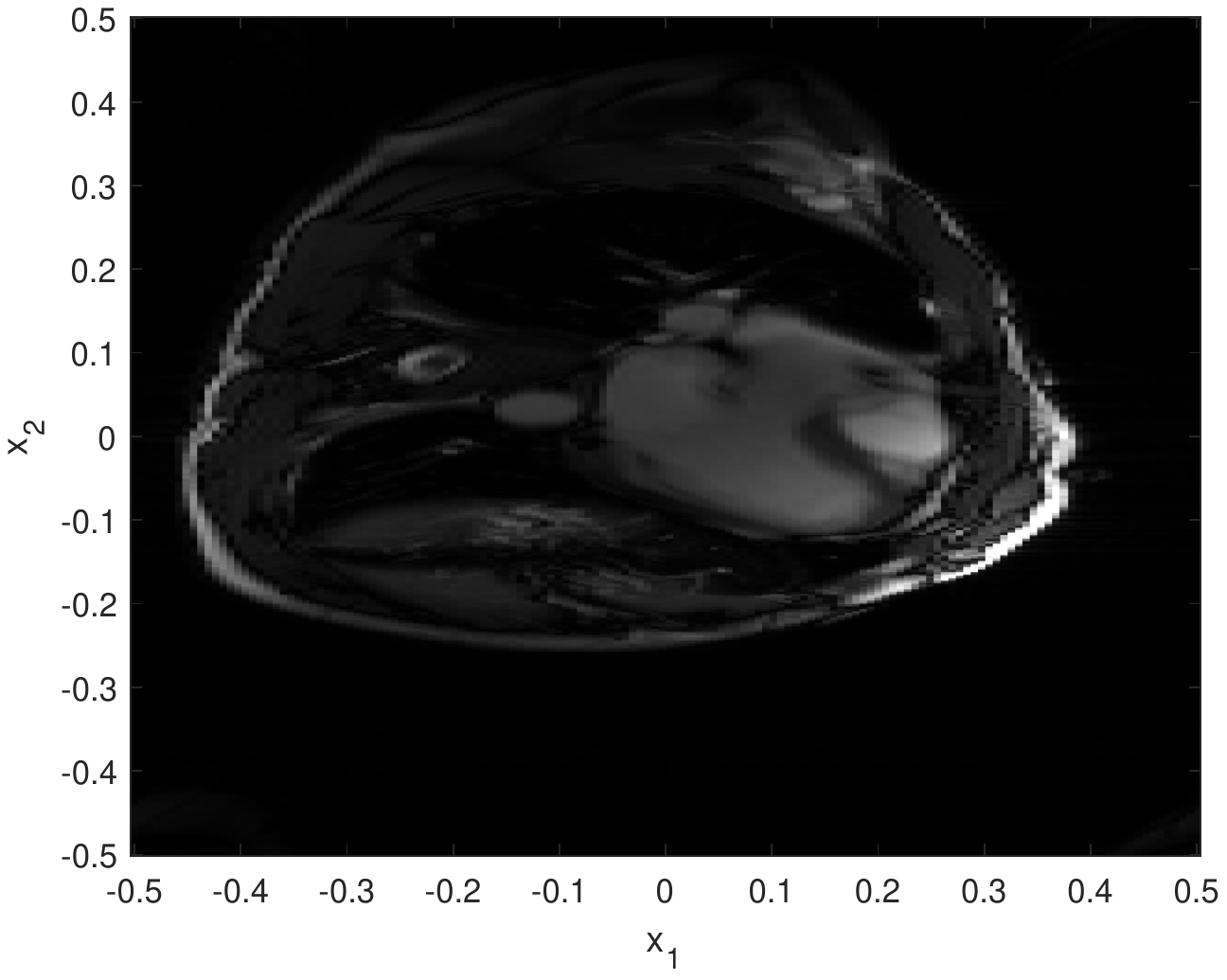}
\subcaption{cardiac (34ch)}
\end{subfigure}
\caption{SOS images.}
\label{images}
\end{figure}

Here we discuss the data sets that will be used to test our reconstruction method. We consider four images for testing. These  include, the brain image of figure \ref{bp}, and the phantom of figure \ref{rp}. We also consider the spine image data of \cite{PULSAR} and the CINE cardiac data from \url{https://ocmr.info/}. The SOS images corresponding to each data set are presented in figure \ref{images}.
The details of the data sets considered are given below:
\begin{itemize}
\item Brain image - the image resolution is $n=m=200$, and the number of coils is $K=8$. The data is downloaded from the BART toolbox \cite{BART}.
\item Real phantom - the image resolution is $n=256$, $m=340$, and the number of coils is $K=16$. The data is downloaded from \url{https://mr.usc.edu/download/data/}. To generate the data, a physical phantom was imaged on a 3T MRI scanner using a turbo spin-echo sequence. The data was acquired with a $220\text{mm}\times 292\text{mm}$ field of view on a $256\times 340$ Cartesian sampling grid. As per our agreement for use of this data, we acknowledge NSF support, specifically NSF grant CCF-1350563.
\item Spine image - the image resolution is $n=m=256$, and the number of coils is $K=4$. The data is downloaded from the PULSAR toolbox \cite{PULSAR}.
\item Cardiac image - the image resolution is $n=384$, $m=144$, and the number of coils is $K=34$. The data is downloaded from the CINE database \url{https://ocmr.info/} and read using the code of \cite{OCMR}.
\end{itemize}
For further information on the data sets listed see the associated references. The data sets are chosen to reflect a wide range of image resolution, number of coils, and practical application. For example, the real phantom is an example image with sharp edges (high-frequency Fourier components), and the brain is a classic example from {\it{in vivo}} medical MRI.

\subsection{Methods for comparison}
\label{comp_methods}
Here we discuss the methods from the literature that AC will be compared against. AC is variation of  ESPIRiT whereby the image is reconstructed slice-by-slice on 1-D lines. With this in mind, we choose to compare against three other ESPIRiT variations from the literature, namely CG ESPIRiT \cite{ESP,SENSE}, $L^1$ ESPIRiT \cite{feng,ESP,ESP1}, and TV ESPIRiT \cite{feng}. The aim of the methods listed above is to minimize the functional
\begin{equation}
\|\mathcal{F}_TS{\bf{F}}-\vb\|_2^2+\alpha\mathcal{R}({\bf{F}}),
\end{equation}
where $\mathcal{F}_T$ is the truncated Fourier operator, $S$ are the sensitivity weightings, ${\bf{F}}$ is the unknown image, and $\vb$ is the data. The regularization penalty $\mathcal{R}$ depends on the method as described in the list below:
\begin{itemize}
\item CG ESPIRiT (denoted ESP for short) - $\mathcal{R}({\bf{F}})=\|{\bf{F}}\|^2_2$, or ESPIRiT with Tikhonov regularization.
\item $L^1$ ESPIRiT (denoted $L^1$ ESP) - $\mathcal{R}({\bf{F}})=\|\Psi{\bf{F}}\|_1$, where $\Psi$ is a pre-specified sparsifying transform (e.g., wavelet). We set $\Psi$ as a translation invariant Daubechies wavelet, as is done in the BART examples \cite{BART}.
\item TV ESPIRiT (denoted TV ESP)  - $\mathcal{R}({\bf{F}})=\sqrt{\|\nabla \Re{\bf{F}}\|^2_2+\|\nabla \Im{\bf{F}}\|^2_2+\beta^2}$. This is the same smoothed TV penalty of \eqref{TV}, although generalized to 2-D images. We use the same smoothed TV idea for AC and TV ESP, for fairness.
\end{itemize}
To implement CG ESPIRiT and $L^1$ ESPIRiT we use the Matlab code supplied as part of the BART toolbox \cite{BART}. To implement TV ESPIRiT, we use the algorithm of \cite{feng}. Specifically, using the notation of \cite[equation (1)]{feng}, we set the wavelet sparsity parameter $\mu=0$, and replace the conventional TV penalty with smoothed TV as defined above. 

\subsection{Performance metrics}
Here we define the metrics that will be used to measure performance. Let ${\bf{F}}$ and ${\bf{F}}_{\epsilon}$ denote vectors of ground truth and estimated image pixel values within a Region Of Interest (ROI), which encapsulates the nonzero region (support) of the image. Then we define the relative least squares error
\begin{equation}
\epsilon=\frac{\|{\bf{F}}-{\bf{F}}_{\epsilon}\|_2}{\|{\bf{F}}\|_2}.
\end{equation}
We also consider the structural similarity index, commonly used to evaluate reconstruction quality in MRI \cite{SSIM}. Let $X$ and $Y$ respectively denote an $N\times N$ neighborhood of a ground truth and estimate image, and let $\vx$ and $\vy$ denote their corresponding vectorized forms. Then, the structural similarity shared by $X$ and $Y$ is defined as
\begin{equation}
\label{ssim}
\text{SSIM}(X,Y)=\frac{(2\mu_{\vx}\mu_{\vy}+c_1)(2\sigma_{\vx\vy}+c_2)}{(\mu_{\vx}^2+\mu_{\vy}^2+c_1)(\sigma_{\vx}^2+\sigma_{\vy}^2+c_2)},
\end{equation}
where $\mu_{\vx}$ denotes the mean value of $\vx$, $\sigma_{\vx}$ denotes the standard deviation, and $\sigma_{\vx\vy}$ denotes the covariance of $\vx$ and $\vy$. $c_1$ and $c_2$ are small-value constants included to stabilize the division. Specifically we set $c_1=.0001$, and $c_2=.0009$. SSIM varies between $0$ and $1$, with 1 indicating perfect structural similarity and 0, no similarity. We calculate SSIM on every $3\times 3$ window in the ROI and take the average to evaluate structural similarity. We denote the mean structural similarity on the ROI by $\text{SSIM}_{\mu}$. As a ``proper" ground truth image is not available for the data sets considered, we use the SOS image calculated from the complete data to calculate the performance metrics. Note, the SOS image is not a ground truth since it contains measurement noise.

\subsection{Hyperparameter selection}
\label{hyperparam}
Here we discuss selection of hyperparameters. The sensitivity maps are calculated using ESPIRiT, using the hyperparameters (e.g., kernel window size, eigenvalue threshold) specified in the BART examples \cite{BART}. The number of ESPIRiT maps is set at $p=2$ throughout. The parameters used to calculate the sensitivity maps are fixed throughout this paper. The splitting parameter (as defined in \cite{feng}) for TV ESPIRiT and $L^1$ ESPIRiT is set at 0.4. For TV ESPIRiT and AC, $\beta=0.01$ is kept fixed throughout. The smoothing parameter, $\alpha$, is chosen to give the best performance in terms of $\text{SSIM}_{\mu}$, for all methods compared against. 

We emphasize that the hyperparameters discussed above (e.g., $\alpha,\beta,p$) were chosen heuristically, for each method considered. Readers of this paper may wish to consider hyperparameter selection methods, such as the discrepancy principle \cite{DP1}.


\section{Results}
In this section, we test our method on the images introduced in section \ref{test_images}, when the locations of the missing $\textbf{k}$-space lines are selected at random from a uniform distribution. 
Random subsampling is used to simulate motion corruption, as is done also in \cite{motion1,motion2}. For example, random subsampling can be used to simulate artifacts due to rotation and translation of the head in brain MRI \cite[figure 2 (a)]{motion2}.
In all examples conducted, we retain 32 central lines of $\textbf{k}$-space, for calibration. Retention of the central $\textbf{k}$-space lines is also done in the simulations of \cite{motion1}.

\begin{figure}[!h]
\centering
\begin{subfigure}{0.24\textwidth}
\includegraphics[width=0.9\linewidth, height=3.2cm, keepaspectratio]{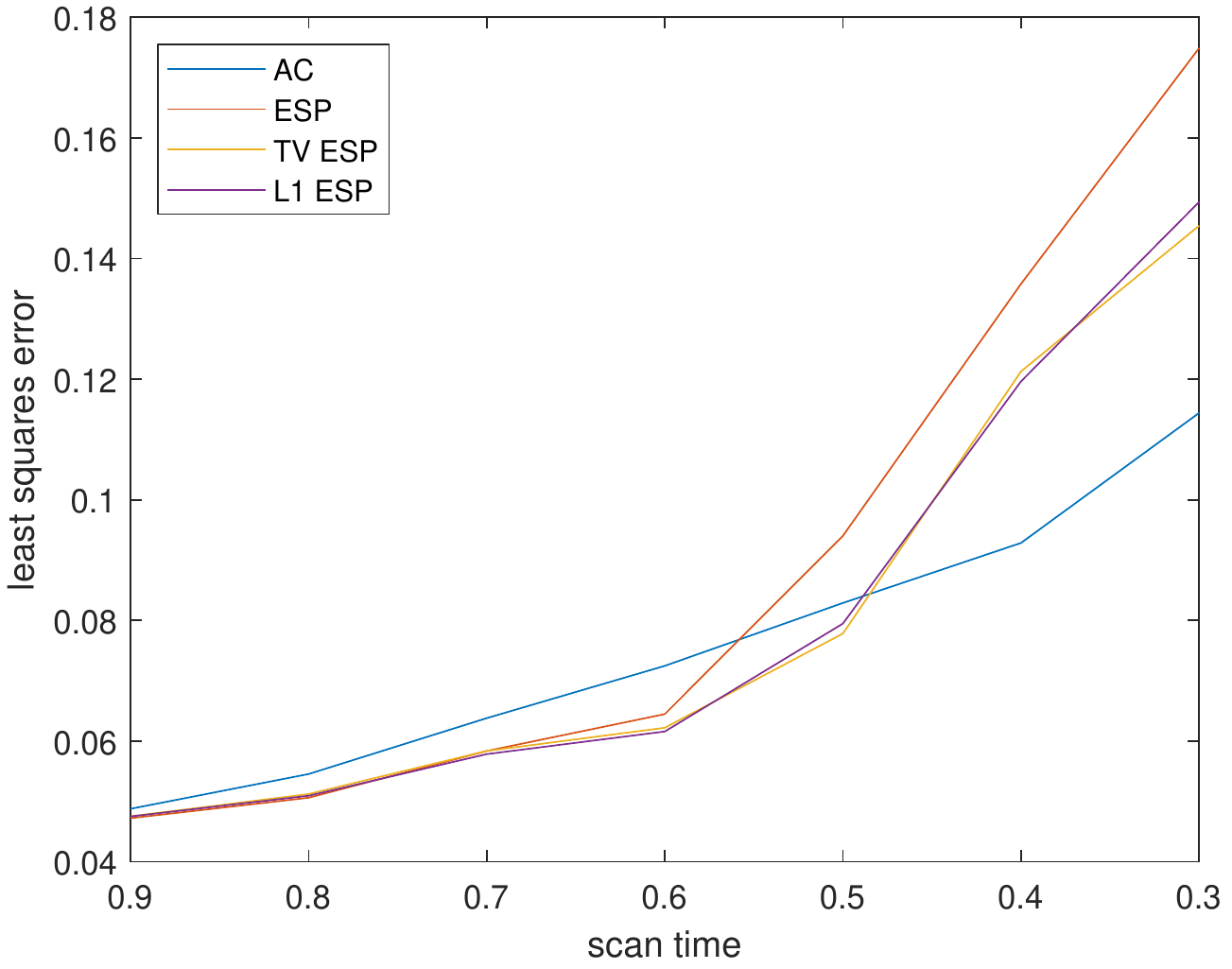}
\end{subfigure}
\begin{subfigure}{0.24\textwidth}
\includegraphics[width=0.9\linewidth, height=3.2cm, keepaspectratio]{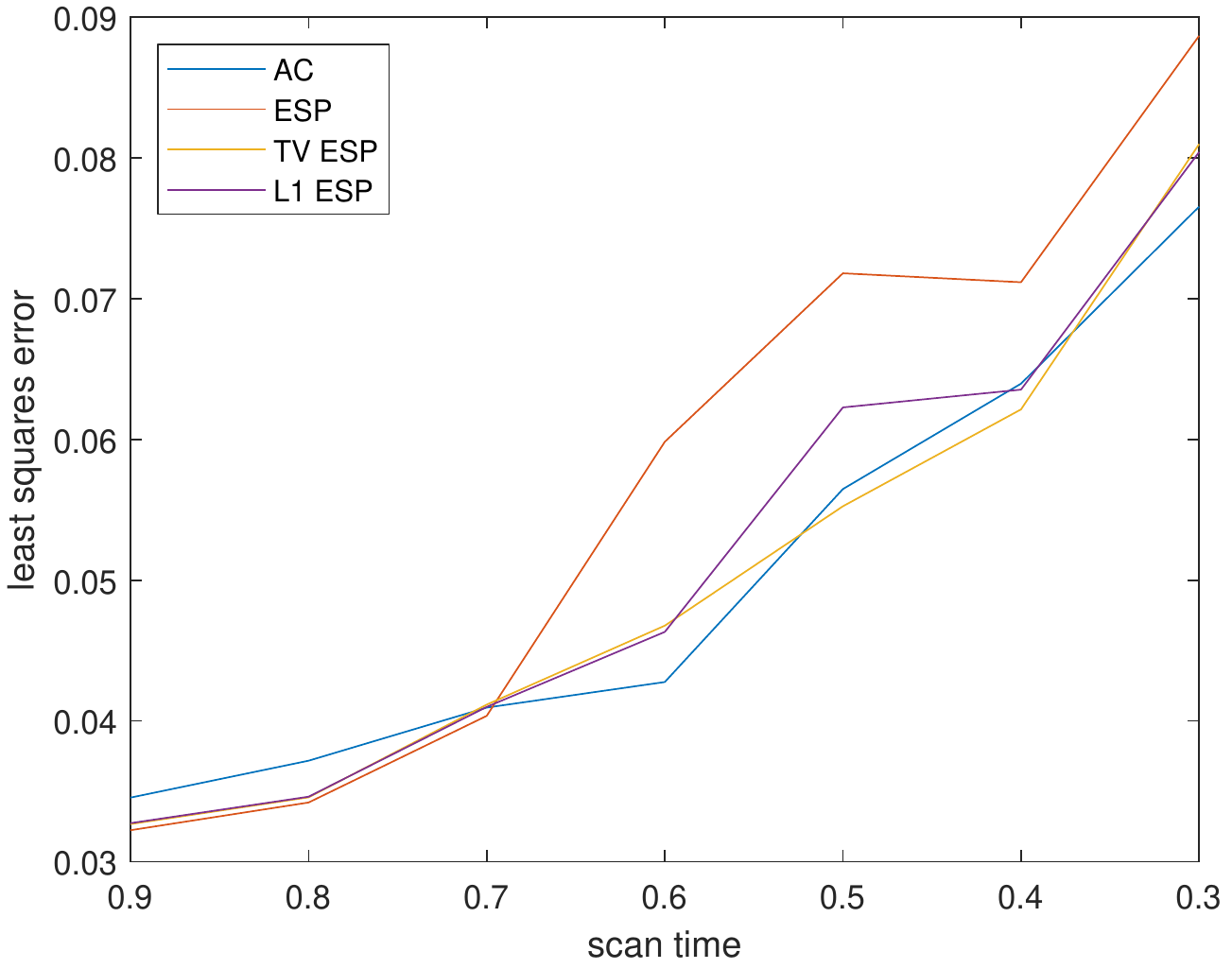}
\end{subfigure}
\begin{subfigure}{0.24\textwidth}
\includegraphics[width=0.9\linewidth, height=3.2cm, keepaspectratio]{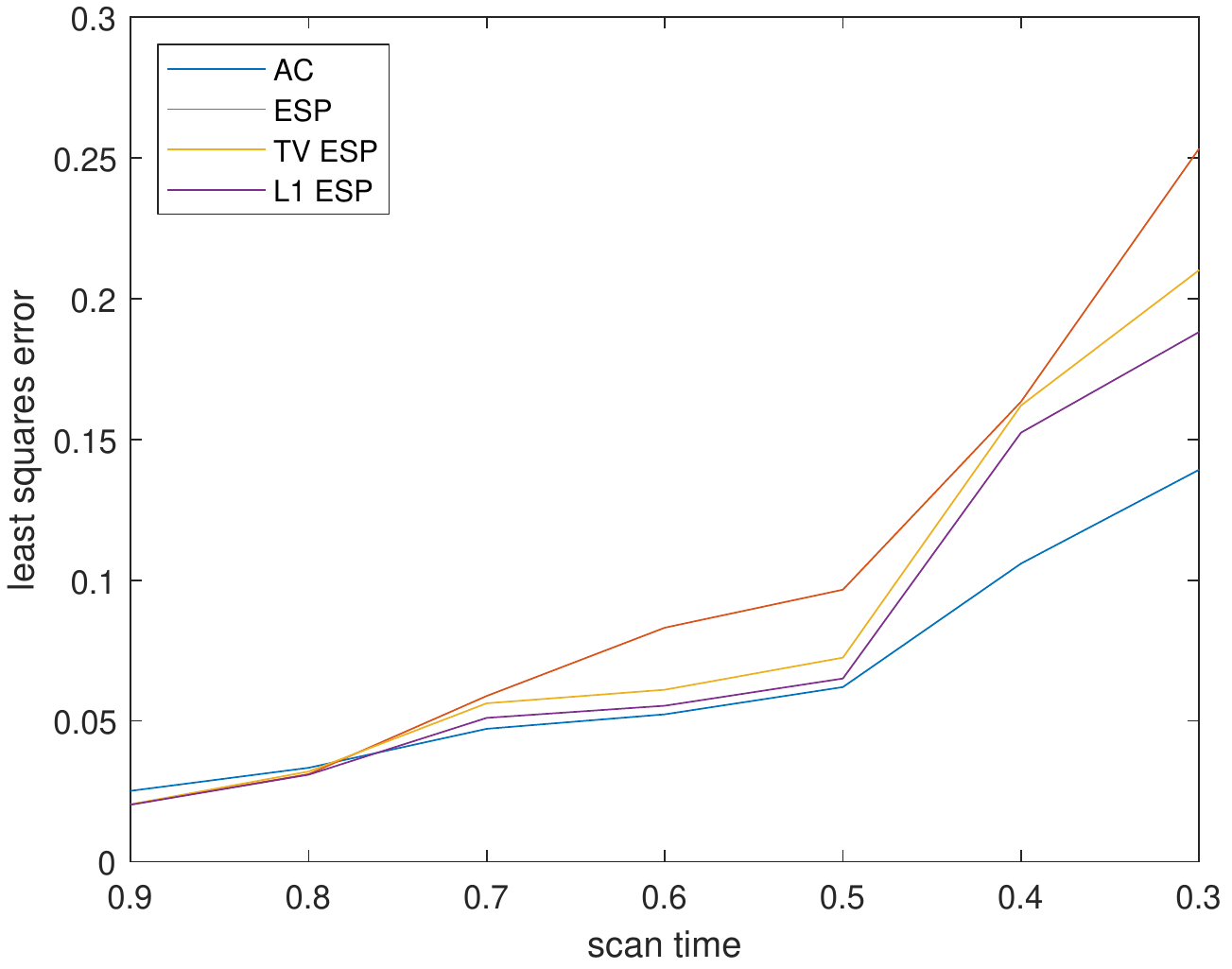}
\end{subfigure}
\begin{subfigure}{0.24\textwidth}
\includegraphics[width=0.9\linewidth, height=3.2cm, keepaspectratio]{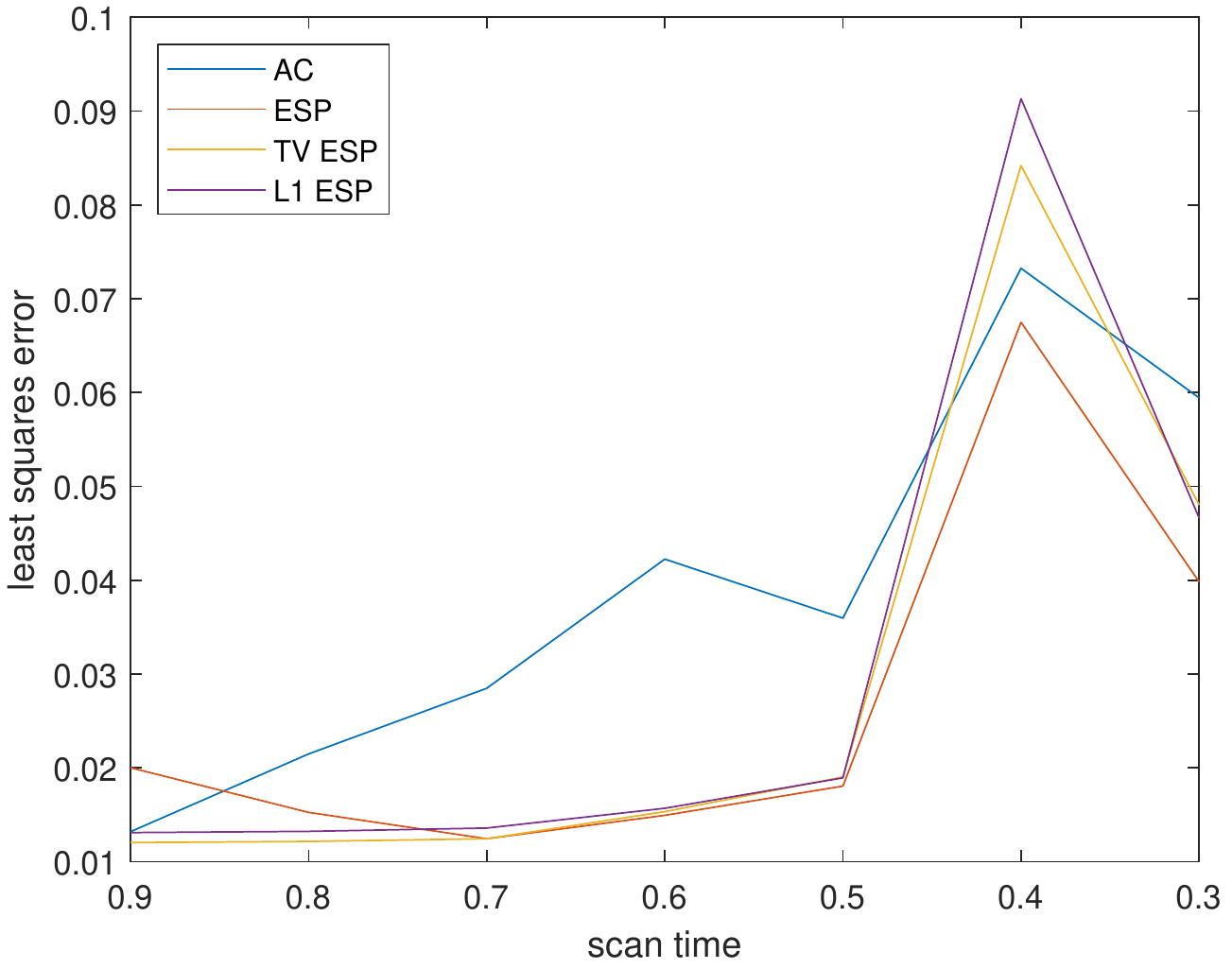}
\end{subfigure}
\begin{subfigure}{0.24\textwidth}
\includegraphics[width=0.9\linewidth, height=3.2cm, keepaspectratio]{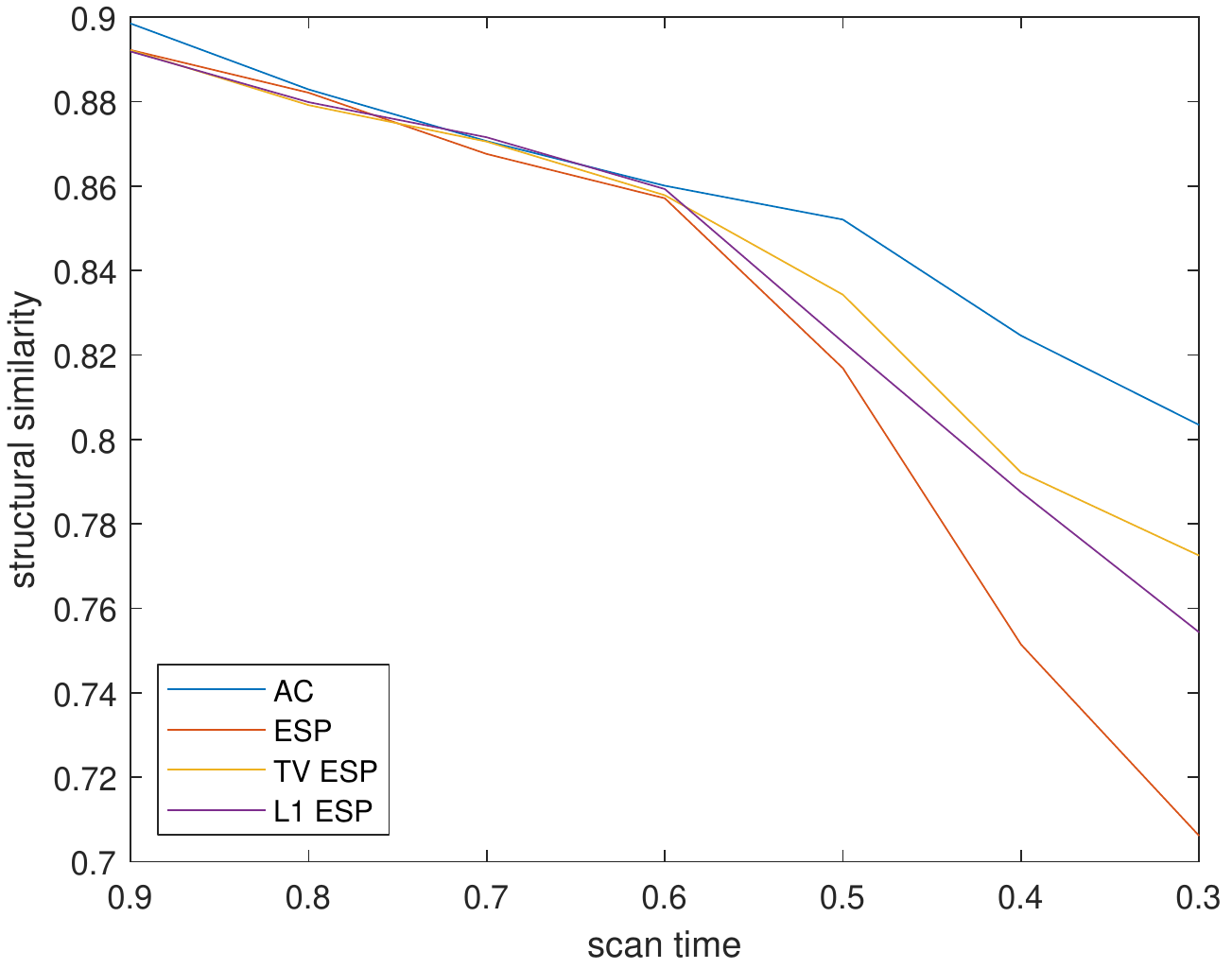}
\subcaption*{phantom (16ch)}
\end{subfigure}
\begin{subfigure}{0.24\textwidth}
\includegraphics[width=0.9\linewidth, height=3.2cm, keepaspectratio]{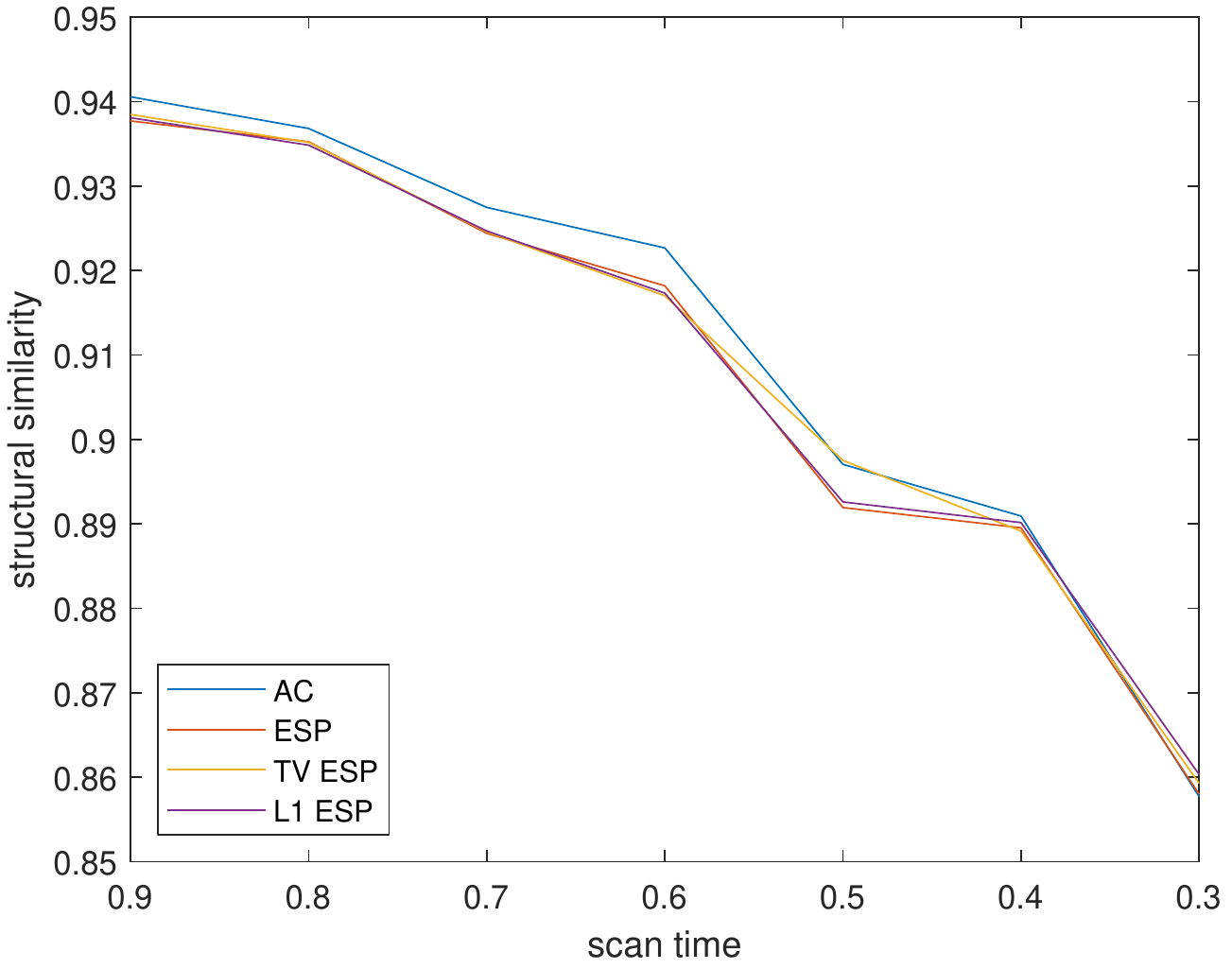}
\subcaption*{brain (8ch)}
\end{subfigure}
\begin{subfigure}{0.24\textwidth}
\includegraphics[width=0.9\linewidth, height=3.2cm, keepaspectratio]{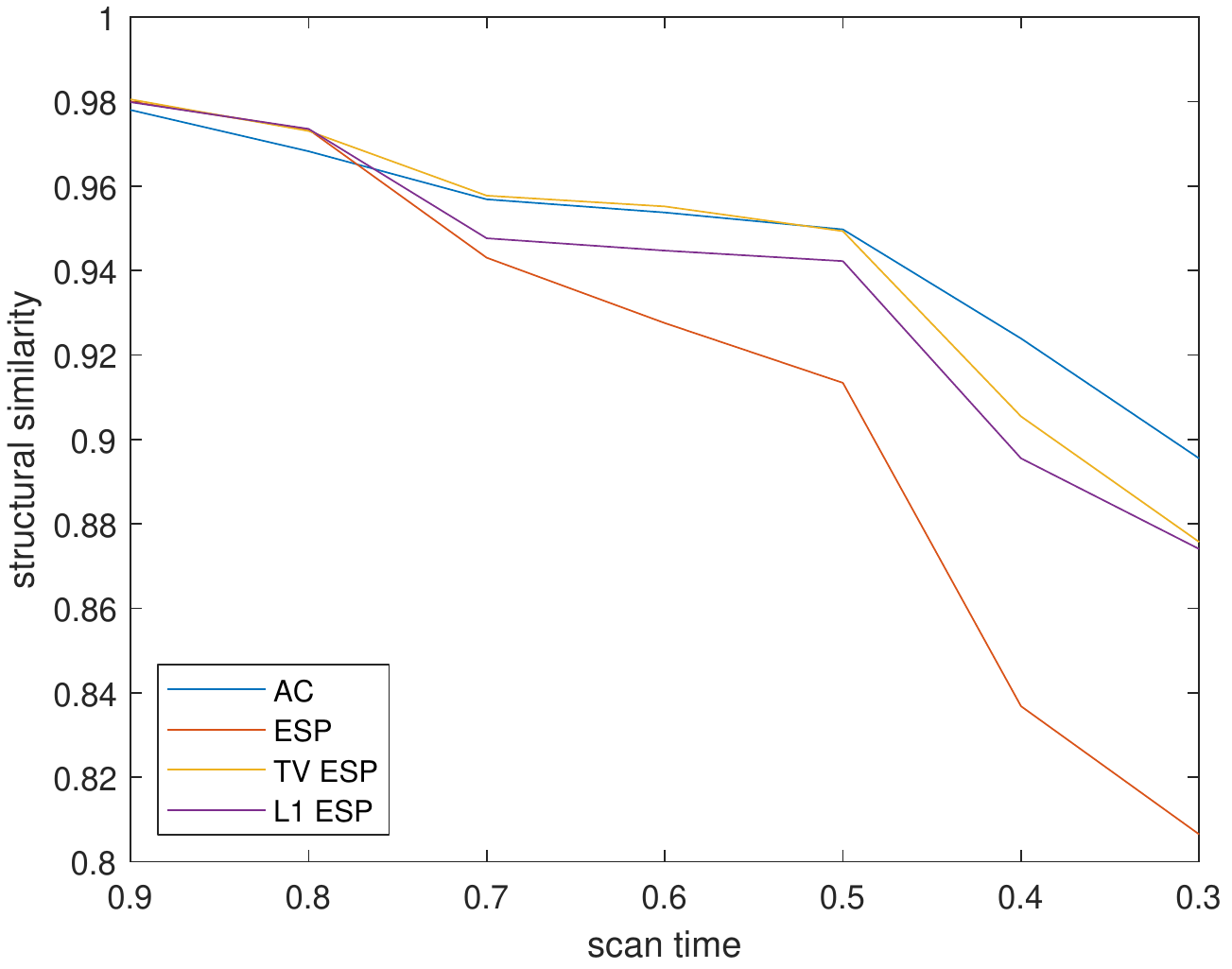}
\subcaption*{spine (4ch)}
\end{subfigure}
\begin{subfigure}{0.24\textwidth}
\includegraphics[width=0.9\linewidth, height=3.2cm, keepaspectratio]{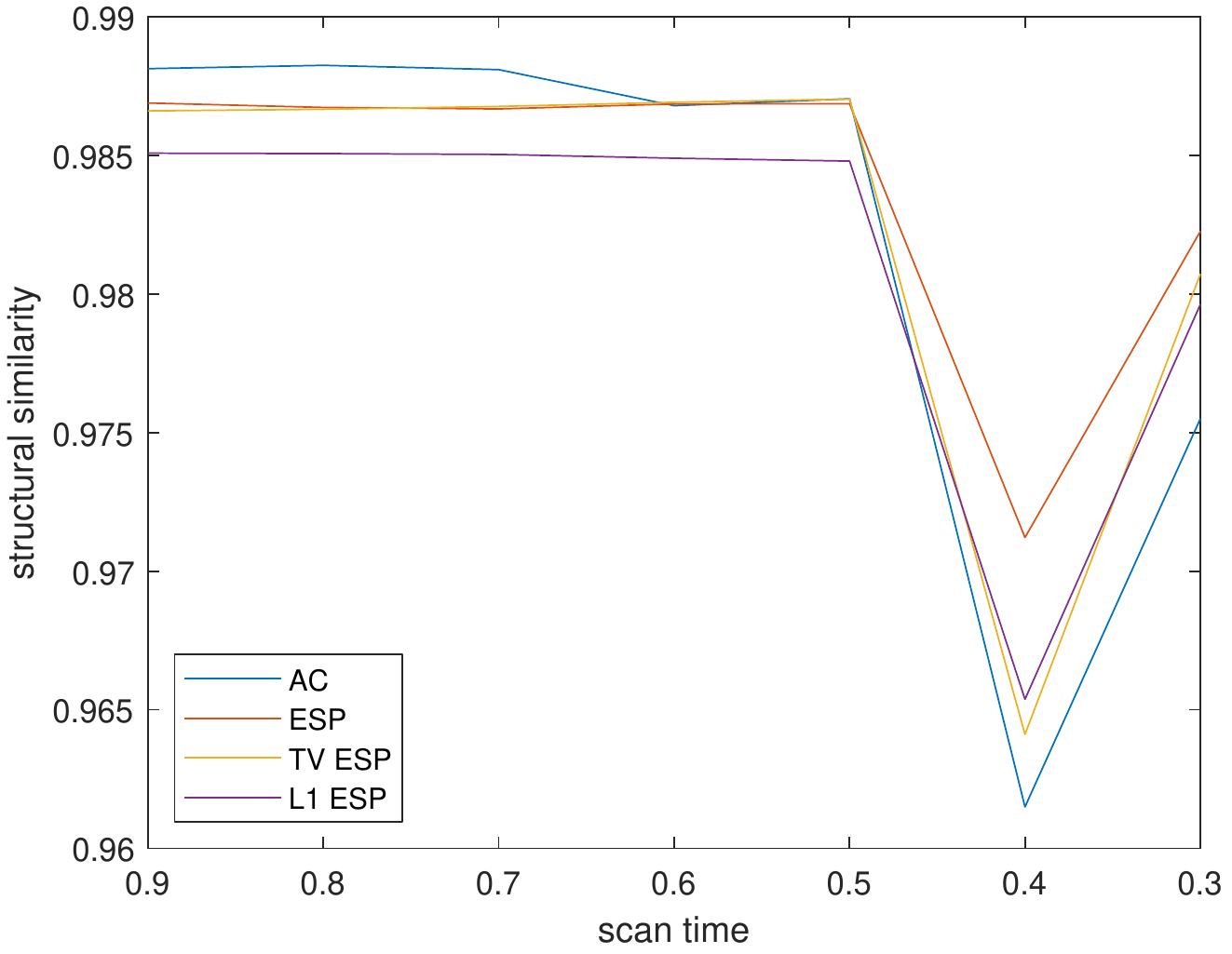}
\subcaption*{cardiac (34ch)}
\end{subfigure}
\caption{Errors against scan time curves for each method and image considered. Top row - least squares error. Bottom row - structural similarity.}
\label{error curves}
\end{figure}

In figure \ref{error curves}, we plot the least squares error ($\epsilon$) and mean structural similarity ($\text{SSIM}_{\mu}$) against the scan time for each method considered. The scan time is as defined in section \ref{subsamp}.
When the scan time exceeds $50\%$, AC is shown to offer competitive performance when compared to ESP, $L^1$ ESP, and TV ESP, in terms of structural similarity, for all images considered (see the bottom row of figure \ref{error curves}). For the real phantom and spine images, AC offers greater structural similarity and reduced error (in the order of a few percentage points) in comparison to its competitors when the scan percentage is $50\%$ or less. In general, the structural similarity offered by AC is greater than or equal to that of ESP, TV ESP and $L^1$ ESP, across all scan times and images considered, except in the case of the cardiac image with $\text{scan time}<50\%$ where the differences in structural similarity are $\approx 0.5\%$ and all methods perform almost equally as well in terms of structural similarity.
\begin{figure}[!h]
\centering
\begin{subfigure}{0.24\textwidth}
\includegraphics[width=0.9\linewidth, height=3.2cm, keepaspectratio]{phantom_GT}
\end{subfigure}
\hspace{-0.5cm}
\begin{subfigure}{0.24\textwidth}
\includegraphics[width=0.9\linewidth, height=3.2cm, keepaspectratio]{brain_GT}
\end{subfigure}
\hspace{-0.5cm}
\begin{subfigure}{0.24\textwidth}
\includegraphics[width=0.9\linewidth, height=3.2cm, keepaspectratio]{spine_GT}
\end{subfigure}
\hspace{-0.5cm}
\begin{subfigure}{0.24\textwidth}
\includegraphics[width=0.9\linewidth, height=3.2cm, keepaspectratio]{cardiac_GT}
\end{subfigure}
\begin{subfigure}{0.24\textwidth}
\includegraphics[width=0.9\linewidth, height=3.2cm, keepaspectratio]{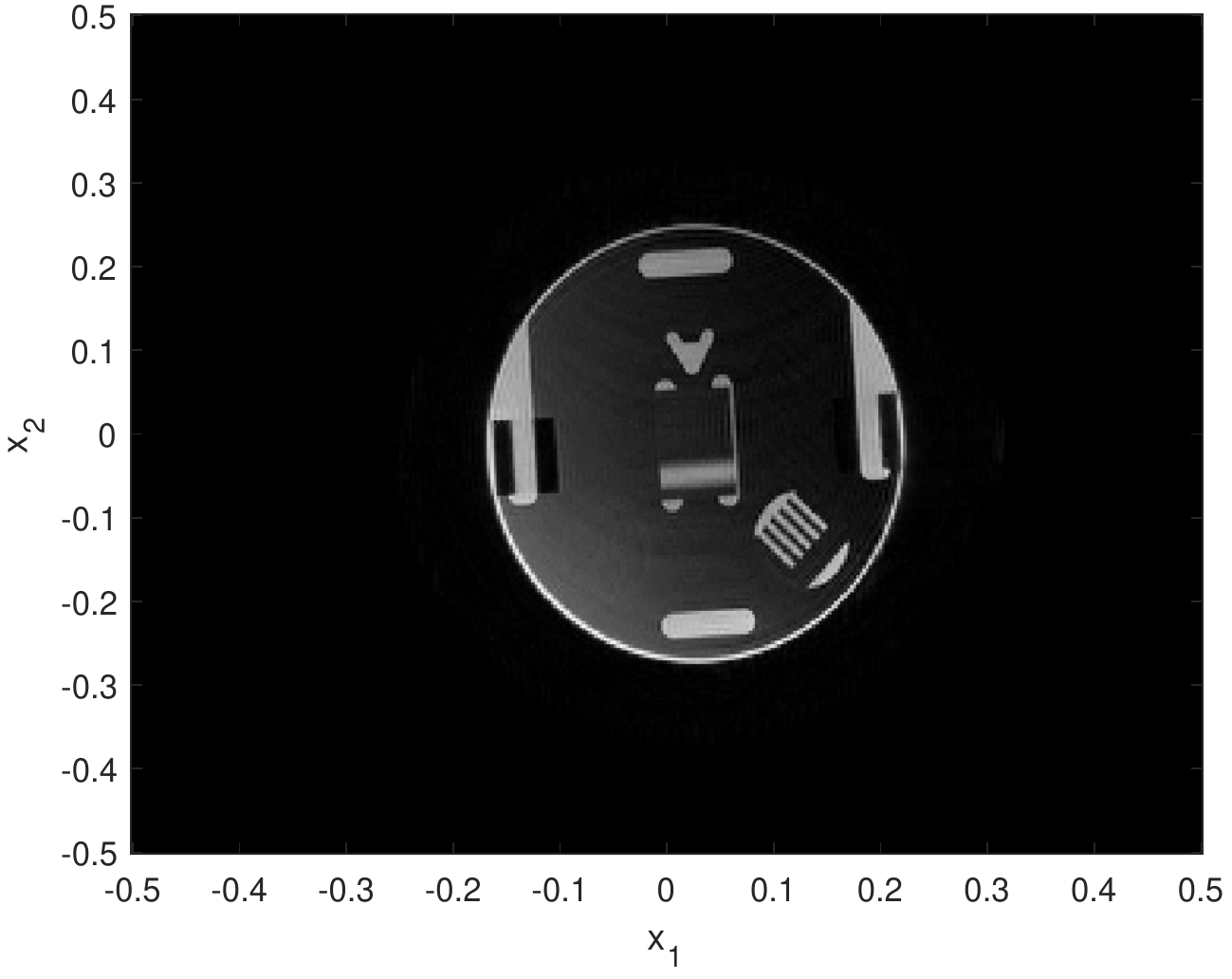}
\end{subfigure}
\hspace{-0.5cm}
\begin{subfigure}{0.24\textwidth}
\includegraphics[width=0.9\linewidth, height=3.2cm, keepaspectratio]{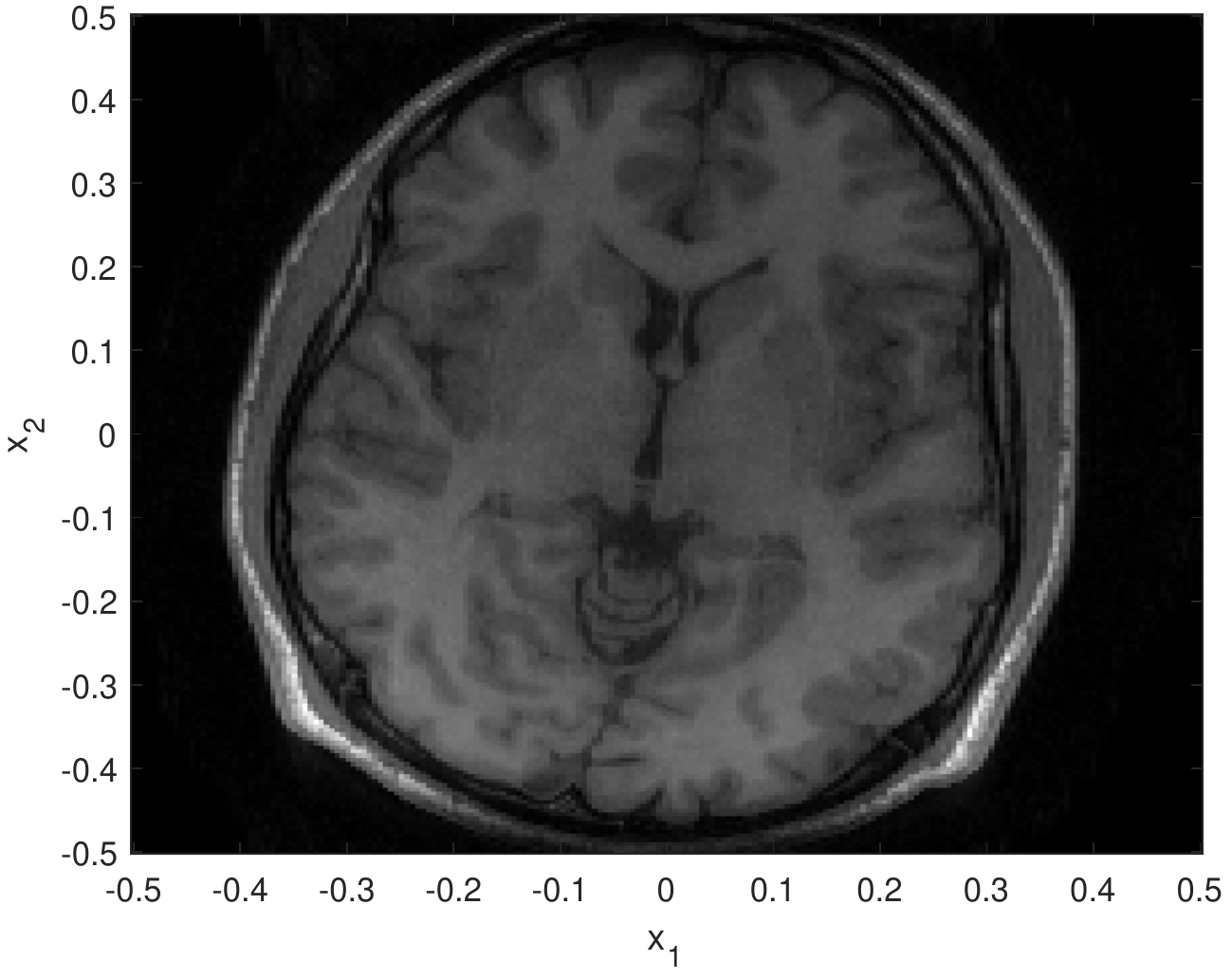}
\end{subfigure}
\hspace{-0.5cm}
\begin{subfigure}{0.24\textwidth}
\includegraphics[width=0.9\linewidth, height=3.2cm, keepaspectratio]{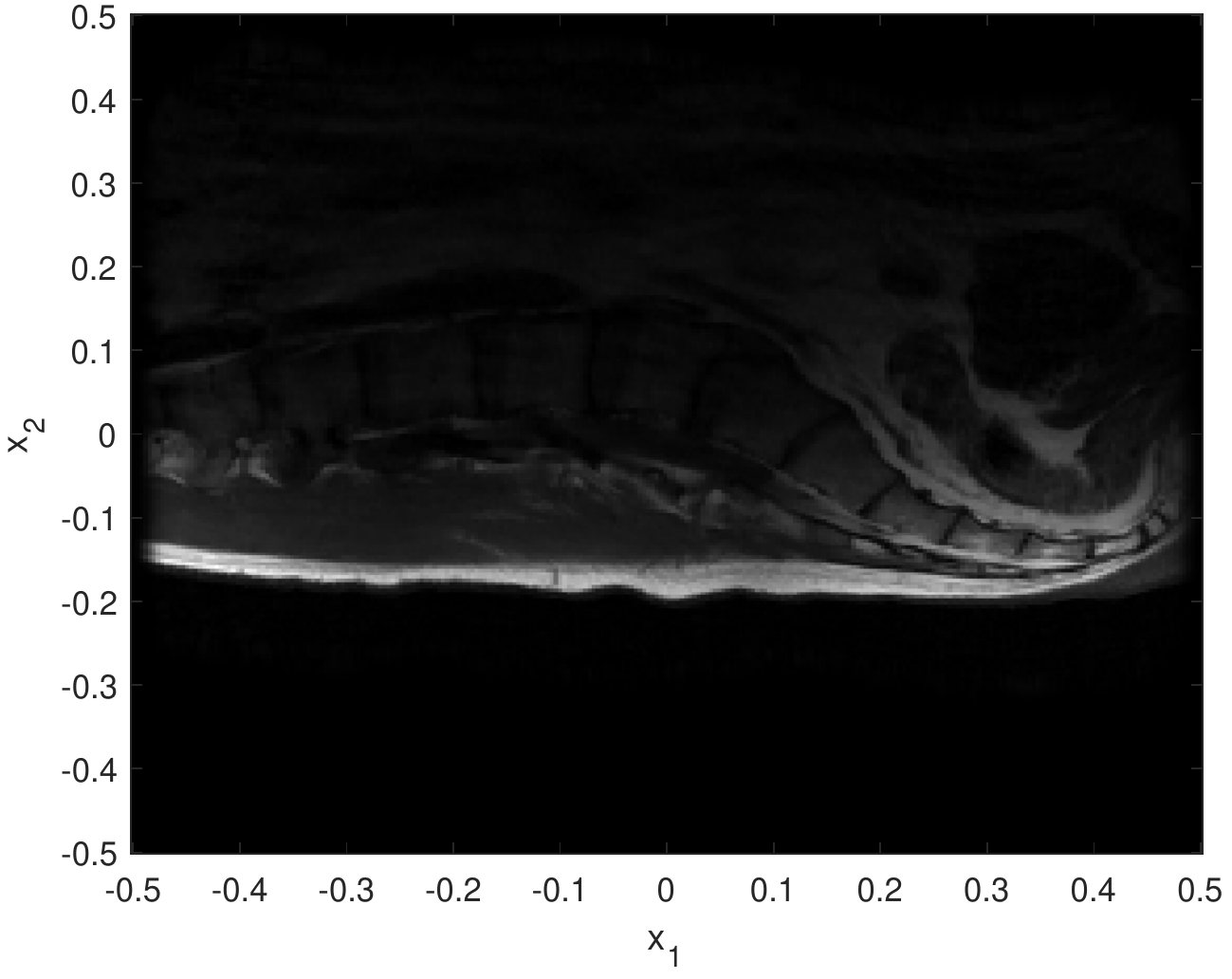}
\end{subfigure}
\hspace{-0.5cm}
\begin{subfigure}{0.24\textwidth}
\includegraphics[width=0.9\linewidth, height=3.2cm, keepaspectratio]{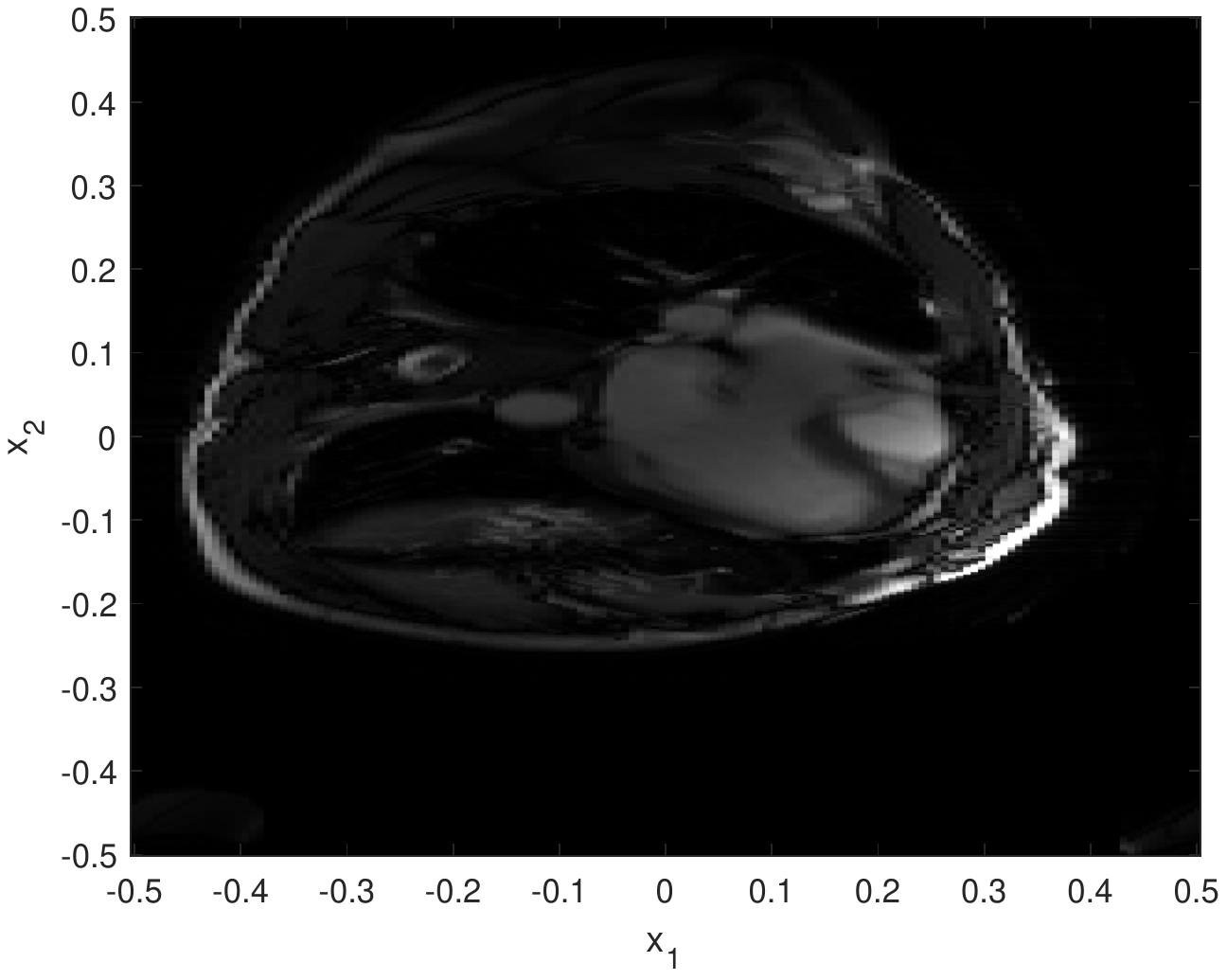}
\end{subfigure}
\begin{subfigure}{0.24\textwidth}
\includegraphics[width=0.9\linewidth, height=3.2cm, keepaspectratio]{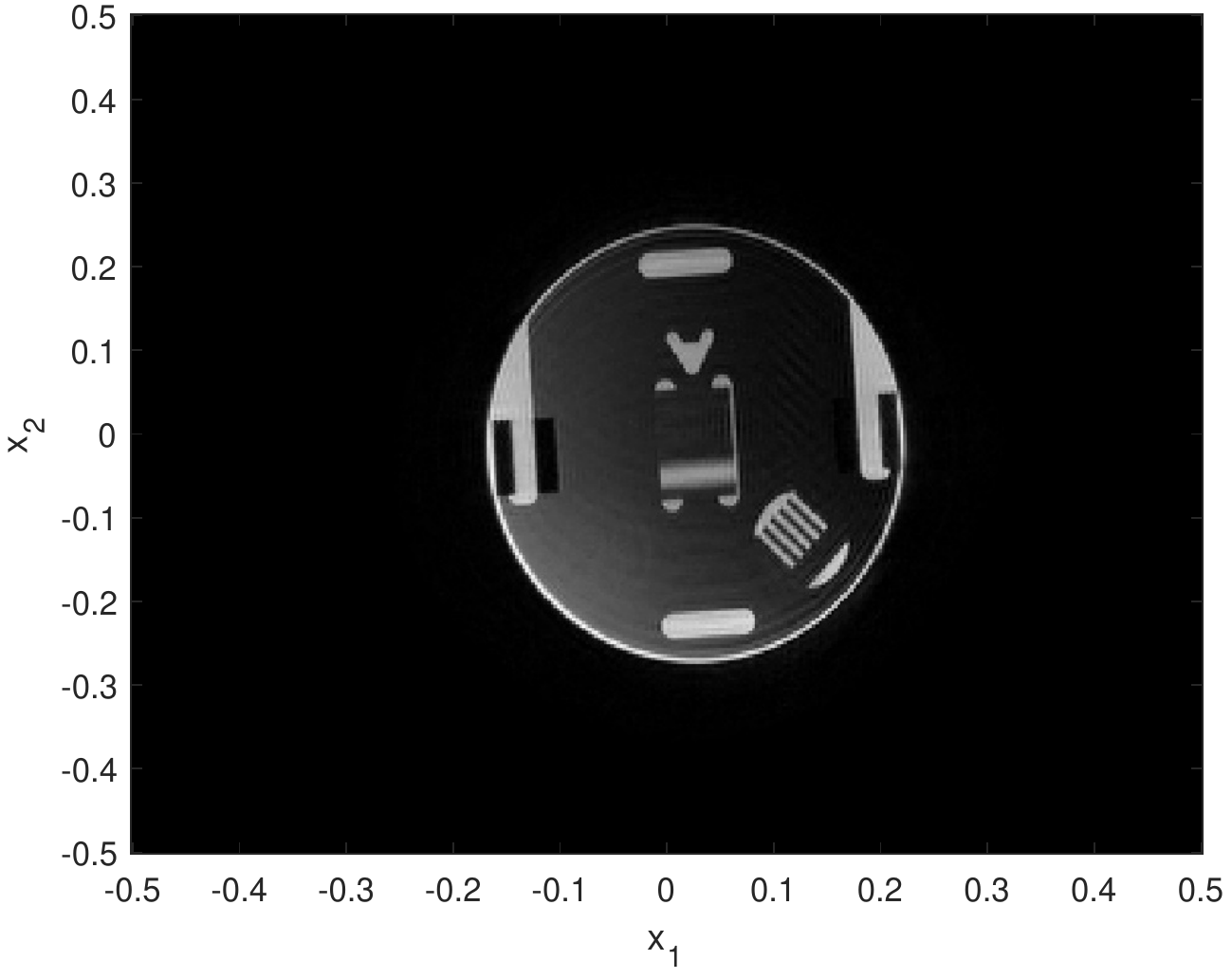}
\end{subfigure}
\hspace{-0.5cm}
\begin{subfigure}{0.24\textwidth}
\includegraphics[width=0.9\linewidth, height=3.2cm, keepaspectratio]{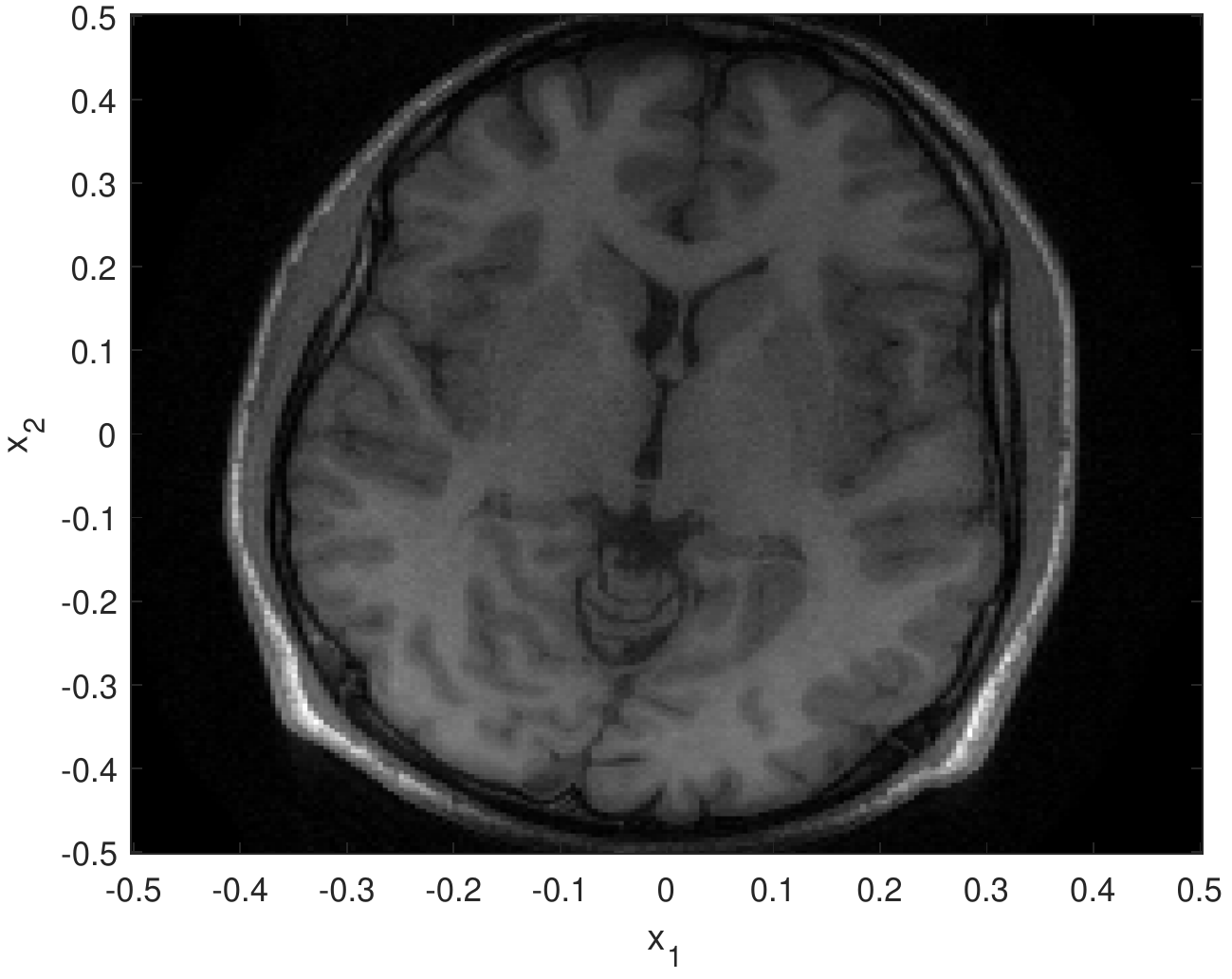}
\end{subfigure}
\hspace{-0.5cm}
\begin{subfigure}{0.24\textwidth}
\includegraphics[width=0.9\linewidth, height=3.2cm, keepaspectratio]{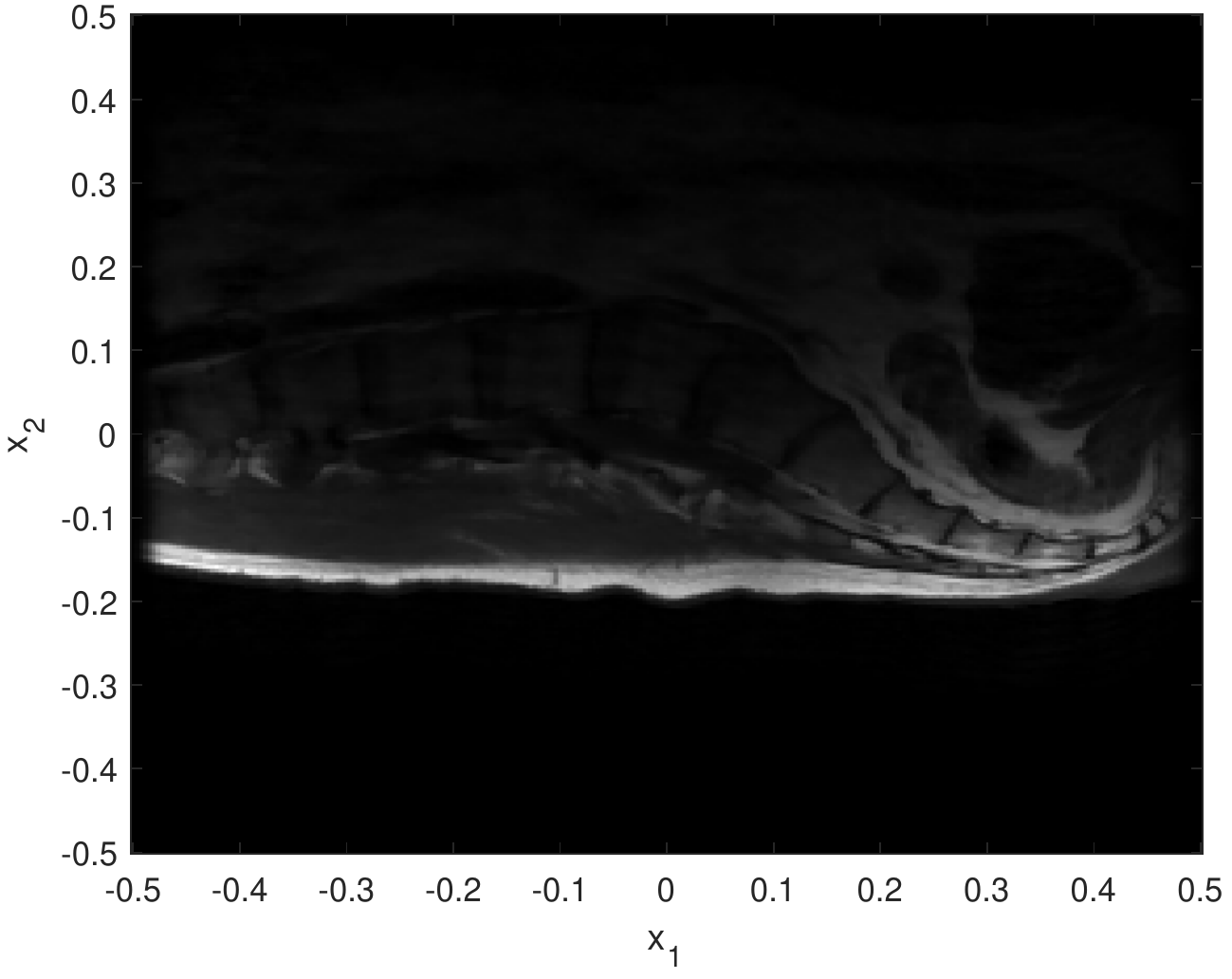}
\end{subfigure}
\hspace{-0.5cm}
\begin{subfigure}{0.24\textwidth}
\includegraphics[width=0.9\linewidth, height=3.2cm, keepaspectratio]{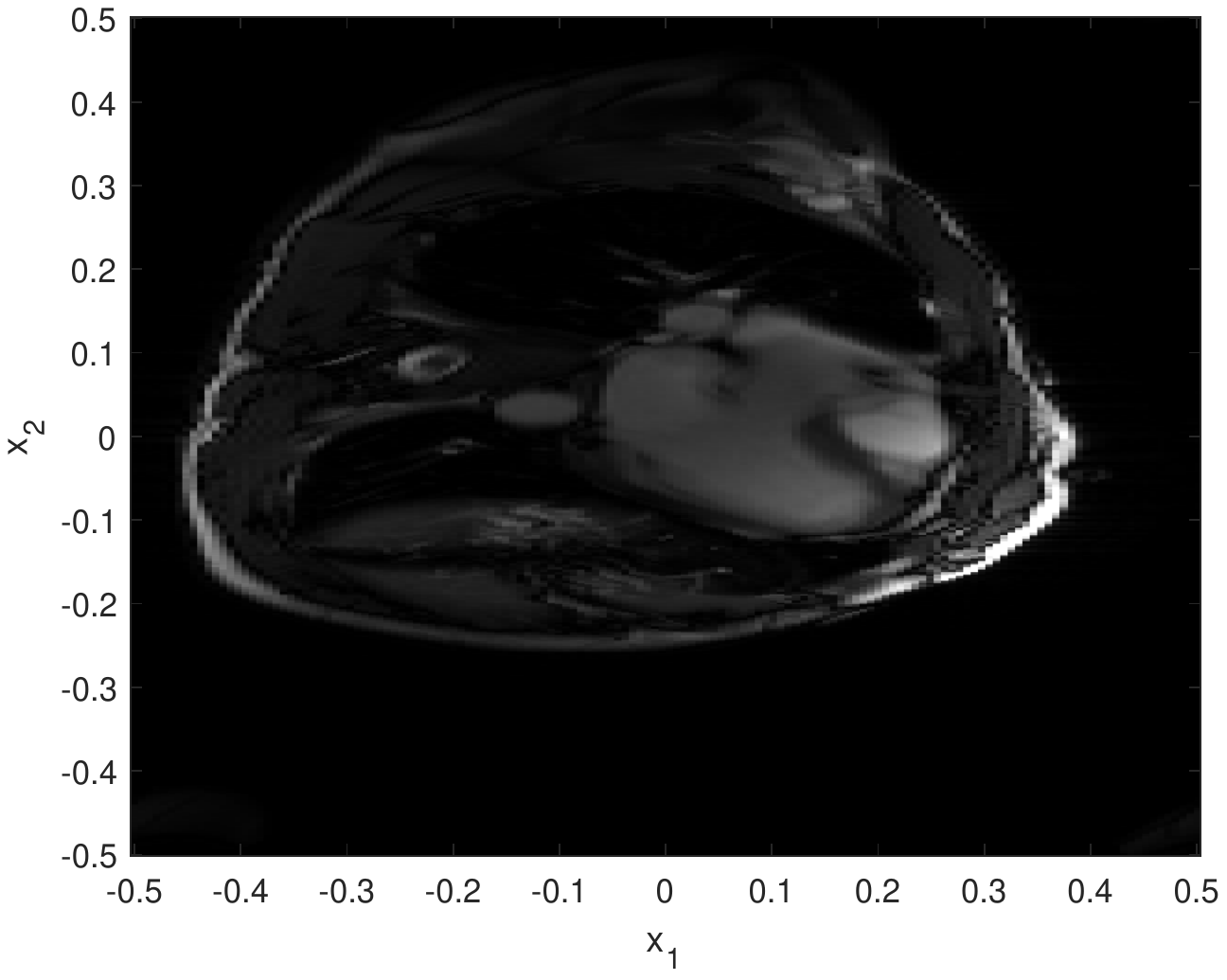}
\end{subfigure}
\begin{subfigure}{0.24\textwidth}
\includegraphics[width=0.9\linewidth, height=3.2cm, keepaspectratio]{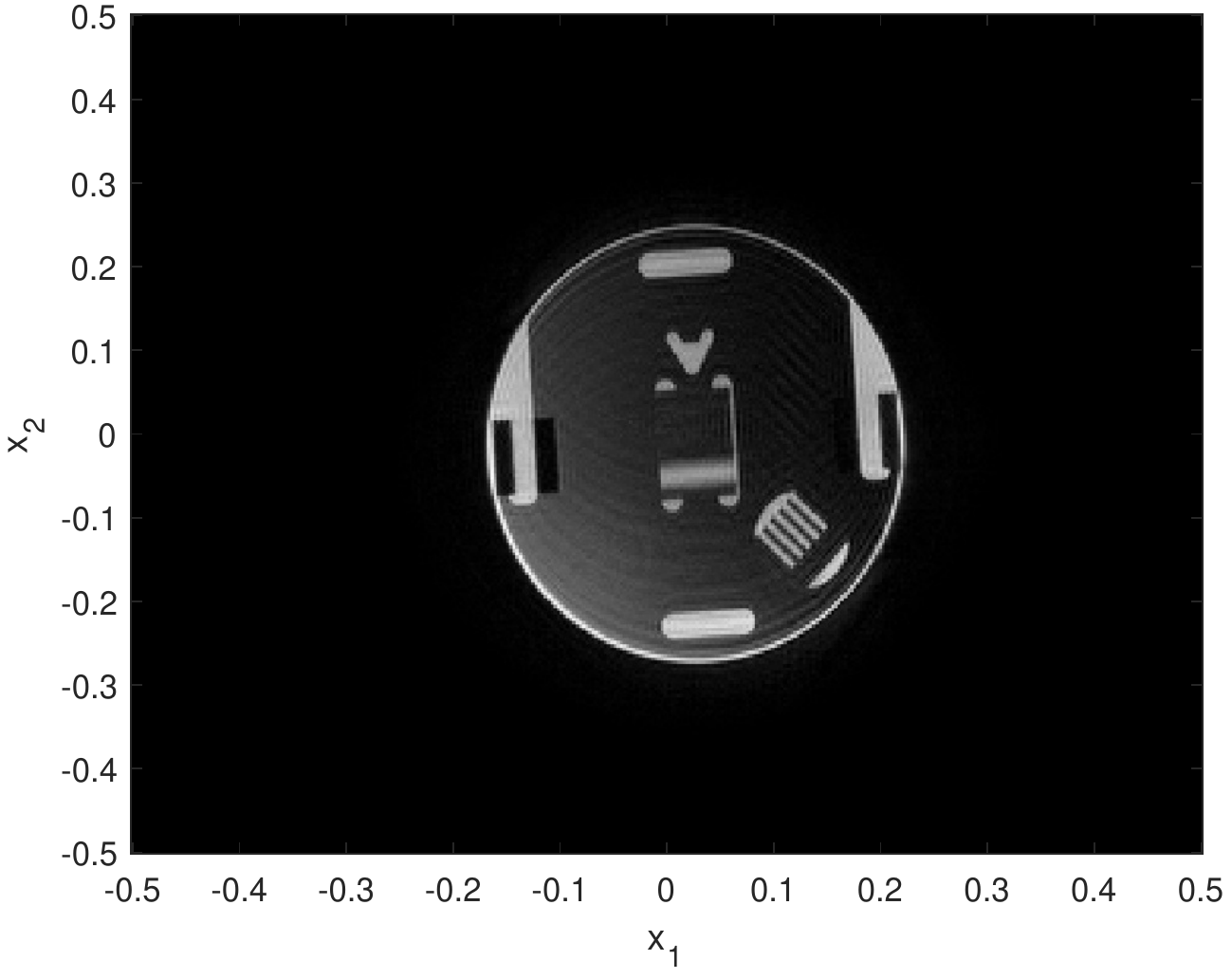}
\end{subfigure}
\hspace{-0.5cm}
\begin{subfigure}{0.24\textwidth}
\includegraphics[width=0.9\linewidth, height=3.2cm, keepaspectratio]{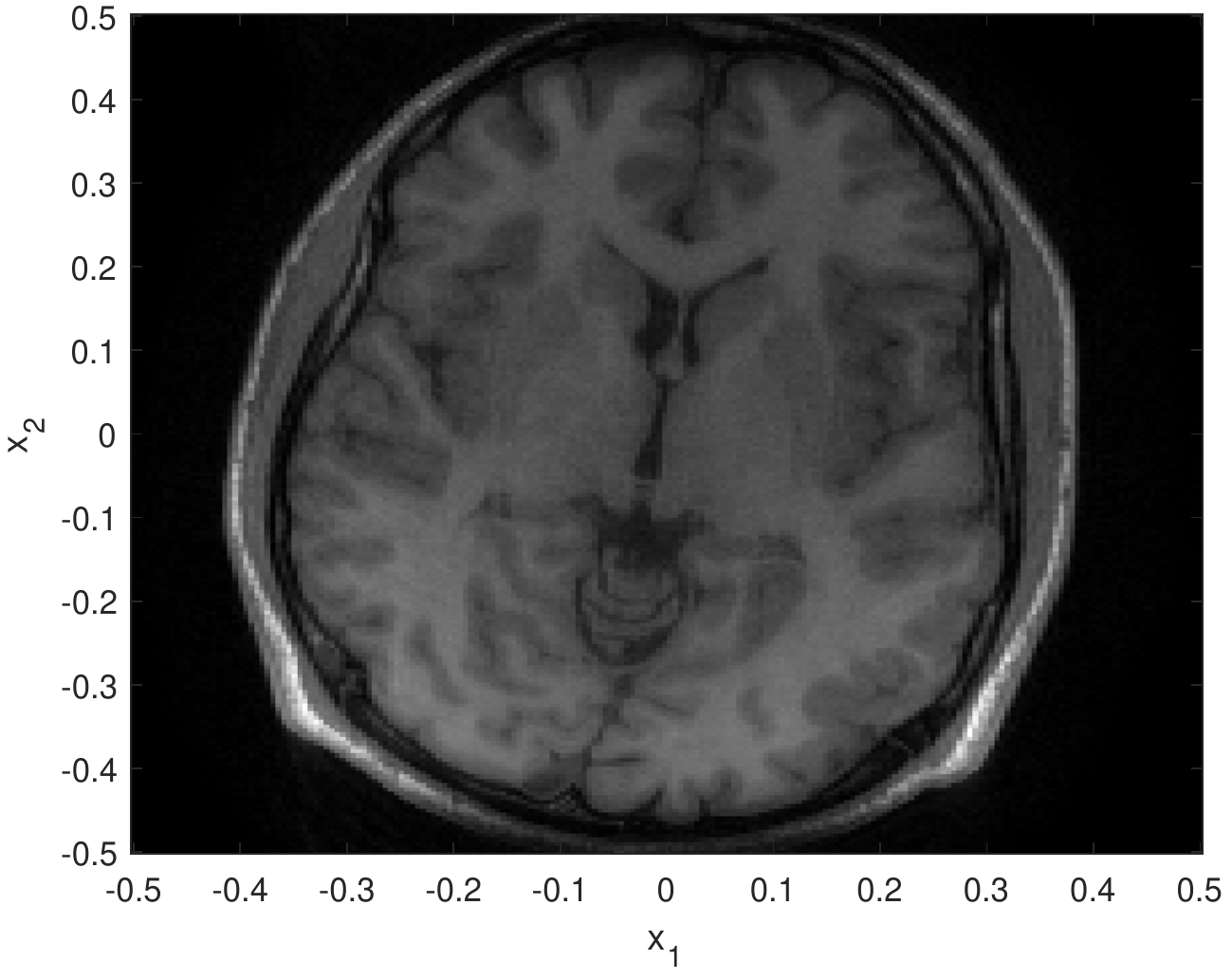}
\end{subfigure}
\hspace{-0.5cm}
\begin{subfigure}{0.24\textwidth}
\includegraphics[width=0.9\linewidth, height=3.2cm, keepaspectratio]{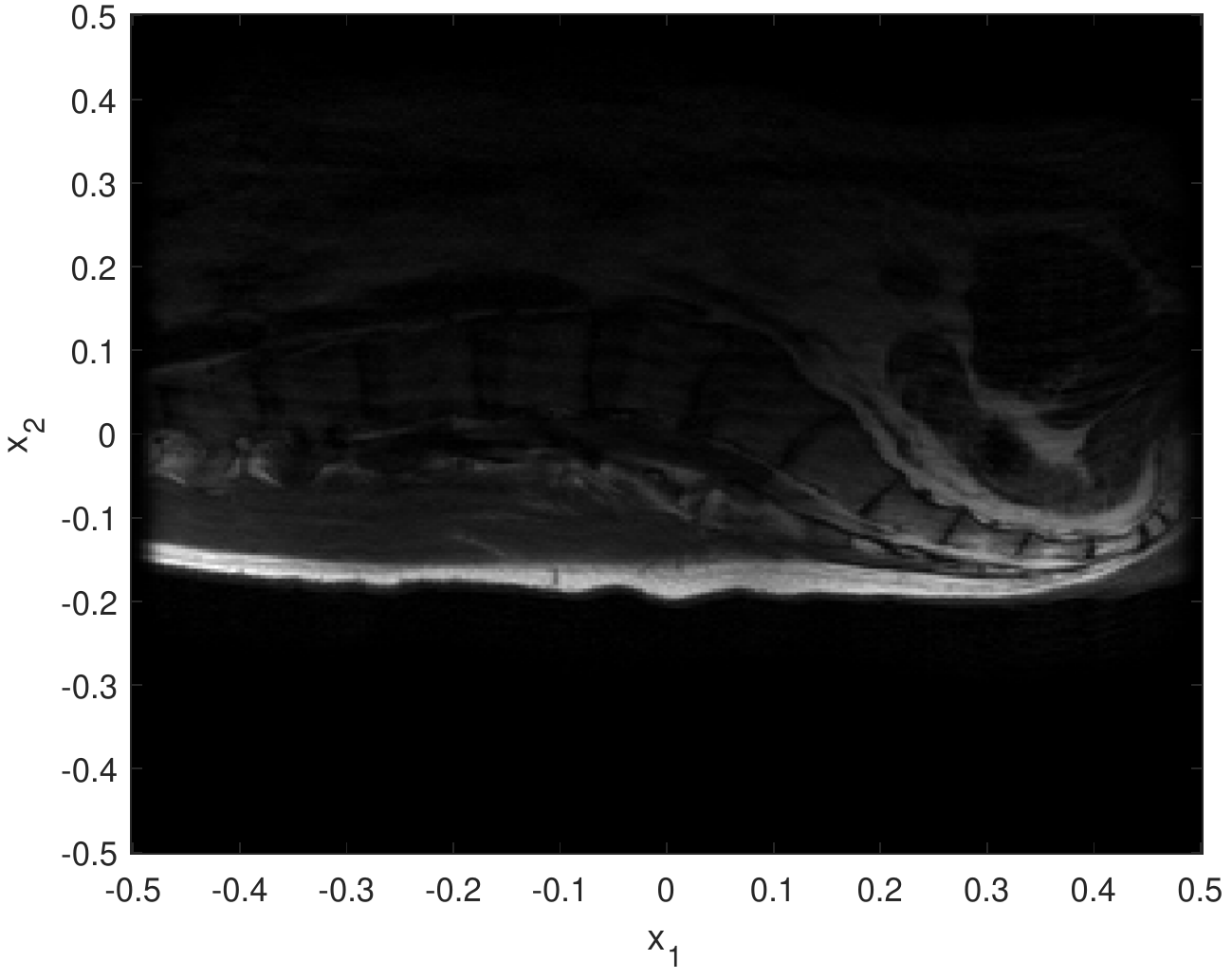}
\end{subfigure}
\hspace{-0.5cm}
\begin{subfigure}{0.24\textwidth}
\includegraphics[width=0.9\linewidth, height=3.2cm, keepaspectratio]{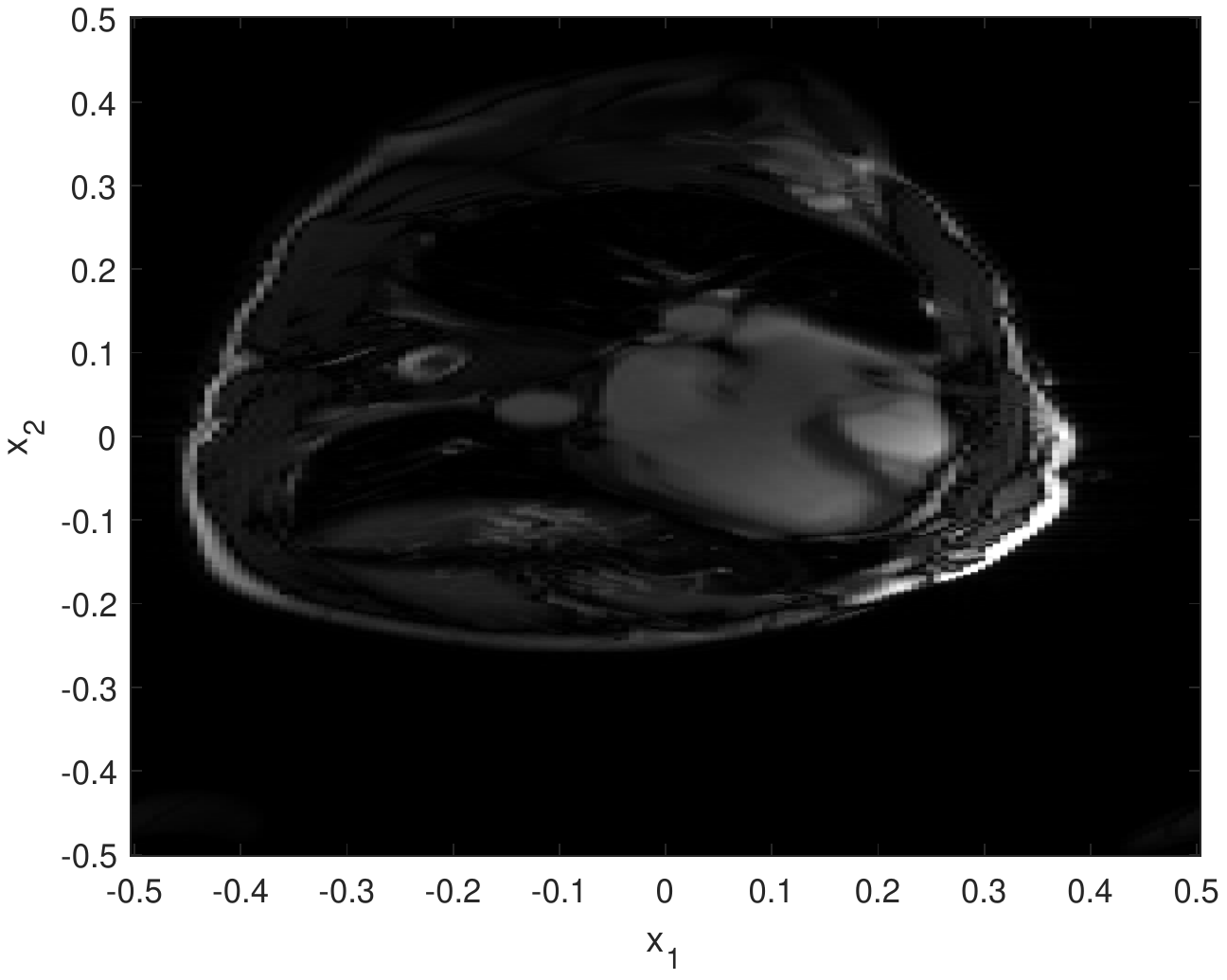}
\end{subfigure}
\begin{subfigure}{0.24\textwidth}
\includegraphics[width=0.9\linewidth, height=3.2cm, keepaspectratio]{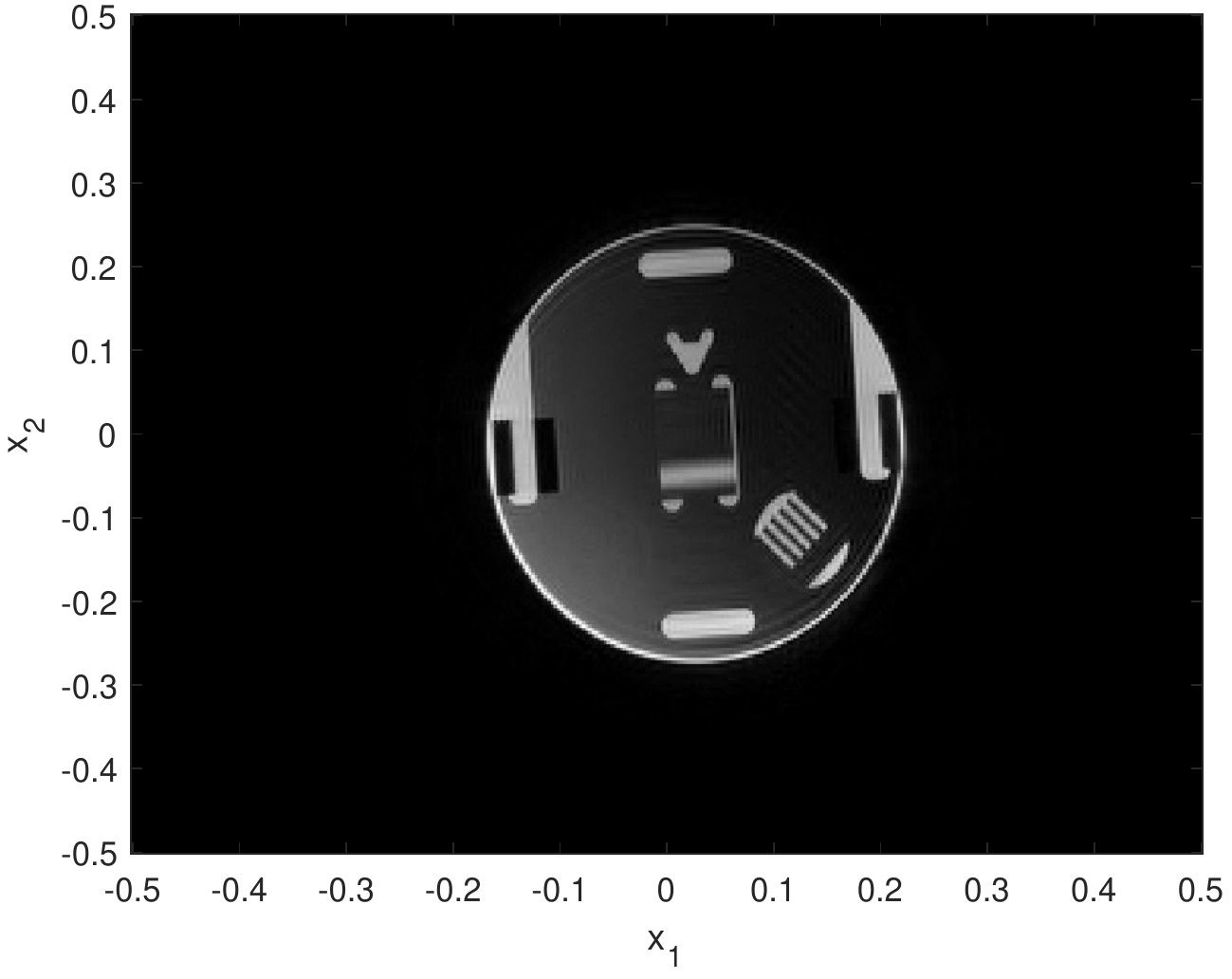}
\subcaption*{phantom ($50\%$)}
\end{subfigure}
\hspace{-0.5cm}
\begin{subfigure}{0.24\textwidth}
\includegraphics[width=0.9\linewidth, height=3.2cm, keepaspectratio]{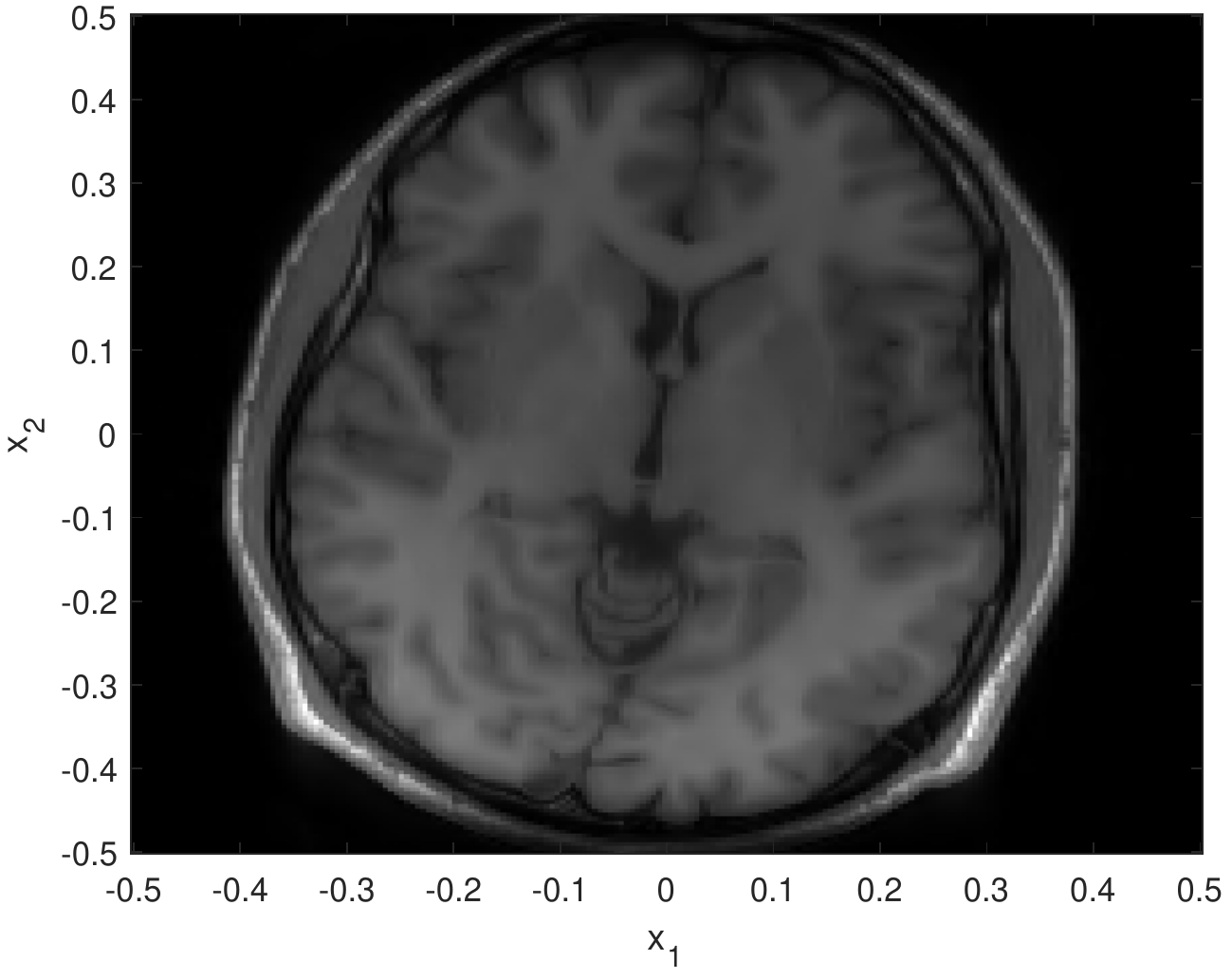}
\subcaption*{brain ($60\%$)}
\end{subfigure}
\hspace{-0.5cm}
\begin{subfigure}{0.24\textwidth}
\includegraphics[width=0.9\linewidth, height=3.2cm, keepaspectratio]{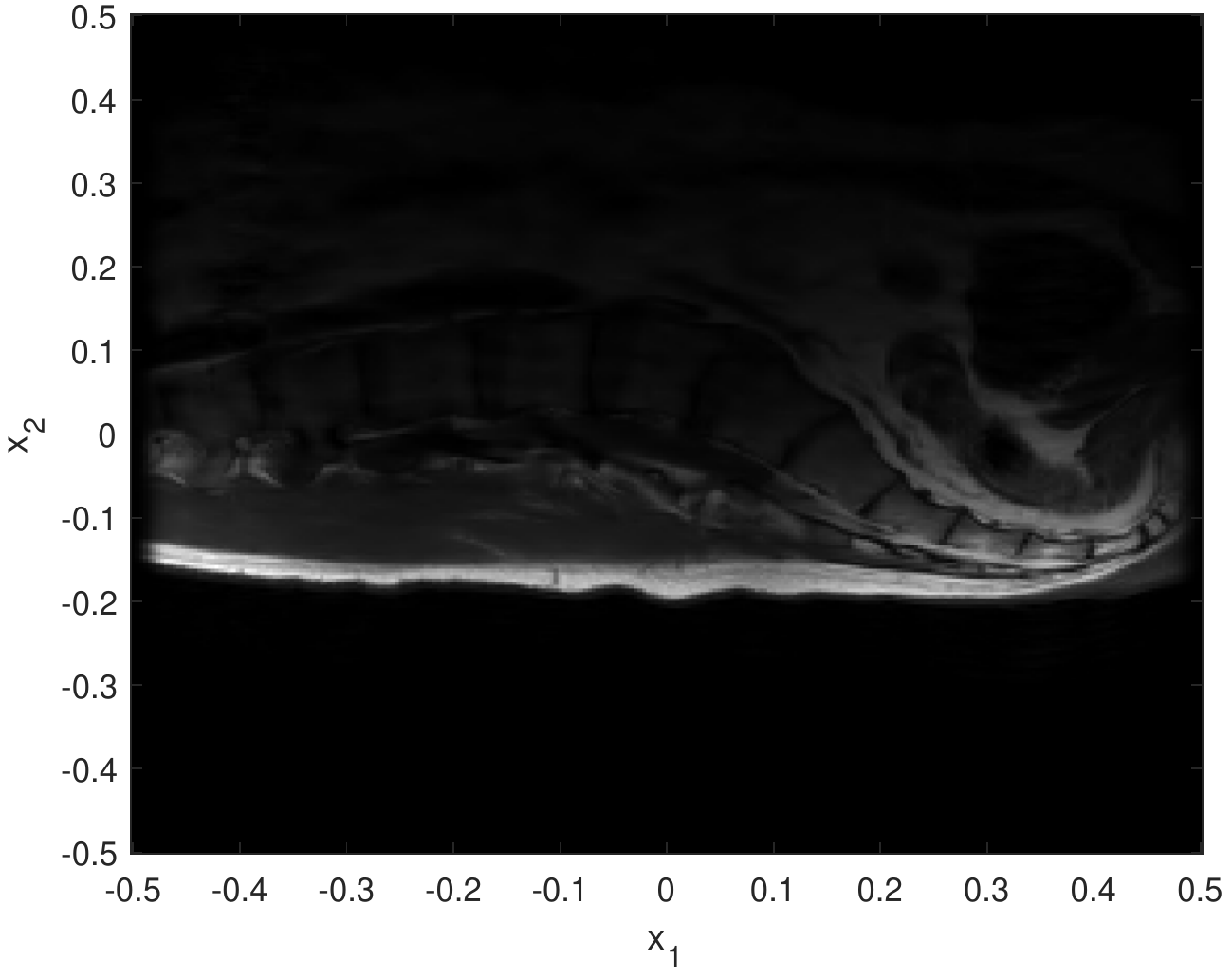}
\subcaption*{spine ($60\%$)}
\end{subfigure}
\hspace{-0.5cm}
\begin{subfigure}{0.24\textwidth}
\includegraphics[width=0.9\linewidth, height=3.2cm, keepaspectratio]{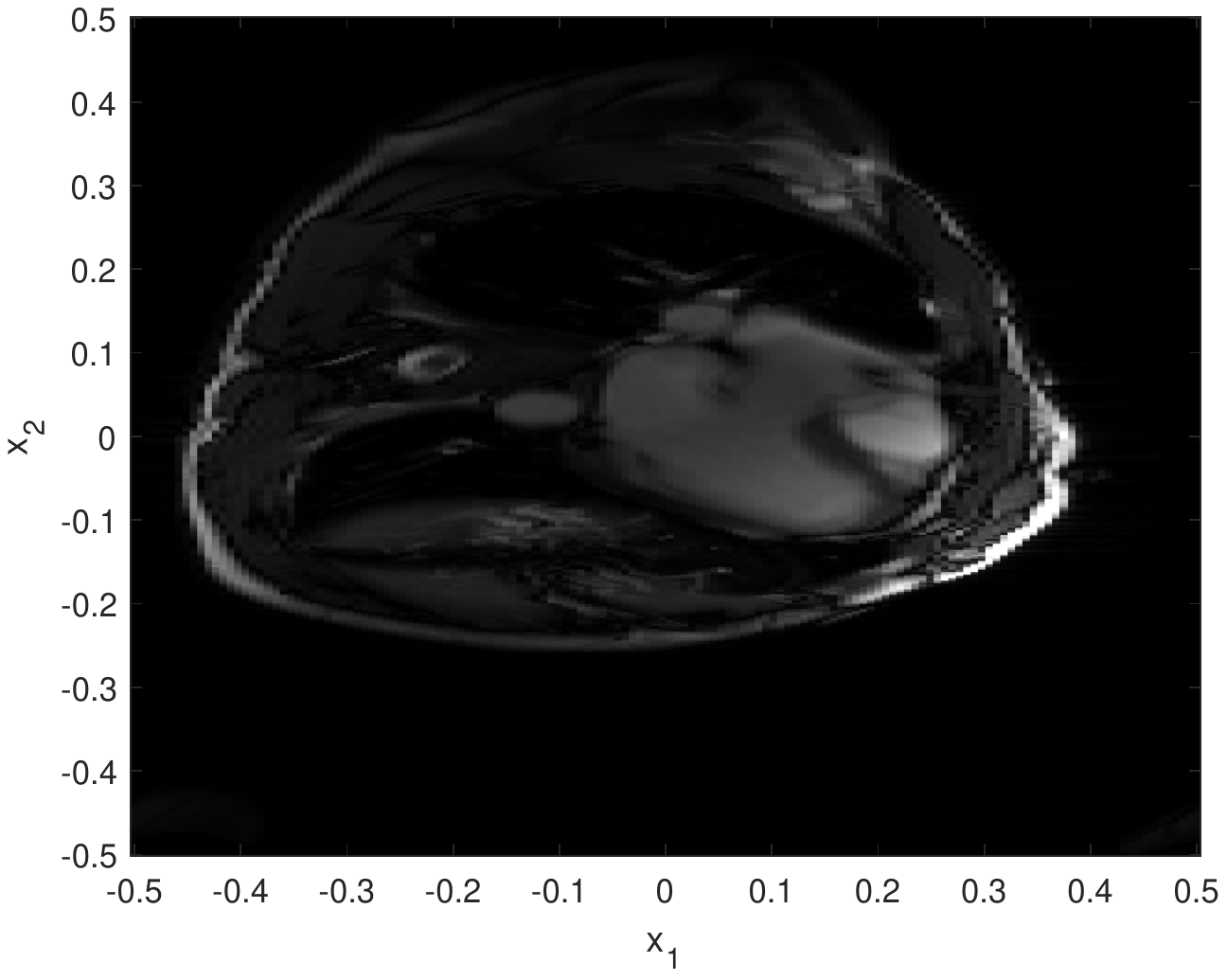}
\subcaption*{cardiac ($70\%$)}
\end{subfigure}
\caption{Example image reconstructions. Row 1 - SOS images. Row 2 - AC reconstructions. Row 3 - TV ESP reconstructions. Row 4 - ESP reconstructions. Row 5 - $L^1$ ESP reconstructions. The scan times in each case are given in parenthesis on the bottom row.}
\label{phantom_recon}
\end{figure}

In regards to the real phantom image, while the structural similarity offered by AC is optimal, the methods of the literature slightly outperform AC in terms of $\epsilon$ when the scan time is $>50\%$ (see the top left of figure \ref{error curves}). For the brain image, AC offers highly competitive performance in terms of $\epsilon$ and $\text{SSIM}_{\mu}$ score across all scan times. In the case of the cardiac image, the methods of the literature outperform AC in terms of $\epsilon$ when the scan time is greater than $40\%$. All methods are competitive in terms of $\epsilon$ in the cardiac example when $\text{scan time}\leq 40\%$. To summarize the results of figure \ref{error curves}, AC demonstrates optimal overall performance in terms structural similarity across a range of scan times with random undersampling, when compared to similar methods ESP, TV ESP, and $L^1$ ESP. In particular, AC is more robust in cases when the data is most limited and the scan time is $<50\%$, e.g., for the spine and real phantom images. The least squares error offered by AC is largely competitive with ESP, TV ESP, and $L^1$ ESP (e.g., for the phantom, brain, and spine images), and in only limited examples (i.e., the cardiac example with $\text{scan time}>40\%$) did we see a more significant increase in error using AC when compared to ESP, TV ESP, and $L^1$ ESP. 
\begin{figure}[!h]
\centering
\begin{subfigure}{0.24\textwidth}
\includegraphics[width=0.9\linewidth, height=3.2cm, keepaspectratio]{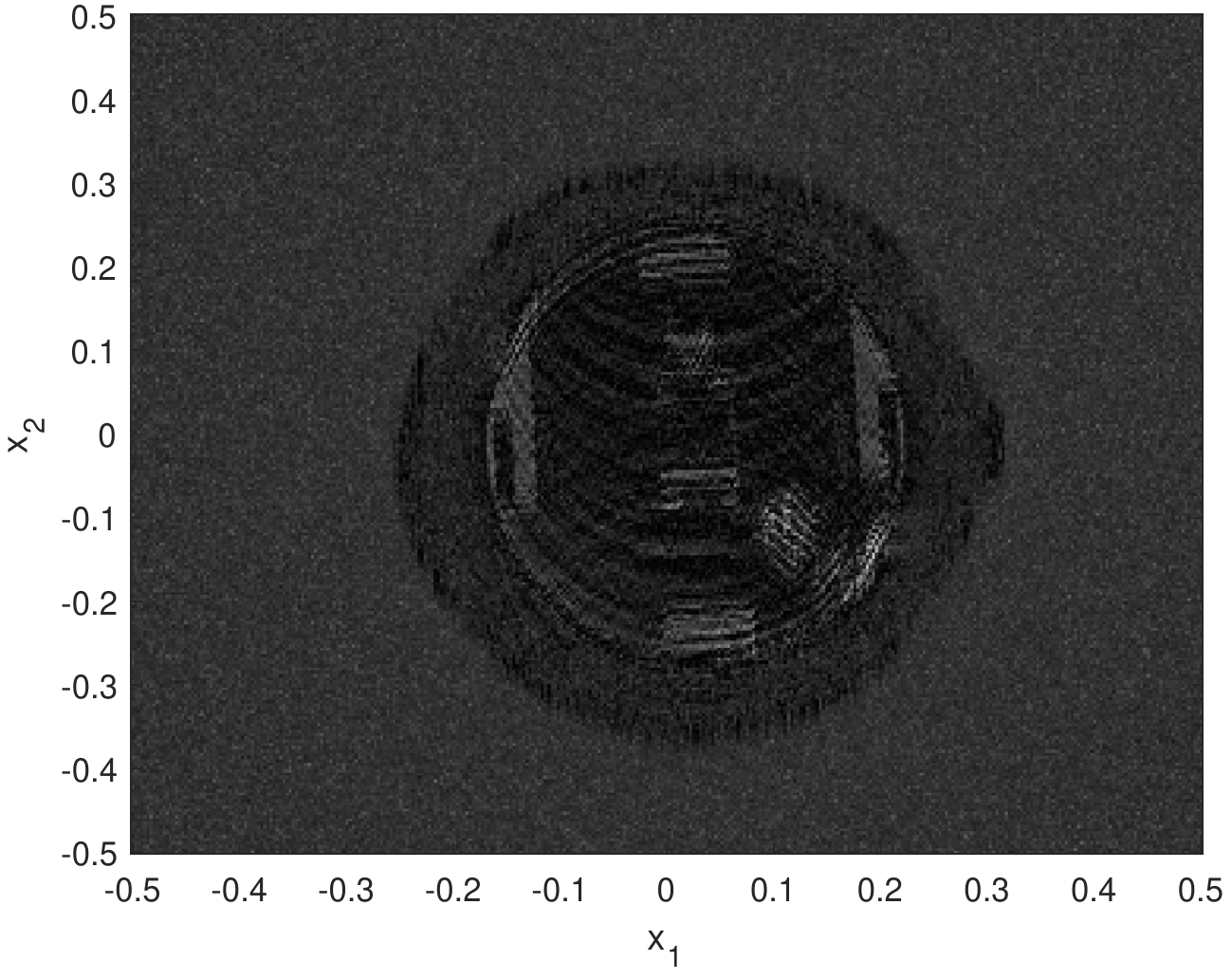}
\end{subfigure}
\hspace{-0.5cm}
\begin{subfigure}{0.24\textwidth}
\includegraphics[width=0.9\linewidth, height=3.2cm, keepaspectratio]{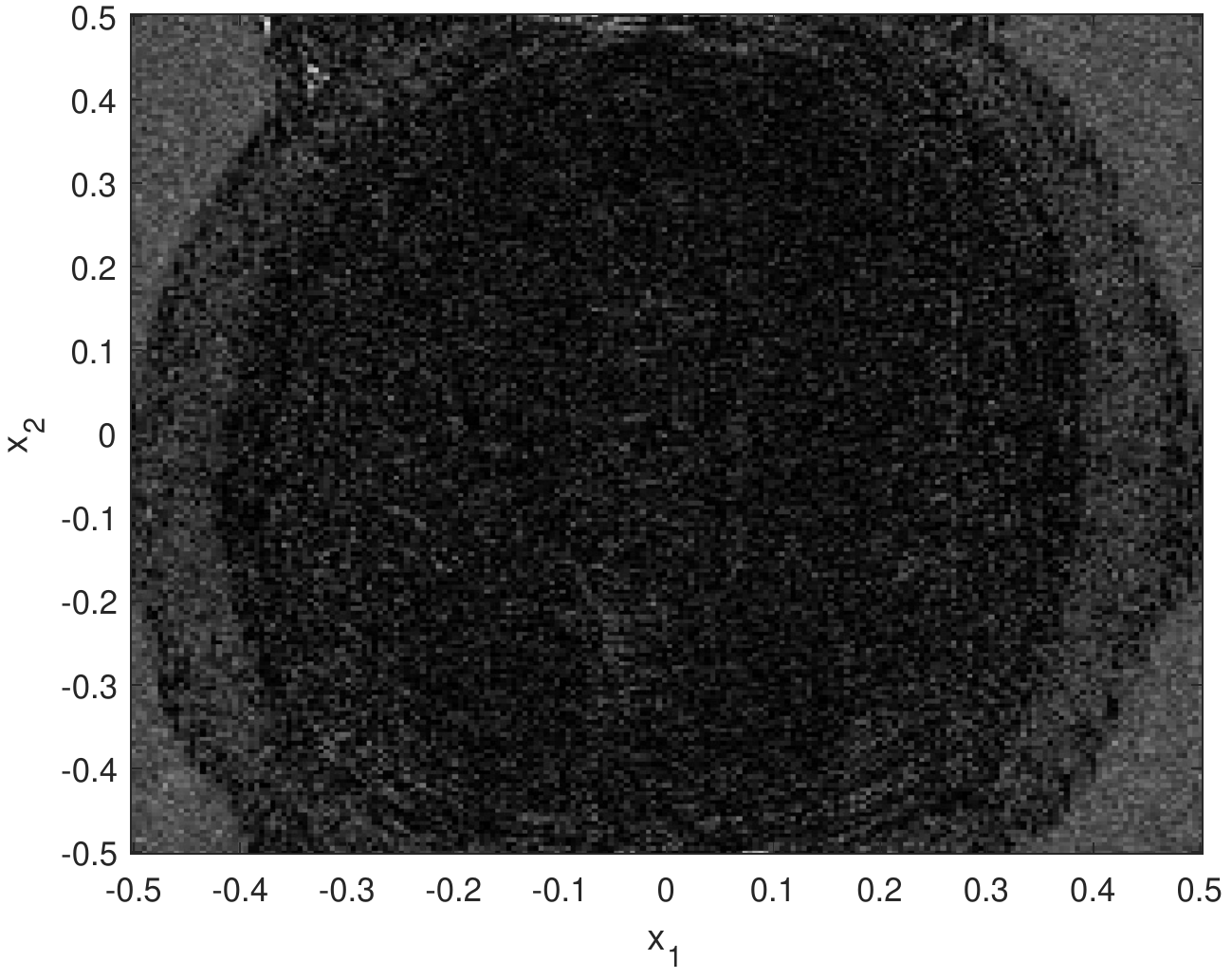}
\end{subfigure}
\hspace{-0.5cm}
\begin{subfigure}{0.24\textwidth}
\includegraphics[width=0.9\linewidth, height=3.2cm, keepaspectratio]{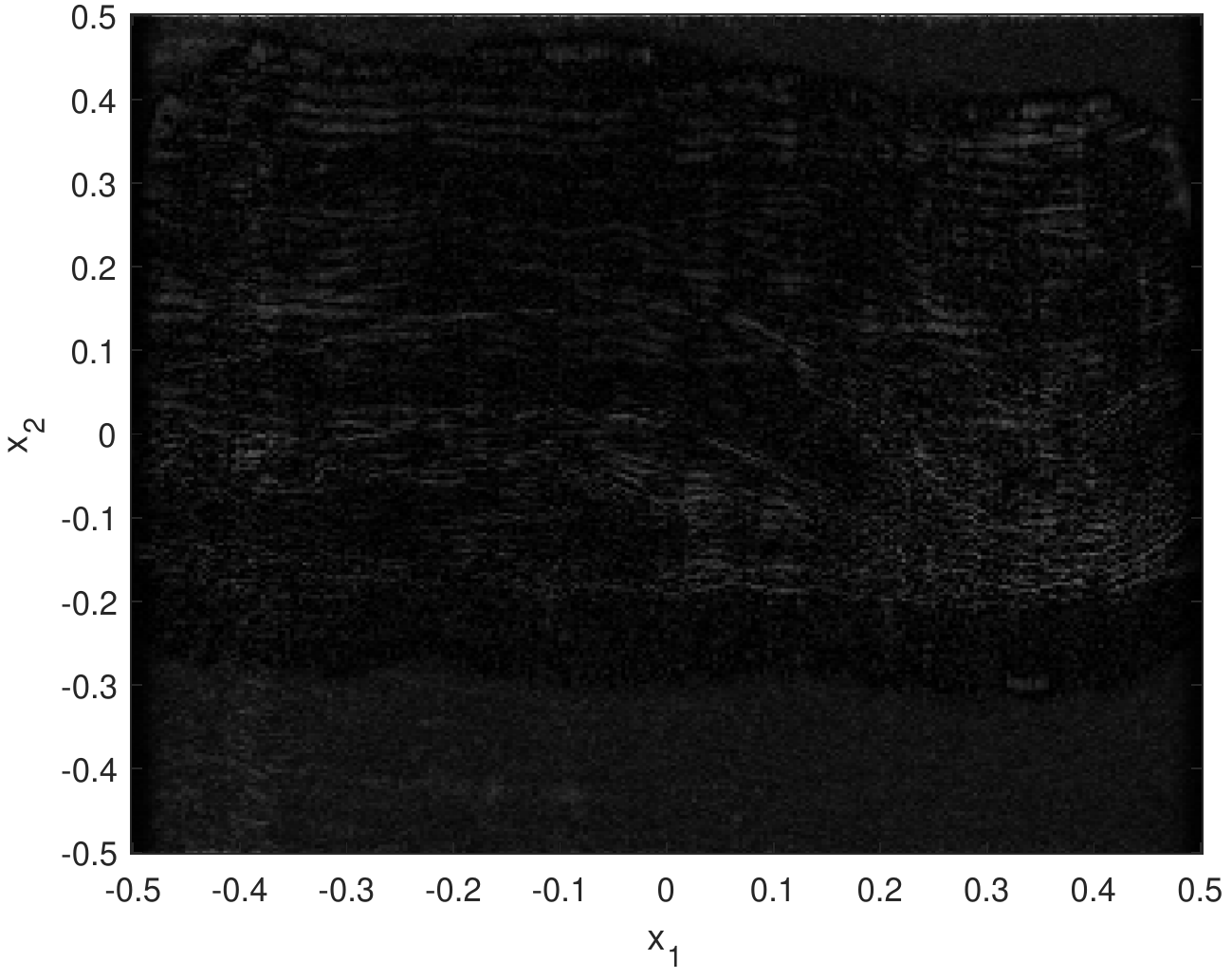}
\end{subfigure}
\hspace{-0.5cm}
\begin{subfigure}{0.24\textwidth}
\includegraphics[width=0.9\linewidth, height=3.2cm, keepaspectratio]{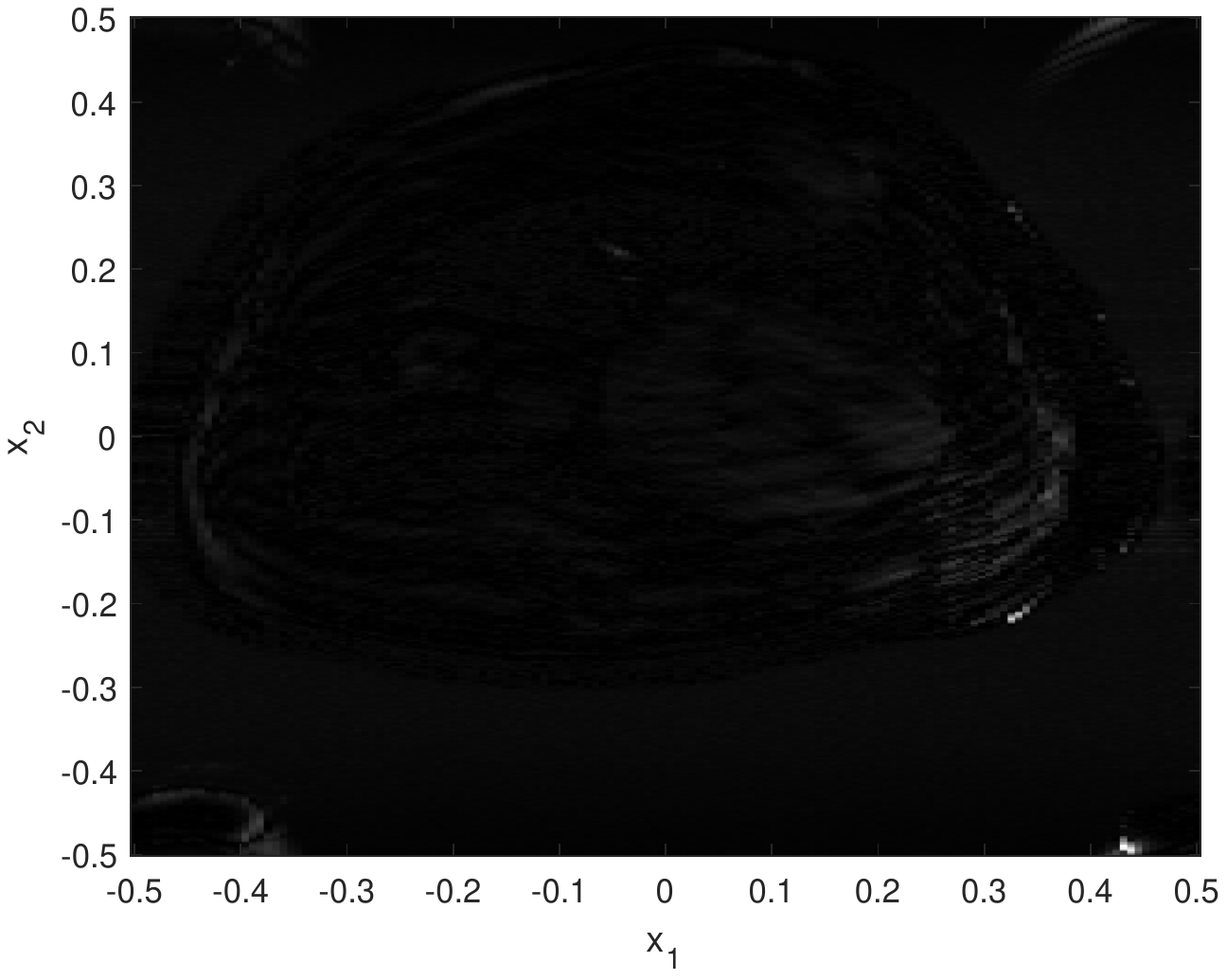}
\end{subfigure}
\begin{subfigure}{0.24\textwidth}
\includegraphics[width=0.9\linewidth, height=3.2cm, keepaspectratio]{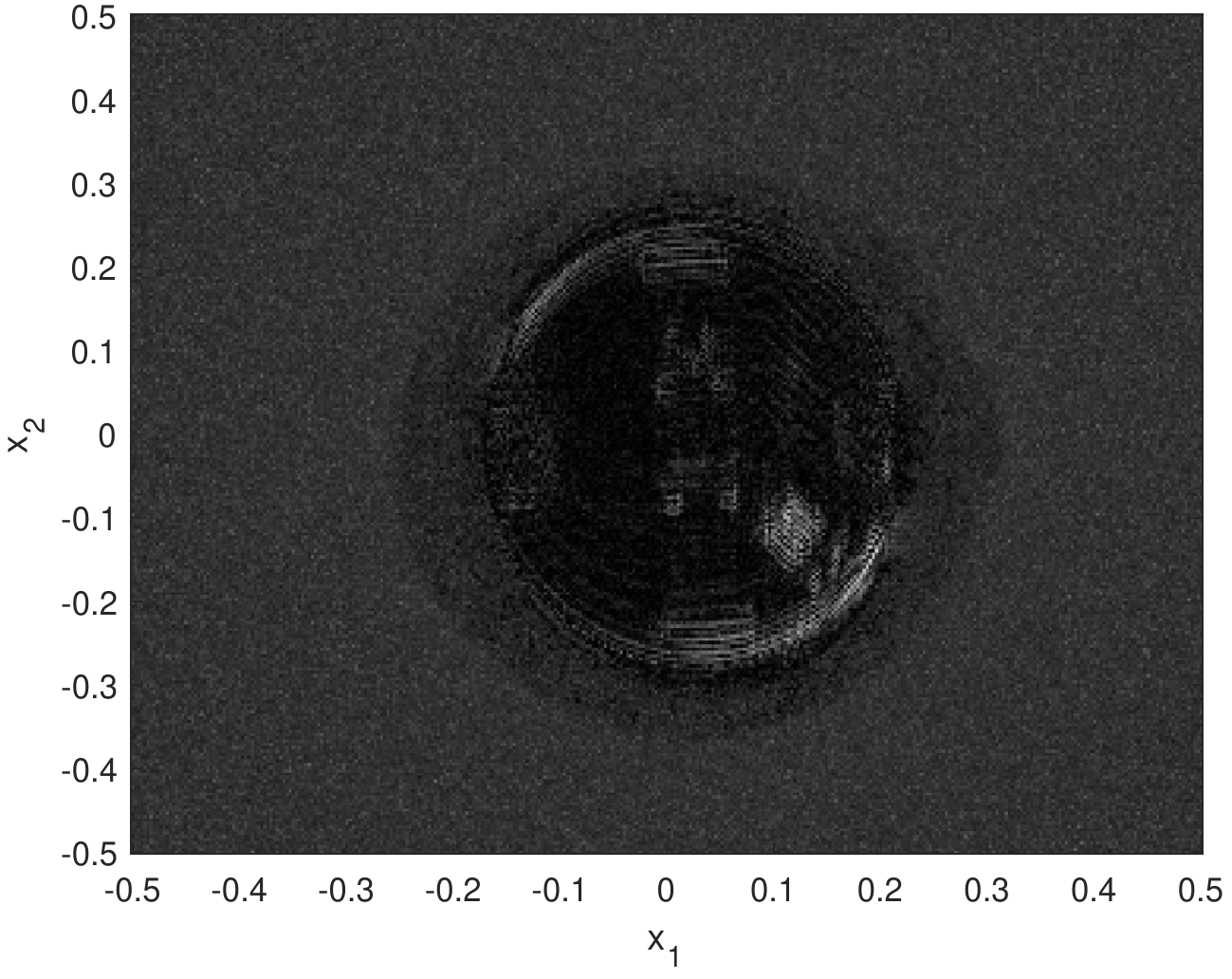}
\end{subfigure}
\hspace{-0.5cm}
\begin{subfigure}{0.24\textwidth}
\includegraphics[width=0.9\linewidth, height=3.2cm, keepaspectratio]{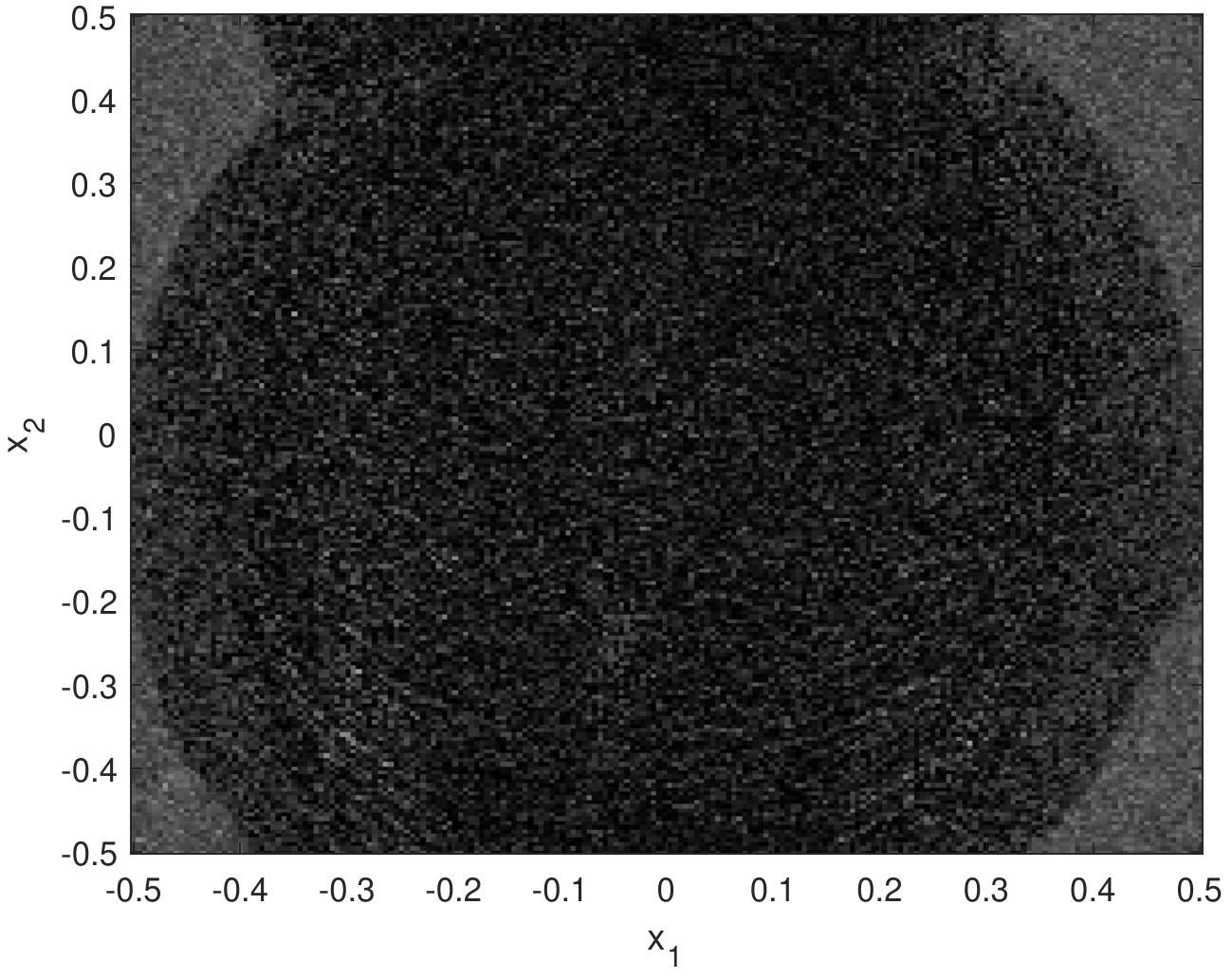}
\end{subfigure}
\hspace{-0.5cm}
\begin{subfigure}{0.24\textwidth}
\includegraphics[width=0.9\linewidth, height=3.2cm, keepaspectratio]{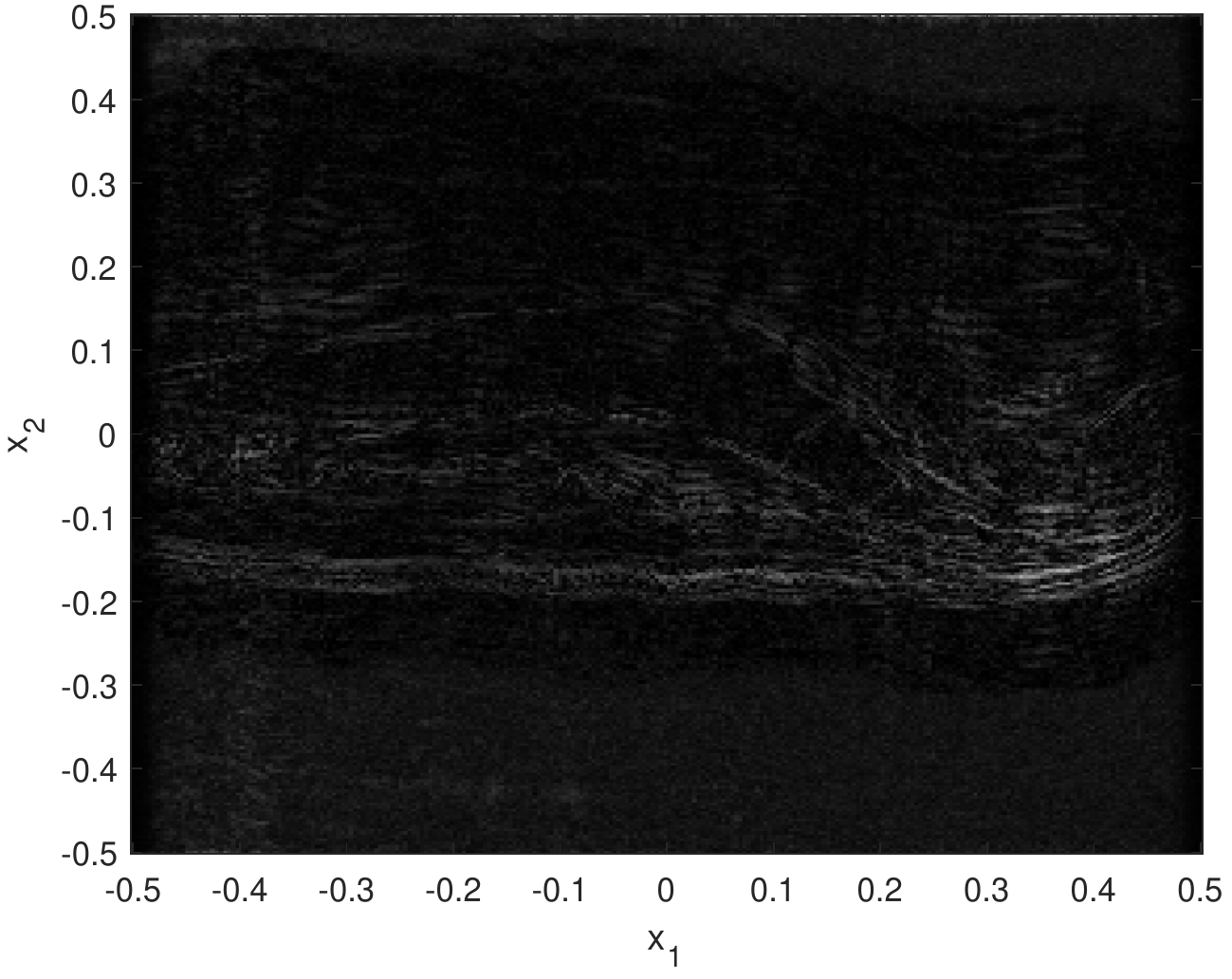}
\end{subfigure}
\hspace{-0.5cm}
\begin{subfigure}{0.24\textwidth}
\includegraphics[width=0.9\linewidth, height=3.2cm, keepaspectratio]{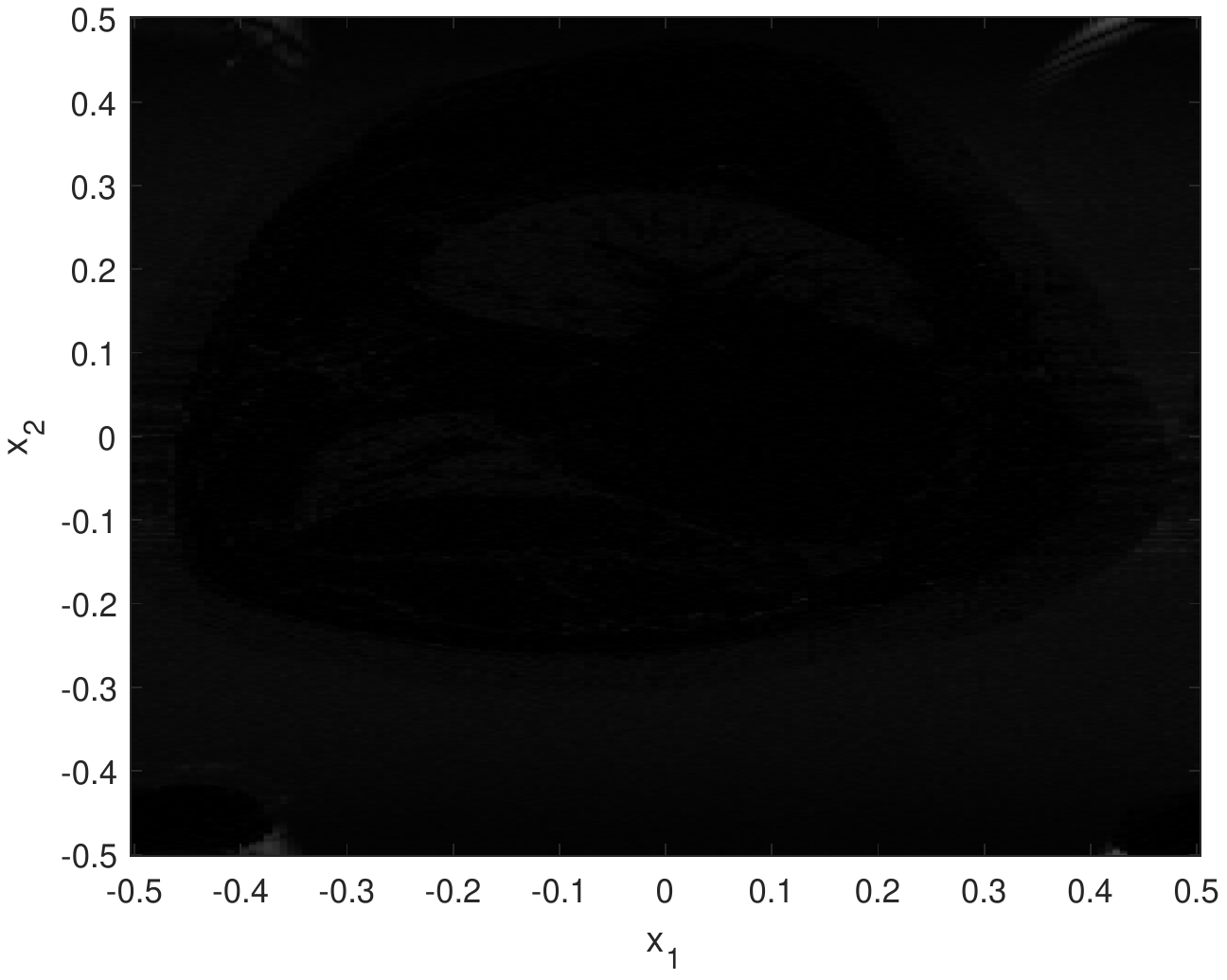}
\end{subfigure}
\begin{subfigure}{0.24\textwidth}
\includegraphics[width=0.9\linewidth, height=3.2cm, keepaspectratio]{random_real_phantom_6}
\end{subfigure}
\hspace{-0.5cm}
\begin{subfigure}{0.24\textwidth}
\includegraphics[width=0.9\linewidth, height=3.2cm, keepaspectratio]{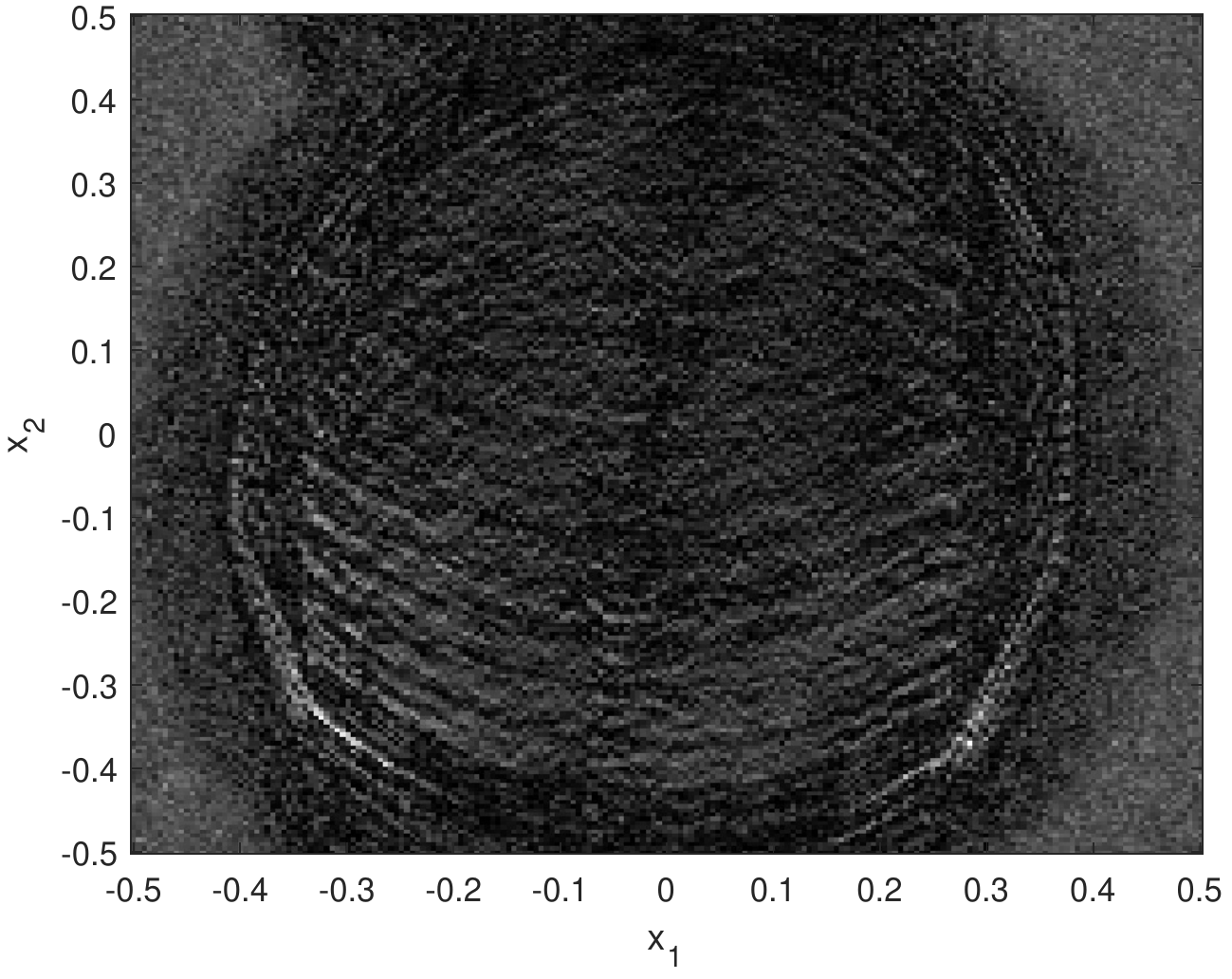}
\end{subfigure}
\hspace{-0.5cm}
\begin{subfigure}{0.24\textwidth}
\includegraphics[width=0.9\linewidth, height=3.2cm, keepaspectratio]{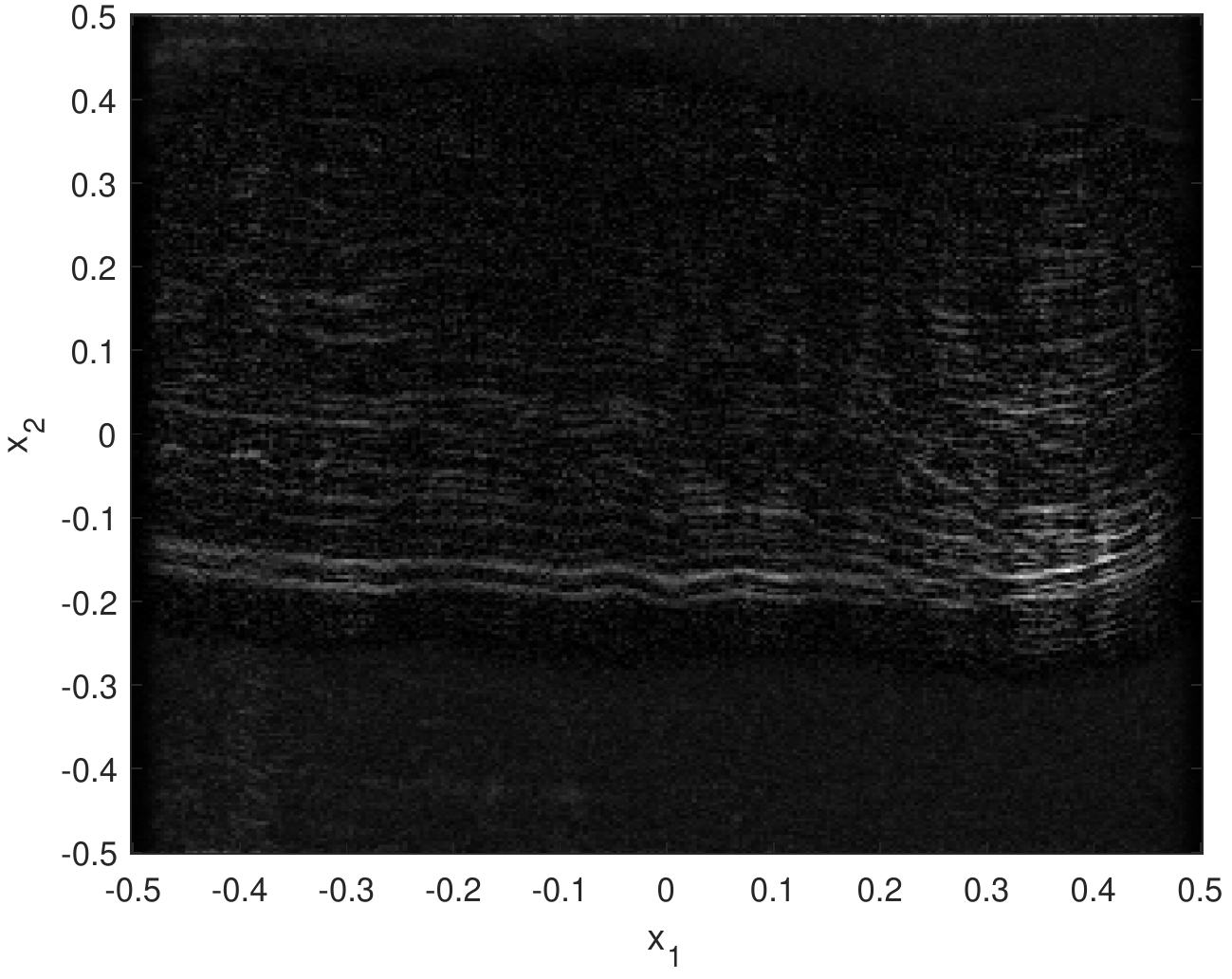}
\end{subfigure}
\hspace{-0.5cm}
\begin{subfigure}{0.24\textwidth}
\includegraphics[width=0.9\linewidth, height=3.2cm, keepaspectratio]{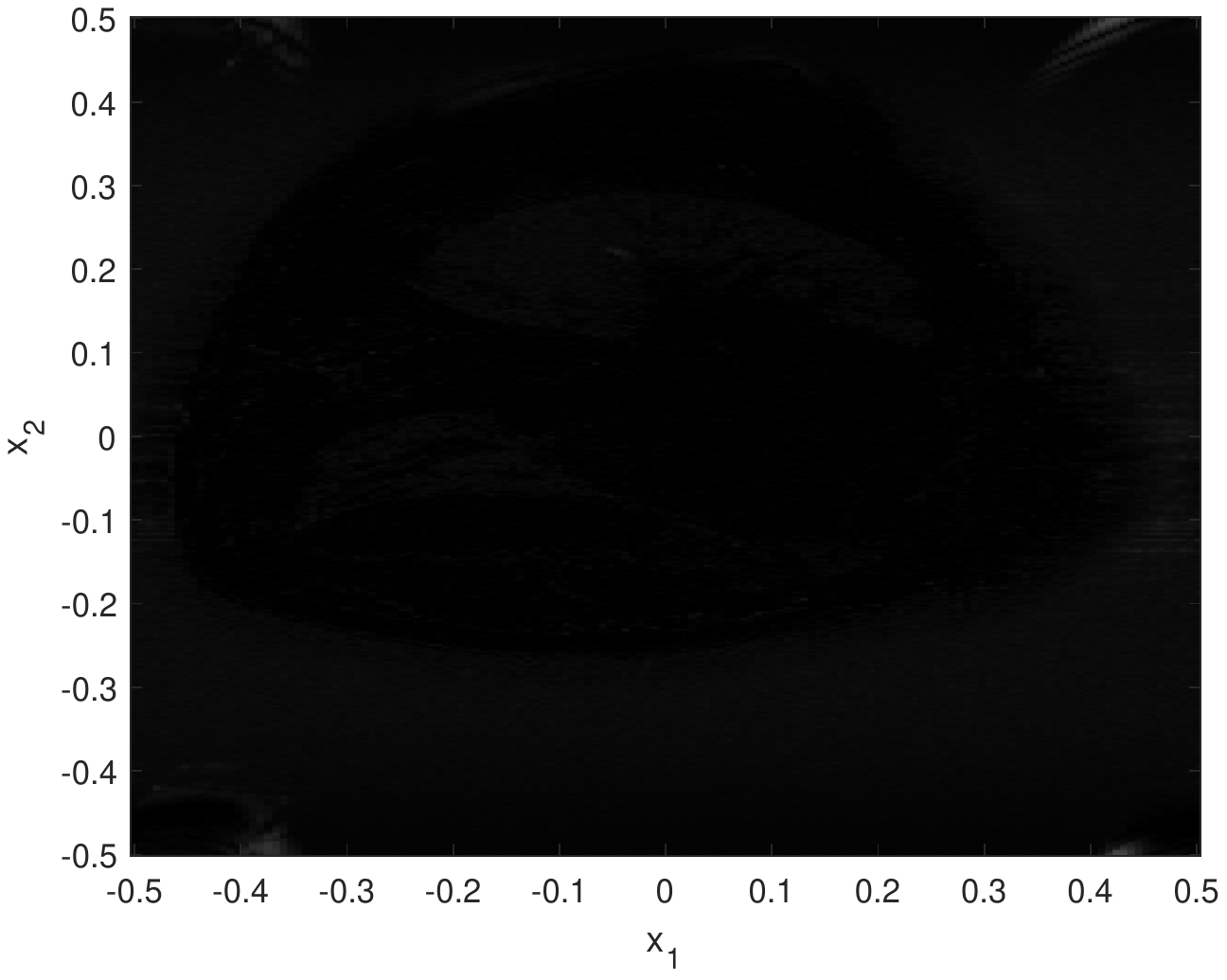}
\end{subfigure}
\begin{subfigure}{0.24\textwidth}
\includegraphics[width=0.9\linewidth, height=3.2cm, keepaspectratio]{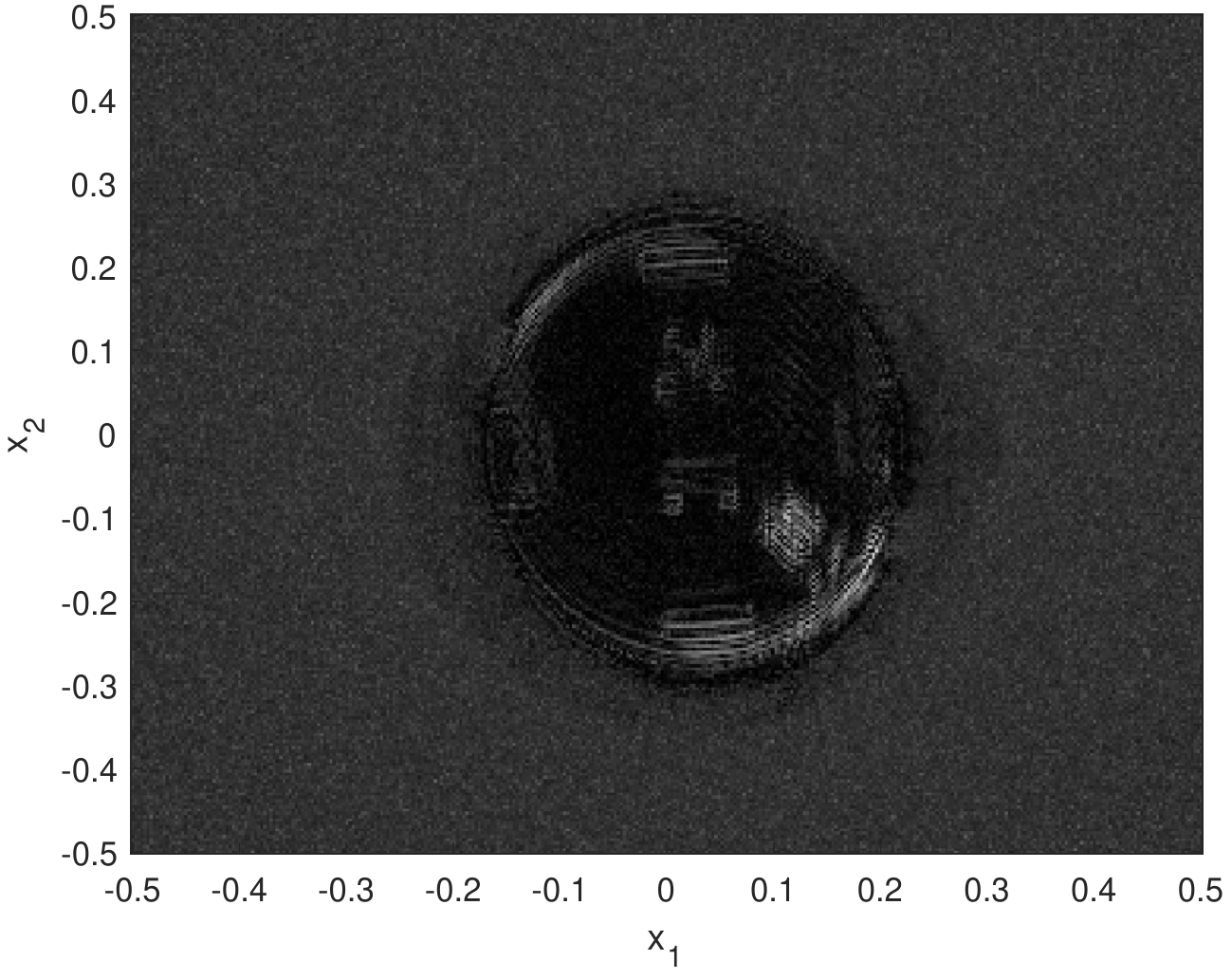}
\subcaption*{phantom ($50\%$)}
\end{subfigure}
\hspace{-0.5cm}
\begin{subfigure}{0.24\textwidth}
\includegraphics[width=0.9\linewidth, height=3.2cm, keepaspectratio]{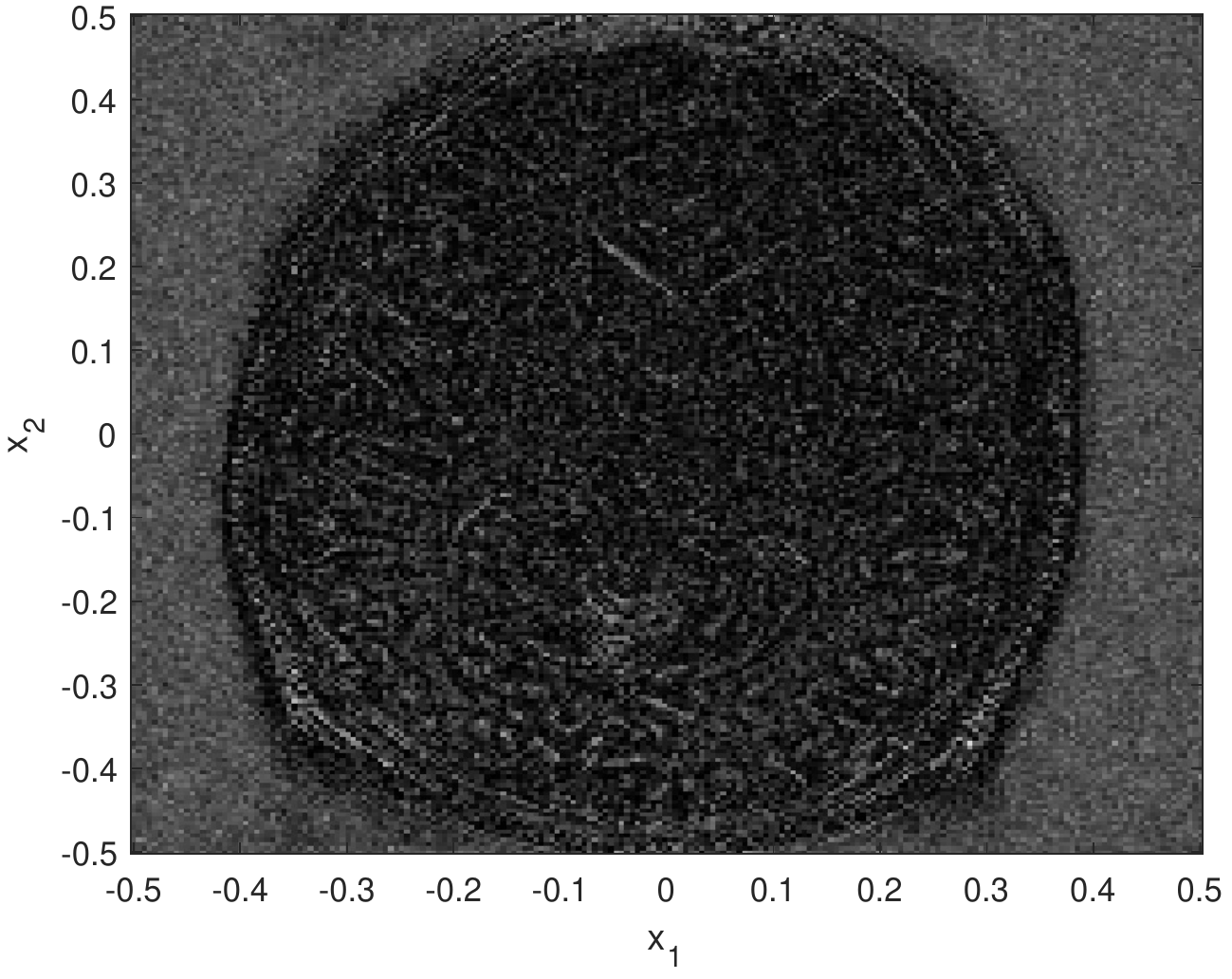}
\subcaption*{brain ($60\%$)}
\end{subfigure}
\hspace{-0.5cm}
\begin{subfigure}{0.24\textwidth}
\includegraphics[width=0.9\linewidth, height=3.2cm, keepaspectratio]{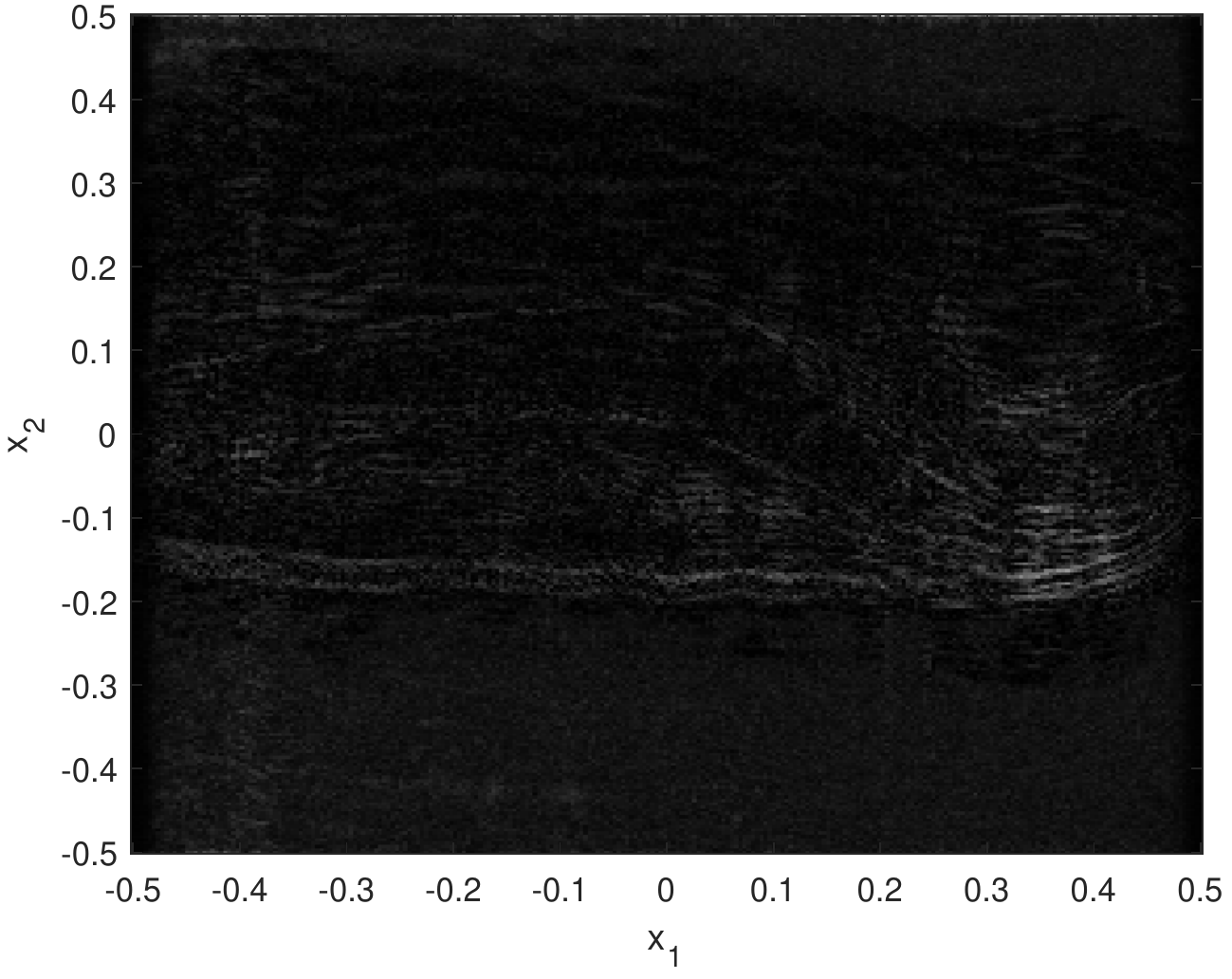}
\subcaption*{spine ($60\%$)}
\end{subfigure}
\hspace{-0.5cm}
\begin{subfigure}{0.24\textwidth}
\includegraphics[width=0.9\linewidth, height=3.2cm, keepaspectratio]{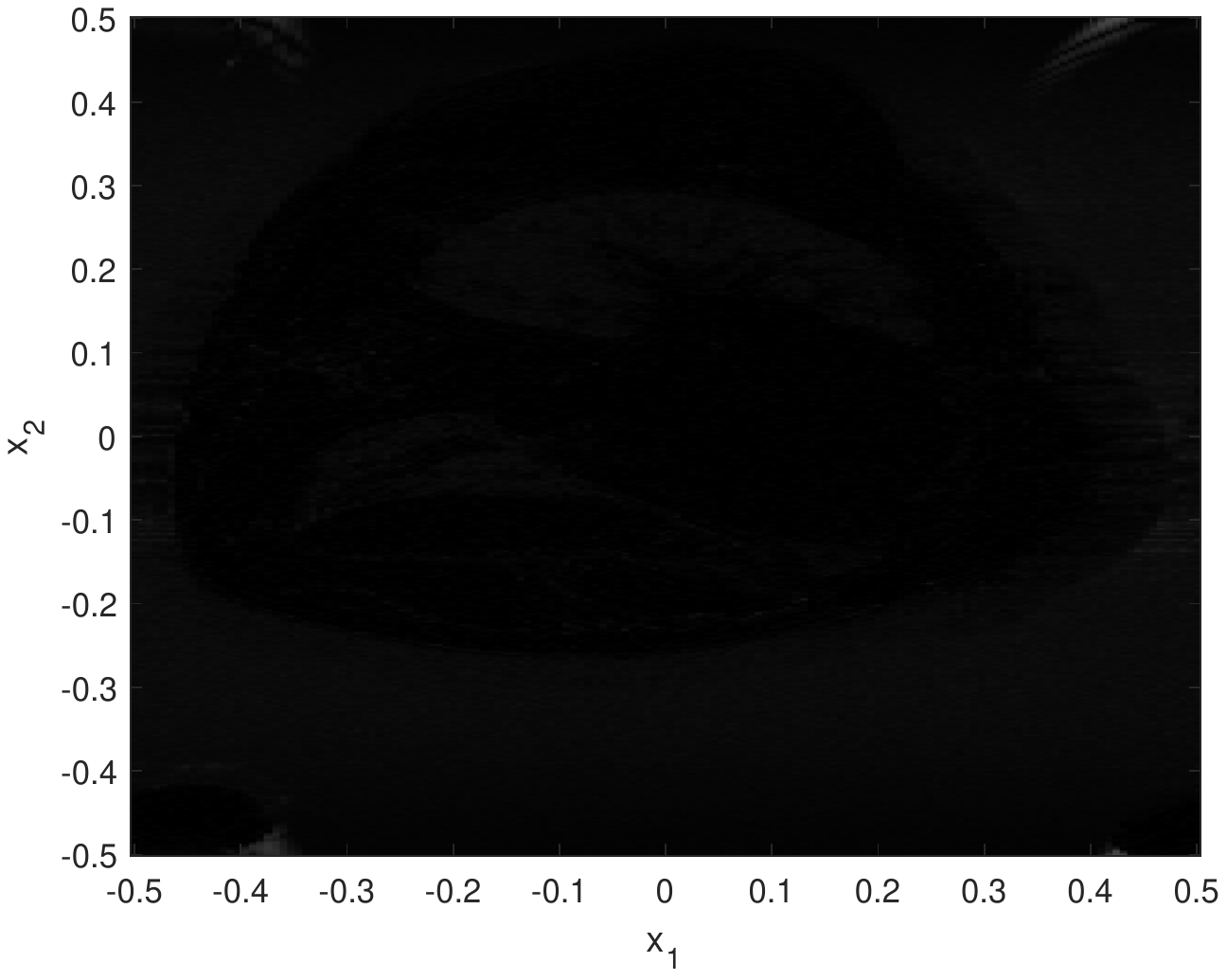}
\subcaption*{cardiac ($70\%$)}
\end{subfigure}
\caption{Absolute error images corresponding to the results of figure \ref{phantom_recon}. Row 1 - AC errors. Row 2 - TV ESP errors. Row 3 - ESP errors. Row 4 - $L^1$ ESP errors. The scan times in each case are given in parenthesis on the bottom row.}
\label{error_recon}
\end{figure}

To facilitate this conversation further, we show some of the image reconstructions which correspond to specific points on the error curves of figure \ref{error curves}, and give a side-by-side reconstruction quality comparison of all methods. We also present the absolute error images, which are calculated as $|F-F_{\epsilon}|$, where $F$ is the SOS image, and $F_{\epsilon}$ is the reconstruction. See figures \ref{phantom_recon} and \ref{error_recon}. Note, all image reconstructions compared throughout this paper are on the same grayscale (colorbar). The reconstruction examples are chosen to reflect scenarios where AC performs optimally, and sub-optimally in terms of $\epsilon$ and $\text{SSIM}_{\mu}$, when compared to the methods of the literature. 

In the first column of figure \ref{phantom_recon}, we present real phantom image reconstructions at $50\%$ scan time. AC, TV ESP, and $L^1$ ESP offer competitive image quality, and the sharp edges (jump discontinuities) of the image are recovered with high resolution. In the first column of figure \ref{error_recon}, we show the corresponding absolute error images which highlight some of the more specific differences in performance between each method. For example, the AC errors appear more uniformly spread across the ROI and lower in magnitude when compared to ESP, TV ESP, and $L^1$ ESP, where the error is more concentrated and greater in magnitude, particularly towards the bottom-right and top-left corners of the phantom. Thus, the ESP, TV ESP and $L^1$ ESP reconstructions induce stronger artifacts in more localized regions (e.g., the bottom-right region), which may help to explain the difference in structural similarity scores when compared to AC. In this example, AC outperforms ESP, TV ESP and $L^1$ ESP in terms of structural similarity, whereas all methods perform comparably in terms of overall error $\epsilon$. Thus it appears the improvement in structural similarity using AC is likely due to more uniform spread of the error and greater suppression of high intensity artifacts.

In the second column of figure \ref{phantom_recon}, we show image reconstructions of the brain image at $60\%$ scan time. In this example, AC and TV ESP offer the highest image quality and are competitive in terms of performance. On close inspection of the absolute error images in the second column of figure \ref{error_recon}, the errors in the AC and TV ESP reconstructions are similar, except AC contains more noise outside the brain, and TV ESP slightly greater noise within the brain. In the ESP reconstruction, there are artifacts, leading to significantly increased error and reduced structural similarity when compared to the other methods. In the $L^1$ ESP reconstruction, the artifacts appear more concentrated near the image edges, although this does not cause any significant reduction in structural similarity when compared to AC and TV ESP. The $\epsilon$ scores are also similar.

In the third column of figure \ref{phantom_recon}, we present reconstructions of the spine cross-section at $60\%$ scan time. In this example, AC offers the best image quality and performance in terms of $\epsilon$ and $\text{SSIM}_{\mu}$, which is further evidenced by the absolute error images in the third column of figure \ref{error_recon}. The overall presence of artifacts is reduced in the AC reconstruction, when compared to ESP, TV ESP and $L^1$ ESP, particularly along the high density tissue behind the spine (the bottom outline of the image) and towards the top of the spine near the neck. As in the previous example, ESP does not perform well, and there are significant artifacts in the reconstruction. TV ESP and $L^1$ ESP offer similar image quality.

For our final example, in the fourth column of figure \ref{phantom_recon}, we present cardiac image reconstructions at $70\%$ scan time. This is an example where AC under-performs when compared to ESP, TV ESP and $L^1$ ESP. All methods offer high image quality with $\epsilon<3\%$, and $\text{SSIM}_{\mu}>98\%$. There are no observable artifacts, on this scale, in the  ESP, TV ESP and $L^1$ ESP absolute error images in figure \ref{error_recon}. There are some mild artifacts in the AC reconstruction which are most focused towards the bottom-right corner of the image. The artifacts slightly deform the image contrast in the AC reconstruction and cause an increase in $\epsilon$, but do not effect the structure of the reconstruction. This is evidenced by the structural similarity scores, which are comparable across all methods in this example.

\section{Discussion and Conclusions}
We have introduced a new reconstruction methodology (denoted AC), based on the analytic continuation ideas of \cite{natterer}, for limited-data, multi coil MRI, whereby the regions of missing $\textbf{k}$-space are lines parallel to either $k_1$ or $k_2$. 
Specifically, we showed that the limited-data MRI problem could be reduced to a set of one-dimensional Fredholm integral equations. In section \ref{SVD}, we presented an SVD analysis of the Fredholm operators, which gave insight to the problem stability. For example, we investigated how the number of coils effected problem stability in section \ref{coils}, and showed that, as the number of coils increased, the condition number of the operator decreased, which suggests a more stable inversion. 

In section \ref{results}, we compared AC against three similar methods from the literature, namely ESP, TV ESP, and $L^1$ ESP. AC was shown to offer optimal performance in terms of structural similarity across a range of image examples, with varying resolution, numbers of coils (e.g., $K=4,8,16,34$), application (e.g., brain, spine and cardiac imaging), and amounts of missing data. In particular, AC was more robust to data limitations (e.g., less coils, lower scan time) when compared to  ESP, TV ESP, and $L^1$ ESP.

In section \ref{hyperparam}, we discussed the selection of hyperparameters for ESP, TV ESP, $L^1$ ESP, and AC. The hyperparameter selection was done heuristically, and required some trial and error. In future work, we aim to investigate the effectiveness of AC when the hyperparameters are chosen using hyperparameter selection methods, such as the 
discrepancy principle \cite{DP1}.

In further work, we aim to extend the SVD analysis ideas of section \ref{SVD}, e.g., in conjunction with $g$-factor analysis \cite{SENSE,G1}, to determine optimal sub-sampling schemes for different imaging applications (e.g., {\it{in vivo}} medical MRI). 

\section*{Acknowledgments}
As per our agreement for use of the real phantom data of \url{https://mr.usc.edu/download/data/}, we acknowledge NSF support, specifically NSF grant CCF-1350563. We would like to thank Prof. Andre Vanderkouwe and Prof. Robert Frost for their engaging and helpful discussion, and for their help in communicating this work to an MRI audience. The authors also wish to acknowledge funding support from the Massachusetts Life Sciences Center Bits to Bytes Program, and Abcam, Inc.

\end{document}